\documentclass[10pt]{article}
\usepackage{mathrsfs}
\usepackage{amsfonts,mathrsfs,amssymb,amsthm,mathptm}
\usepackage{amsmath,amscd}
\usepackage{mathptm,pslatex}
\usepackage{color}
\usepackage[dvips]{graphicx}
\usepackage[all]{xy}
\usepackage{graphicx}
\usepackage{fancyhdr}

\begin{document}
\newtheorem{Def}{Definition}[section]
\newtheorem{Bsp}[Def]{Example}
\newtheorem{Prop}[Def]{Proposition}
\newtheorem{Theo}[Def]{Theorem}
\newtheorem{Lem}[Def]{Lemma}
\newtheorem{Koro}[Def]{Corollary}
\theoremstyle{definition}
\newtheorem{Rem}[Def]{Remark}

\newcommand{\add}{{\rm add}}
\newcommand{\Coker}{{\rm Coker}}
\newcommand{\Cone }{{\rm Cone}}
\newcommand{\gd}{{\rm gl.dim}}
\newcommand{\dm}{{\rm dom.dim}}
\newcommand{\DTr}{{\rm DTr}}
\newcommand{\E}{{\rm E}}
\newcommand{\Mor}{{\rm Morph}}
\newcommand{\End}{{\rm End}}
\newcommand{\Ext}{{\rm Ext}}
\newcommand{\Fac}{{\rm Fac}}
\newcommand{\fd} {{\rm fin.dim}}
\newcommand{\fld}{{\rm flat.dim}}
\newcommand{\Gen}{{\rm Gen}}
\newcommand{\Hom}{{\rm Hom}}
\newcommand{\id}{{\rm inj.dim}}
\newcommand{\ind}{{\rm ind}}
\newcommand{\Img  }{{\rm Im}}
\newcommand{\Ker  }{{\rm Ker}}
\newcommand{\ol}{\overline}
\newcommand{\overpr}{$\hfill\square$}
\newcommand{\pd}{{\rm proj.dim}}
\newcommand{\projdim}{\pd}
\newcommand{\rad}{{\rm rad}}
\newcommand{\rsd}{{\rm res.dim}}
\newcommand{\rd} {{\rm rep.dim}}
\newcommand{\soc}{{\rm soc}}
\renewcommand{\top}{{\rm top}}

\newcommand{\cpx}[1]{#1^{\bullet}}
\newcommand{\D}[1]{{\mathscr D}(#1)}
\newcommand{\Dz}[1]{{\mathscr D}^+(#1)}
\newcommand{\Df}[1]{{\mathscr D}^-(#1)}
\newcommand{\Db}[1]{{\mathscr D}^b(#1)}
\newcommand{\C}[1]{{\mathscr C}(#1)}
\newcommand{\Cz}[1]{{\mathscr C}^+(#1)}
\newcommand{\Cf}[1]{{\mathscr C}^-(#1)}
\newcommand{\Cb}[1]{{\mathscr C}^b(#1)}
\newcommand{\K}[1]{{\mathscr K}(#1)}
\newcommand{\Kz}[1]{{\mathscr K}^+(#1)}
\newcommand{\Kf}[1]{{\mathscr  K}^-(#1)}
\newcommand{\Kb}[1]{{\mathscr K}^b(#1)}
\newcommand{\modcat}{\ensuremath{\mbox{{\rm -mod}}}}
\newcommand{\Modcat}{\ensuremath{\mbox{{\rm -Mod}}}}
\newcommand{\gm}{{\rm _{\Gamma_M}}}
\newcommand{\gmr}{{\rm _{\Gamma_M^R}}}
\newcommand{\stmodcat}[1]{#1\mbox{{\rm -{\underline{mod}}}}}
\newcommand{\StHom}{{\rm \underline{Hom}}}
\newcommand{\pmodcat}[1]{#1\mbox{{\rm -proj}}}
\newcommand{\imodcat}[1]{#1\mbox{{\rm -inj}}}
\newcommand{\Pmodcat}[1]{#1\mbox{{\rm -Proj}}}
\newcommand{\Imodcat}[1]{#1\mbox{{\rm -Inj}}}
\newcommand{\opp}{^{\rm op}}
\newcommand{\otimesL}{\otimes^{\rm\mathbb L}}
\newcommand{\rHom}{{\rm\mathbb R}{\rm Hom}\,}

\def\vez{\varepsilon}\def\bz{\bigoplus}  \def\sz {\oplus}
\def\epa{\xrightarrow} \def\inja{\hookrightarrow}

\newcommand{\lra}{\longrightarrow}
\newcommand{\lraf}[1]{\stackrel{#1}{\lra}}
\newcommand{\ra}{\rightarrow}
\newcommand{\dk}{{\rm dim_{_{k}}}}
\newcommand{\colim}{{\rm colim\, }}
\newcommand{\limt}{{\rm lim\, }}
\newcommand{\Add}{{\rm Add }}
\newcommand{\Prod}{{\rm Prod }}
\newcommand{\Tor}{{\rm Tor}}
\newcommand{\Cogen}{{\rm Cogen}}

{\Large \bf
\begin{center}

Ringel modules and homological subcategories
\end{center}}

\smallskip
\centerline{\textbf{Hongxing Chen} and
\textbf{Changchang Xi}$^*$}

\renewcommand{\thefootnote}{\alph{footnote}}
\setcounter{footnote}{-1} \footnote{ $^*$ Corresponding author.
Email: xicc@cnu.edu.cn; Fax: 0086 10 58808202.}
\renewcommand{\thefootnote}{\alph{footnote}}
\setcounter{footnote}{-1} \footnote{2010 Mathematics Subject
Classification: Primary 18E30, 16G10, 13B30; Secondary 16S10,
13E05.}
\renewcommand{\thefootnote}{\alph{footnote}}
\setcounter{footnote}{-1} \footnote{Keywords: Cotilting modules;
Derived category; Homological subcategory; Recollement; Ringel
modules; Tilting modules.}

\begin{abstract}
Given a good $n$-tilting module $T$ over a ring $A$, let $B$ be the
endomorphism ring of $T$, it is an open question whether the kernel
of the left-derived functor $T\otimesL_B-$ between the derived
module categories of $B$ and $A$ could be realized as the derived
module category of a ring $C$ via a ring epimorphism $B\ra C$ for
$n\ge 2$. In this paper, we first provide a uniform way to deal with
the above question both for tilting and cotilting modules by
considering a new class of modules called Ringel modules, and then
give criterions for the kernel of $T\otimesL_B-$ to be equivalent to
the derived module category of a ring $C$ with a ring epimorphism
$B\ra C$. Using these characterizations, we display both a positive
example of $n$-tilting modules from noncommutative algebra, and a
counterexample of $n$-tilting modules from commutative algebra to
show that, in general, the open question may have a negative answer.
As another application of our methods, we consider the dual question
for cotilting modules, and get corresponding criterions and
counterexamples. The case of cotilting modules, however, is much
more complicated than the case of tilting modules.
\end{abstract}

{\footnotesize\tableofcontents}

\section{Introduction \label{sect1}}

As is well known, tilting theory has had significant applications in
many branches of mathematics (see \cite{HHK}), and the key
objectives in this theory are tilting modules, or more generally,
tilting complexes or objects. Given a good tilting module $T$ over a
ring $A$, let $B$ be the endomorphism ring of $T$, if $T$ is
classical, then a beautiful theorem of Happel says that the derived
module category $\D B$ of $B$ is triangle equivalent to the derived
module category $\D A$ of $A$ (see \cite{H}). Thus one can use
derived invariants to understand homological, geometric and
numerical properties of $A$ through $B$, or conversely, of $B$
through $A$. This theorem also tells that one cannot get new derived
categories from classical tilting modules. For infinitely generated
tilting modules,  Bazzoni, Mantese and Tonolo recently show a
remarkable result: $\D A$ can be regarded as a full subcategory or a
quotient category of $\D B$ (see \cite{Bz2}). Moreover, it is proved
in \cite{CX1} that if the projective dimension of $T$ is at most
$1$, then there is a homological ring epimorphism $\lambda: B\ra C$
of rings such that the kernel of the total left-derived functor
$T\otimes_B^{\mathbb L}-$, as a full triangulated subcategory of
$\D{B}$, can be realized as the derived module category $\D{C}$ of
$C$. Thus, for (infinitely generated) good tilting modules of
projective dimension at most $1$, Happel's theorem now has a new
appearance and can be featured as a recollement of derived module
categories:
$$\xymatrix@C=1.2cm{\D{C}\ar[r]^-{D(\lambda_*)}
&\D{B}\ar[r]^-{{} _AT\otimesL_B-}\ar@/^1.2pc/[l]\ar@/_1.3pc/[l]
&\D{A} \ar@/^1.2pc/[l]\ar@/_1.3pc/[l]}\vspace{0.2cm}$$ However, for
tilting modules of higher projective dimension, the existence of the
above recollement is unknown (see the first open question in
\cite{CX1}). On the one hand, the argument used in \cite{CX1}
actually does not work any more for the general case because the
proof there involves a two-term complex which depends on the
projective dimension. Thus some new ideas are necessary for
attacking the general situation. On the other hand, neither positive
examples nor counterexamples to this general case are known to
experts. So, it is quite mysterious whether the above recollement
still exists for a good tilting module of projective dimension at
least $2$.

In the present paper, we shall consider this question in detail. In
fact, our discussion is implemented in the framework of Ringel
modules (see Definition \ref{rm}). This provides us a way to deal
with the above question uniformly for higher tilting and cotilting
modules. We first provide characterizations of when the kernel of
the functor $T\otimes_B^{\mathbb L}-$ can be realized as the derived
module category of a ring $C$ with a homological ring epimorphism
$B\ra C$, and then use these criterions to give positive and
negative examples to the above question for tilting modules of
projective dimension bigger than $1$. Finally, as another
application of our criterions, we shall consider the above question
for cotilting modules.

Before stating our main results precisely, we first introduce
notation and recall some definitions.

Let $A$ be a ring with identity, and let $n$ be a natural number. A
left $A$-module $T$ is called an $n$-\emph{tilting} $A$-module (see
\cite{ct}) if the following three conditions are satisfied:

$(T1)$ There is an exact sequence
$$
0\lra P_n \lra \cdots \lra P_1\lraf{\sigma} P_0\lraf{\pi} T\lra 0$$
of $A$-modules such that all $P_i$ are projective, that is, the
projective dimension of $T$ is at most $n$;

$(T2)$ $\Ext^j_A(T, T^{(I)})=0$ for all $j\geq 1$ and nonempty sets
$I$, where $T^{(I)}$ denotes the direct sum of $I$ copies of $T$;

$(T3)$ There is an exact sequence
$$
0\lra {}_AA \lraf{\omega} T_0 \lra T_1\lra \cdots \lra T_n\lra 0$$
of $A$-modules such that $T_i$ is isomorphic to a direct  summand of
a direct sum of copies of $T$ for all $0\leq i\leq n$.

An $n$-tilting module $T$ is said to be  \emph{good} if $(T3)$ can
be replaced by

$(T3)'$ there is an exact sequence
$$ 0\lra {}_AA \lraf{\omega} T_0 \lra
T_1\lra \cdots \lra T_n\lra 0$$ of $A$-modules such that $T_i$ is
isomorphic to a direct summand of a finite direct sum of copies of
$T$ for all $0\leq i\leq n$. A good $n$-tilting module $T$ is said
to be \emph{classical} if the modules $P_i$ in $(T1)$ are finitely
generated (see \cite{bb, hr}).

For any given tilting $A$-module $T$ with $(T1)$-$(T3)$, the module
$T':=\bigoplus_{i=0}^n T_i$ is a good $n$-tilting module which is
equivalent to the given one, that is, $T$ and $T'$ generate the same
tilting class in the category of $A$-modules (see \cite{Bz2}).

Let $T$ be an $n$-tilting $A$-module and $B$ the endomorphism ring
of $_AT$. In general, the total right-derived functor $\rHom_A(T,-)$
does not define a triangle equivalence between the (unbounded)
derived category $\D{A}$ of $A$ and the derived category $\D{B}$ of
$B$. However, if $_AT$ is good, then $\rHom_A(T,-)$ is fully
faithful and induces a triangle equivalence between the derived
category $\D{A}$ and the Verdier quotient of $\D{B}$ modulo the
kernel $\Ker(T\otimesL_{B}-)$ of the total left-derived functor
$T\otimesL_{B}-$ (see \cite[Theorem 2.2]{Bz2}). Furthermore, the
functor $\rHom_A(T,-):\D{A}\to\D{B}$ is an equivalence if and only
if $T$ is a classical tilting module if and only if
$\Ker(T\otimesL_{B}-)$ vanishes (see \cite{Bz2}). From this point of
view, the category $\Ker(T\otimesL_{B}-)$ measures the difference
between the derived categories $\D{A}$ and $\D{B}$.


Motivated by the main result in \cite{CX1}, we introduce the
following notion. A full triangulated subcategory $\mathcal X$ of
$\D{B}$ is said to be \emph{homological} if there is a homological
ring epimorphism $B\to C$ of rings such that the restriction functor
$\D{C}\to \D{B}$ induces a triangle equivalence from $\D{C}$ to
$\mathcal{X}$. Thus, if the projective dimension of a good tilting
module $_AT$ is at most $1$, then the subcategory
$\Ker(T\otimesL_{B}-)$ of $\D{B}$ is homological. Now, in terms of
homological subcategories, our question can be restated as follows:

\medskip
{\bf Question.} \emph{Is the full triangulated  subcategory
$\Ker(T\otimesL_{B}-)$ of $\D{B}$ always homological for any good
$n$-tilting $A$-module $T$ with $n\geq 2$? Here, $B$ is the
endomorphism ring of the module $T$.}

\medskip
Let us first give several characterizations for
$\Ker(T\otimesL_{B}-)$ to be homological.

\begin{Theo}\label{main-result}
Suppose that $A$ is a ring and $n$ is a natural number. Let $T$ be a
good $n$-tilting $A$-module, and let $B$ be the endomorphism ring of
$_AT$. Then the following are equivalent:

$(1)$ The full triangulated subcategory $\Ker(T\otimesL_{B}-)$ of
$\D{B}$ is homological.

$(2)$ The category consisting of the $B$-modules $Y$ with
$\Tor_m^B(T,Y)=0$ for all $m\geq0$ is an  abelian subcategory of the
category of all $B$-modules.

$(3)$ The $m$-th cohomology of the complex $\Hom_A(\cpx{P},
A)\otimes_AT_B$ vanishes for all $m\ge 2$, where the complex
$\cpx{P}$ is a deleted projective resolution of $_AT$.

$(4)$ The kernel $K$ of the homomorphism
$\Coker(\varphi_0)\lra\Coker(\varphi_1)$ induced from $\sigma:
P_1\to P_0$ in $(T1)$ satisfies $\,\Ext^m_{B\opp}(T, K)=0$ for all $
m\geq 0$, where $\varphi_i: \Hom_A(P_i, A)\otimes_AT \lra
\Hom_A(P_i, T)$ is the composition map under the identification of
$_AT_B$ with $\Hom_A(A, T)$ for $i=0, 1$.

\smallskip
In particular, if $n=2$, then $(1)$ holds if and only if
$\,\Ext^2_A(T,\,A)\otimes_AT=0$.
\end{Theo}

We remark that if the category $\Ker(T\otimesL_{B}-)$ is homological
in $\D{B}$, then the generalized  localization $\lambda: B\to B_T$
of $B$ at the module $T_B$ exists (see Definition \ref{genloc}) and
is homological, and therefore there is a recollement of derived
module categories:
$$\xymatrix@C=1.2cm{\D{B_T}\ar[r]^-{D(\lambda_*)}
&\D{B}\ar[r]^-{ {}_AT\otimesL_B-}\ar@/^1.2pc/[l]\ar@/_1.3pc/[l]
&\D{A} \ar@/^1.2pc/[l]\ar@/_1.3pc/[l]}\vspace{0.2cm}$$

\smallskip
\noindent where $D(\lambda_*)$ stands for the restriction functor
induced by $\lambda$. Thus, Theorem \ref{main-result} can be
regarded as a kind of generalization of \cite[Theorem 1.1 (1)]{CX1},
and also gives an explanation why \cite[Theorem 1.1 (1)]{CX1} holds.

As a consequence of  Theorem \ref{main-result}, we have the
following corollary in which $(1)$ extends \cite[Theorem 1.1
(1)]{CX1}, while our new contribution to $(2)$  is the necessity
part of the statement.

\begin{Koro}\label{cor}
Suppose that $A$ is a ring and $n$ is a natural number. Let $T$ be a
good $n$-tilting $A$-module, and let $B$ be the endomorphism ring of
$_AT$.

$(1)$ If $_AT$ decomposes into $M\oplus N$ such that the projective
dimension of $_AM$ is at most $1$ and that the first syzygy of $_AN$
is finitely generated, then the category $\Ker(T\otimesL_{B}-)$ is
homological.

$(2)$ Suppose that $A$ is commutative. If $\Hom_A(T_{i+1}, T_i)=0$
for all $T_i$ in $(T3)'$ with $1\leq i \leq n-1$, then the category
$\Ker(T\otimesL_{B}-)$ is homological if and only if the projective
dimension of $_AT$ is at most $1$, that is, $_AT$ is a $1$-tilting
module.
\end{Koro}

A remarkable consequence of Corollary \ref{cor} is that we can get
an answer to the above-mentioned question. In fact, in Section
\ref{7.1}, we display an example of an $n$-tilting module $T$ for
each $n\geq 2$ and shows that $\Ker(T\otimesL_{B}-)$ is not
homological.

Dually, there is the notion of (good) cotilting modules of finite
injective dimension over arbitrary rings. This notion involves
injective cogenerators of module categories. As is known, there is
no nice duality between infinitely generated tilting and cotilting
modules. This means that methods for dealing with tilting modules
may not work dually with cotilting modules. Nevertheless, we shall
use methods in this paper to deal with cotilting modules with
respect to some ``nice'' injective cogenerators. Our methods cover
particularly cotilting modules over Artin algebras. Here, our main
concern again is when the induced subcategories of derived
categories of the endomorphism rings of good cotilting modules are
homological, or equivalently, the existence of a recollement similar
to \cite[Theorem 1.1 (1)]{CX1}.

Our consideration is focused on (infinitely generated) cotilting
modules over Artin algebras $A$. Let $D$ be the usual duality of an
Artin algebra. The dual module $D(A_A)$ is an injective cogenerator
for the category of $A$-modules, and called the \emph{ordinary
injective cogenerator}. Our main result for cotilting modules is as
follows.

\begin{Theo}\label{coth}
Suppose that $A$ is an Artin algebra. Let $U$ be a good
$1$-cotilting $A$-module with respect to the ordinary injective
cogenerator for the category of $A$-modules. Set $R:=\End_A(U)$ and
$M:=\Hom_A(U, D(A))$. Then the universal localization $\lambda: R\to
R_M$ of $R$ at the module $_RM$ is homological, and there exists a
recollement of derived module categories:
$$\xymatrix@C=1.2cm{\D{R_M}\ar[r]^-{D(\lambda_*)}
&\D{R}\ar[r]\ar@/^1.2pc/[l]\ar@/_1.2pc/[l] &\D{A}
\ar@/^1.2pc/[l]\ar@/_1.2pc/[l]}\vspace{0.2cm}$$ where $D(\lambda_*)$
stands for the restriction functor induced by $\lambda$.
\end{Theo}

As is known, over an Artin algebra, each $1$-cotilting module is
equivalent to the dual of a $1$-tilting right module (see
\cite[Chapter 11, Section 4.15]{HHK}). However, we cannot get
Theorem \ref{coth} from the result \cite[Theorem 1.1 (1)]{CX1}
because the relationship between the endomorphism ring of an
infinitely generated $1$-cotlting module and the one of the
corresponding $1$-tilting right module is unknown.

For a more general formulation of Theorem \ref{coth} on higher
cotilting modules, one may see Corollary \ref{real-cotilt} and the
diagram ($\ddag$) above Corollary \ref{real-cotilt}. For higher
cotiltig modules, we also give conditions and counterexamples for
subcategories from cotilting modules not to be homological, though
additional attention is needed.

\medskip
The contents of this paper are sketched as follows. In Section 2, we
fix notation, recall some definitions and prove some homological
formulas. In Section \ref{sect3}, we introduce bireflective and
homological subcategories in derived categories of rings, and
discuss when bireflective subcategories are homological. In Section
\ref{sect4}, we introduce a new class of modules, called Ringel
modules, and establish a crucial result, Proposition
\ref{realization}, which is used not only to decide if a
bireflective subcategory  is homological, but also to investigate
higher tilting and cotilting modules in the later considerations. In
Section \ref{sect5}, we apply the results in previous sections to
good tilting modules and show Theorem \ref{main-result} as well as
Corollary \ref{cor}. At the end of this section, we point out an
example which shows that there do exist higher tilting modules
satisfying the conditions of Corollary \ref{cor} (1). In Section
\ref{sect6}, we first apply our results in Section \ref{sect4} to
cotilting modules in a general setting, and then prove Theorem
\ref{coth} for Artin algebras. It is worth noting that, for
cotilting $A$-modules $U$, recollements of $\D{\End_A(U)}$ may
depend on the choices of injective cogenerators to which the
cotilting modules are referred. In this section, we also give
conditions for the subcategories from cotilting modules not to be
homological. This is a preparation for constructing counterexamples
in the next section. In Section \ref{sect7}, we apply our results in
Section \ref{sect5} to good tilting modules $T$ over commutative
rings, and give a counterexample to show that, in general,
$\Ker(T\otimesL_B-)$ may not be realized as the derived module
category of a ring $C$ with a homological ring epimorphism $B\ra C$.
For higher cotilting modules, the same situation occurs. More
precisely, we shall use results in Section \ref{sect6} to display a
counterexample which demonstrates that, in general, the
corresponding subcategories from cotilting modules cannot be
realizable as derived module categories of rings. This section ends
with a few open questions closely related to the results in this
paper.

\section{Preliminaries\label{sect2}}
In this section, we briefly recall some definitions, basic facts and
notation used in this paper. For unexplained notation employed in
this paper, we refer the reader to \cite{CX1} and the references
therein.

\subsection{Notation}\label{sect2.1}

Let $\mathcal C$ be an additive category.

Throughout the paper, a full subcategory $\mathcal B$ of $\mathcal
C$ is always assumed to be closed under isomorphisms, that is, if
$X\in {\mathcal B}$ and $Y\in\cal C$ with $Y\simeq X$, then
$Y\in{\mathcal B}$.

Let $X$ be an object in $\mathcal{C}$. Denote by $\add(X)$ the full
subcategory of $\mathcal{C}$ consisting of all direct summands of
finite coproducts of copies of $M$. If $\mathcal{C}$ admits small
coproducts (that is, coproducts indexed over sets exist in
${\mathcal C}$), then we denote by $\Add(X)$ the full subcategory of
$\mathcal{C}$ consisting of all direct summands of small coproducts
of copies of $X$. Dually, if $\mathcal{C}$ admits small products,
then we denote by $\Prod(X)$ the full subcategory of $\mathcal{C}$
consisting of all direct summands of small products of copies of
$X$.

Given two morphisms $f: X\to Y$ and $g: Y\to Z$ in $\mathcal C$, we
denote the composite of $f$ and $g$ by $fg$ which is a morphism from
$X$ to $Z$. The induced morphisms $\Hom_{\mathcal
C}(Z,f):\Hom_{\mathcal C}(Z,X)\ra \Hom_{\mathcal C}(Z,Y)$ and
$\Hom_{\mathcal C}(f,Z): \Hom_{\mathcal C}(Y, Z)\ra \Hom_{\mathcal
C}(X, Z)$ are denoted by $f^*$ and $f_*$, respectively.

We denote the composition of a functor $F:\mathcal {C}\to
\mathcal{D}$ between categories $\mathcal C$ and $\mathcal D$ with a
functor $G: \mathcal{D}\to \mathcal{E}$ between categories $\mathcal
D$ and $\mathcal E$ by $GF$ which is a functor from $\mathcal C$ to
$\mathcal E$. Let $\Ker(F)$ and $\Img(F)$ be the kernel and image of
the functor $F$, respectively. In particular, $\Ker(F)$ is closed
under isomorphisms in $\mathcal{C}$. In this note, we require that
$\Img(F)$ is closed under isomorphisms in $\mathcal{D}$.

Suppose that $\mathcal{Y}$ is a full subcategory of $\mathcal{C}$.
Let $\Ker(\Hom_{\mathcal{C}}(-,\mathcal{Y}))$ be the left orthogonal
subcategory with respect to $\mathcal{Y}$, that is, the full
subcategory of $\mathcal{C}$ consisting of the objects $X$ such that
$\Hom_{\mathcal{C}}(X,Y)=0$ for all objects $Y$ in $\mathcal{Y}$.
Similarly, we can define the right orthogonal subcategory
$\Ker(\Hom_{\mathcal{C}}(\mathcal{Y},-))$ of $\cal C$ with respect
to $\mathcal{Y}$.

Let $\C{\mathcal{C}}$ be the category of all complexes over
$\mathcal{C}$ with chain maps, and $\K{\mathcal{C}}$ the homotopy
category of $\C{\mathcal{C}}$. As usual, we denote by
$\Cb{\mathcal{C}}$ the category of bounded complexes over $C$, and
by $\Kb{\mathcal{C}}$ the homotopy category of $\Cb{\mathcal{C}}$.
When $\mathcal{C}$ is abelian, the derived category of $\mathcal{C}$
is denoted by $\D{\mathcal{C}}$, which is the localization of
$\K{\mathcal C}$ at all quasi-isomorphisms. It is well known that
both $\K{\mathcal{C}}$ and $\D{\mathcal{C}}$ are triangulated
categories. For a triangulated category, its shift functor is
denoted by $[1]$ universally.

If $\mathcal{T}$ is a triangulated category with small coproducts,
then, for an object $U$ in $\mathcal{T}$, we denote by ${\rm
Tria}(U)$ the smallest full triangulated subcategory of
$\mathcal{T}$ containing $U$ and being closed under small
coproducts.

Suppose that $\mathcal{T}$ and $\mathcal{T}'$ are triangulated
categories with small coproducts. If $F:\mathcal{T}\ra \mathcal{T}'$
is a triangle functor which commutes with small coproducts, then
$F({\rm Tria}(U))\subseteq {\rm Tria}(F(U))$ for every object $U$ in
$\mathcal{T}$.

\subsection{Homological formulas}

In this paper, all rings considered are assumed to be associative
and with identity, and all ring homomorphisms preserve identity.
Unless stated otherwise, all modules are referred to left modules.

Let $R$ be a ring. We denote by $R\Modcat$ the category of all
unitary left $R$-modules, by $\Omega_R^n$ the $n$-th syzygy operator
of $R\Modcat$ for $n\in\mathbb{N}$, and regard $\Omega_R^0$ as the
identity operator of $R\Modcat$.

If $M$ is an $R$-module and $I$ is a nonempty set, then we denote by
$M^{(I)}$ and $M^I$ the direct sum and product of $I$ copies of $M$,
respectively. If $f: M\ra N$ is a homomorphism of $R$-modules, then
the image of $x\in M$ under $f$ is denoted by $(x)f$ instead of
$f(x)$. The endomorphism ring of the $R$-module $M$ is denoted by
$\End_R(M)$. Thus $M$ becomes a natural $R$-$\End_R(M)$-bimodule.
Similarly, if $N_R$ is a right $R$-module, then, by our convention,
$N$ is a left $(\End(N_R))\opp$- right $R$-bimodule.

As usual, we simply write $\C{R}$, $\K{R}$ and $\D{R}$ for
$\C{R\Modcat}$, $\K{R\Modcat}$ and $\D{R\Modcat}$, respectively, and
identify $R\Modcat$ with the subcategory of $\D{R}$ consisting of
all stalk complexes concentrated in degree zero. Let $\C {\pmodcat
R}$ be the full subcategory of $\C R$ consisting of those complexes
such that all of their terms are finitely generated projective
$R$-modules.

For each $n\in\mathbb{Z}$, we denote by $H^n(-):\D{R}\to R\Modcat$
the $n$-th cohomology functor. A complex $\cpx{X}$ is said to be
\emph{acyclic (or exact)} if $H^n(\cpx{X})=0$ for all
$n\in\mathbb{Z}$.

In the following, we shall recall some definitions and basic facts
about derived functors defined on derived module categories. For
more details and proofs, we refer to \cite{bn,weibel,HHK,CX3}.

Recall that $\K{R}_P$ (respectively, $\K{R}_I$) denotes the smallest
full triangulated subcategory of $\K{R}$ which

(i) contains all the bounded-above (respectively, bounded-below)
complexes of projective (respectively, injective) $R$-modules, and

(ii) is closed under arbitrary direct sums (respectively, direct
products).

Let $\K{R}_C$ be the full subcategory of $\K{R}$ consisting of all
acyclic complexes. Then $(\K{R}{_P}, \K{R}{_C})$ forms a hereditary
torsion pair in $\K{R}$ in the following sense:

\smallskip
$(a)$ Both $\K{R}{_P}$ and $\K{R}{_C}$ are full triangulated
subcategories of $\K{R}$.

$(b)$ $\Hom_{\K R}(\cpx{M}, \cpx{N})=0$ for $\cpx{M}\in\K{R}{_P}$
and $\cpx{N}\in\K{R}{_C}$.

$(c)$ For each $\cpx{X}\in \K{R}$, there exists a distinguished
triangle in $\K{R}$:
$${_p}\cpx{X} \lraf{\alpha_{\cpx X}}
\cpx{X} \lra {_c}\cpx{X}\lra (_p\cpx{X})[1]$$ such that
${_p}\cpx{X}\in \K{R}{_P}$ and  ${_c}\cpx{X}\in\K{R}{_C}.$

In particular, for each complex $\cpx{X}$ in $\K R$, the chain map
${_p}\cpx{X} \lraf{\alpha_{\cpx X}} \cpx{X} $ is a quasi-isomorphism
in $\K R$. The complex $_p\cpx{X}$ is called the \emph{projective
resolution} of $\cpx{X}$ in $\D{R}$. For example, if $X$ is  an
$R$-module, then we can choose $_pX$ to be a deleted projective
resolution of $_RX$.

Note also that the property $(b)$ implies that each
quasi-isomorphism between complexes in $\K{R}{_P}$ is an isomorphism
in $\K{R}$, that is a chain homotopy equivalence in $\K R$.

Dually, the pair $(\K{R}{_C}, \K{R}{_I})$ is a hereditary torsion
pair in $\K{R}$. This means that, for each $\cpx{X}$ in $\D{R}$,
there exists a complex ${_i}\cpx{X}\in \K{R}_I$ together with a
quasi-isomorphism $\beta_{\cpx{X}}: \cpx{X}\to{_i}\cpx{X}$. The
complex $_i\cpx{X}$ is called the \emph{injective resolution} of
$\cpx{X}$ in $\D{R}$.

More important, the composition functors
$$\K{R}_P\hookrightarrow\K{R}\lra\D{R}\quad\mbox{and}\quad
\K{R}_I\hookrightarrow\K{R}\lra\D{R}$$ are equivalences of
triangulated categories, and the canonical localization functor
$q:\K{R}\to \D{R}$ induces an isomorphism
$\Hom_{\K{R}}(\cpx{X},\cpx{Y})\lraf{\simeq}\Hom_{\D{R}}(\cpx{X},\cpx{Y})$
of abelian groups whenever either $\cpx{X}\in\K{R}_P$ or
$\cpx{Y}\in\K{R}_I$.

For a triangle functor $F:\K{R}\to\K{S}$, we define its \emph{total
left-derived functor} ${\mathbb L}F:\D{R}\to\D{S}$ by
$\cpx{X}\mapsto F(_p\cpx{X})$, and its \emph{total right-derived
functor} ${\mathbb R}F:\D{R}\to\D{S}$ by $\cpx{X}\mapsto
F(_i\cpx{X})$. Specially, if $F$ preserves acyclicity, that is,
$F(\cpx{X})$ is acyclic whenever $\cpx{X}$ is acyclic, then $F$
induces a triangle functor $D(F):\D{R}\to\D{S}$ defined by
$\cpx{X}\mapsto F(\cpx{X})$. In this case, up to natural
isomorphism, we have ${\mathbb L}F={\mathbb R}F=D(F)$, and simply
call $D(F)$ the \emph{derived functor} of $F$.

Let $\cpx{M}$ be a complex of $R$-$S$-bimodules. Then, the tensor
functor and the Hom-functor
$$\cpx{M}\cpx{\otimes}_S-:\K{S}\to\K{R}\quad\mbox{and}\quad\cpx{\Hom}_R(\cpx{M},-):\K{R}\to\K{S}$$
form a pair of adjoint triangle functors. For the concise
definitions of the tensor and Hom complex of two complexes, we
refer, for example, to \cite[Section 2.1]{CX3}. For simplicity, if
$Y\in S\Modcat$  and $X\in R\Modcat$, we denote
$\cpx{M}\cpx{\otimes}_SY$ and $\cpx{\Hom}_R(\cpx{M},\,X)$ by
$\cpx{M}\otimes_SY$ and $\Hom_R(\cpx{M},\,X)$, respectively.

Denote by $\cpx{M}\otimesL_S-$ the total left-derived functor of
$\cpx{M}\cpx{\otimes}_S-$, and by ${\mathbb R}\Hom_R(\cpx{M},-)$ the
total right-derived functor of $\cpx{\Hom}_R(\cpx{M},-)$. Note that
$\big(\cpx{M}\otimesL_S-, {\mathbb R}\Hom_R(\cpx{M},-)\big)$ is
still an adjoint pair of triangle functors.

\smallskip
The following result is freely used, but not explicitly stated in
the literature. Here, we will arrange it as a lemma for later
reference. For the idea of its proof, we refer to \cite[Generalized
Existence Theorem 10.5.9]{weibel}.

\begin{Lem}\label{homo}
Let $R$ and $S$ be rings, and let $H:\K{R}\lra \K{S}$ be a triangle
functor.

$(1)$ Define $\mathcal{L}{_H}$ to be the full subcategory of $\K{R}$
consisting of all complexes $\cpx{X}$ such that the chain map
$H(\alpha_{\cpx{X}}): H({_p}\cpx{X})\lra H(\cpx{X})$ is a
quasi-isomorphism in $\K{S}$. Then

$(i)$ $\mathcal{L}{_H}$ is a triangulated subcategory of $\K{R}$
containing $\K{R}{_P}$.

$(ii)$ $\mathcal{L}{_H}\cap{\K{R}{_C}} = \{\cpx{X}\in\K{R}{_C}\mid
H(\cpx{X})\in\K{S}{_C}\}$.

$(iii)$ There exists a commutative diagram of triangle functors:
$$
\xymatrix{\K{R}{_P}\ar[r]^-{\simeq}\ar[d]_-{\simeq} & \D{R}\ar[d]^-{{\mathbb L}H}\\
\mathcal{L}{_H}/\mathcal{L}{_H}\cap{\K{R}{_C}}\ar[r]^-{D(H)} &
\D{S}}
$$
where $\mathcal{L}{_H}/\mathcal{L}{_H}\cap{\K{R}{_C}}$ denotes the
Verdier quotient of $\mathcal{L}{_H}$ by
$\mathcal{L}{_H}\cap{\K{R}{_C}}$, and where $D(H)$ is defined by
$\cpx{X}\mapsto H(\cpx{X})$ for $\cpx{X}\in\mathcal{L}{_H}$.

$(2)$ Define $\mathcal{R}{_H}$ to be the full subcategory of $\K{R}$
consisting of all complexes $\cpx{X}$ such that the chain map
$H(\beta_{\cpx{X}}): H(\cpx{X})\to H({_i}\cpx{X})$ is a
quasi-isomorphism in $\K{S}$. Then

$(i)$ $\mathcal{R}{_H}$ is a triangulated subcategory of $\K{R}$
containing $\K{R}{_I}\,$.

$(ii)$ $\mathcal{R}{_H}\cap{\K{R}{_C}} = \{\cpx{X}\in\K{R}{_C}\mid
H(\cpx{X})\in\K{S}{_C}\}$.

$(iii)$ There exists a commutative diagram of triangle functors:
$$
\xymatrix{\K{R}{_I}\ar[r]^-{\simeq}\ar[d]_-{\simeq} & \D{R}\ar[d]^-{{\mathbb R}H}\\
\mathcal{R}{_H}/\mathcal{R}{_H}\cap{\K{R}{_C}}\ar[r]^-{D(H)} &
\D{S}}
$$
where $\mathcal{R}{_H}/\mathcal{R}{_H}\cap{\K{R}{_C}}$ denotes the
Verdier quotient of $\mathcal{R}{_H}$ by
$\mathcal{R}{_H}\cap{\K{R}{_C}}$, and where $D(H)$ is defined by
$\cpx{X}\mapsto H(\cpx{X})$ for $\cpx{X}\in\mathcal{R}{_H}$.
\end{Lem}

\medskip
Note that if $H$ commutes with arbitrary direct sums, then
$\mathcal{L}{_H}$ is closed under arbitrary direct sums in $\K R$.
Dually, if $H$ commutes with arbitrary direct products, then
$\mathcal{R}{_H}$ is closed under arbitrary direct products in $\K
R$.

From Lemma \ref{homo}, we see that, up to natural isomorphism, the
action of the functor ${\mathbb L}H$ (respectively, ${\mathbb R}H$)
on a complex $\cpx{X}$ in $\mathcal{L}{_H}$ (respectively,
$\mathcal{R}{_H}$) is the same as that of the functor $H$ on
$\cpx{X}$. Based on this point of view, we obtain the following
result which will be applied in our later proofs.

\begin{Koro}   \label{counit}
Let $R$ and $S$ be two rings. Suppose that $(F, G)$ is a pair of
adjoint triangle functors with $F:\K{S}\to \K{R}$ and $G:\K{R}\to
\K{S}$. Let $\theta: FG\to Id_{\K R}$ and $\varepsilon: ({\mathbb
L}F) ({\mathbb R}G) \to Id_{\D R}$ be the counit adjunctions. If
$\cpx{X}\in \mathcal{R}{_G}$ and $G(\cpx{X})\in \mathcal{L}{_F}$,
then there exists a commutative diagram in $\D{R}$:
$$
\xymatrix
{({\mathbb L}F)({\mathbb R}G)(\cpx{X})\ar[r]^-{\varepsilon_{\cpx{X}}}\ar[d]_-{\simeq} & \cpx{X}\ar@{=}[d]\\
F G(\cpx{X}) \ar[r]^-{\theta_{\cpx{X}}} & \cpx{X}}
$$
\end{Koro}

{\it Proof.} It follows from $\cpx{X}\in \mathcal{R}{_G}$ that the
quasi-isomorphism $\beta_{\cpx{X}}: \cpx{X}\to {_i}\cpx{X}$ in
$\K{R}$ induces a quasi-isomorphism $G(\beta_{\cpx{X}}):
G(\cpx{X})\to G({_i}\cpx{X})$ in $\K{S}$. Since $(\K{S}{_P},
\K{S}{_C})$ is a hereditary torsion pair in $\K{S}$, there exists a
homomorphism $_pG(\beta_{\cpx{X}}): {_p}G(\cpx{X})\to
{_p}G({_i}\cpx{X})$ in $\K{S}$ such that the following diagram is
commutative:
$$
\xymatrix{
_pG(\cpx{X})\ar[r]^-{\alpha_{\,G(\cpx{X})}}\ar[d]_-{_pG(\beta_{\cpx{X}})}\;
&\;G(\cpx{X})\ar[d]^-{G(\beta_{\cpx{X}})}\\
{_p}G(_i\cpx{X})\ar[r]^-{\alpha_{\,G(_i\cpx{X})}} & \,G(_i\cpx{X})}
$$
Note that  $_pG(\beta_{\cpx{X}})$ is a quasi-isomorphism in $\K{S}$
since all the other chain maps in the above diagram are
quasi-isomorphisms. By the property $(b)$ related to the pair
$(\K{S}{_P}, \K{S}{_C})$, we know that $_pG(\beta_{\cpx{X}})$ is an
isomorphism in $\K{S}$, and therefore the chain map
$F(_pG(\beta_{\cpx{X}})): F({_p}G(\cpx{X}))\lra
F({_p}G({_i}\cpx{X}))$ is an isomorphism in $\K{R}$.

Now, we can easily construct the following commutative diagram in
$\K{R}$:
$$
\xymatrix{&
F(_pG(\cpx{X}))\ar[r]^-{F(\alpha_{\,G(\cpx{X})})}\ar[d]_-{F(_pG(\beta_{\cpx{X}}))}^-{\simeq}
& \;FG(\cpx{X}) \ar[r]^-{\theta_{\cpx{X}}} \ar[d]^-{FG(\beta_{\cpx{X}})} &\cpx{X}\ar[d]^-{\beta_{\cpx{X}}}\\
({\mathbb L}F)({\mathbb R}G)(\cpx{X})\ar@{=}[r] &
F({_p}G({_i}\cpx{X}))\ar[r]^-{F(\alpha_{\,G({_i}\cpx{X})})}\;\,&
\;\,FG({_i}\cpx{X})\ar[r]^-{\theta {_{_i\cpx{X}}}} & {_i}\cpx{X} }
$$
Since $G(\cpx{X})\in\mathcal{L}{_F}$ by assumption, the chain map
$F(\alpha_{\,G(\cpx{X})})$ is a quasi-isomorphism in $\K{R}$, and is
an isomorphism in $\D{R}$. Clearly, the quasi-isomorphism
$\beta_{\cpx{X}}$ is an isomorphism in $\D{R}$.

Furthermore, the counit $\varepsilon_{\cpx{X}}: ({\mathbb
L}F)({\mathbb R}G)(\cpx{X})\lra \cpx{X}$ is actually given by the
composite of the following homomorphisms in $\D R$:
$$\xymatrix{
({\mathbb L}F)({\mathbb R}G)(\cpx{X})\ar@{=}[r] &
F({_p}G({_i}\cpx{X}))\ar[r]^-{F(\alpha_{G({_i}\cpx{X})})}\;&
\;FG({_i}\cpx{X})\ar[r]^-{\theta_{_i\cpx{X}}} & _i\cpx{X}
\ar[r]^-{\;(\beta_{\cpx{X}})^{-1}} & \cpx{X}.}
$$

Define $$\tau=
\big(F(_pG(\beta_{\cpx{X}}))\big)^{-1}\,F(\alpha_{\,G(\cpx{X})}):\;
({\mathbb L}F)({\mathbb R}G)(\cpx{X})\lra F G(\cpx{X})$$ which is an
isomorphism in $\D{R}$. It follows that there exists a commutative
diagram in $\D{R}$:
$$ \xymatrix
{({\mathbb L}F)({\mathbb R}G)(\cpx{X})\ar[r]^-{\varepsilon_{\cpx{X}}}\ar[d]^-{\tau} & \cpx{X}\ar@{=}[d]\\
F G(\cpx{X}) \ar[r]^-{\theta_{\cpx{X}}} & \cpx{X}}
$$
This finishes the proof. $\square$

\smallskip
As a preparation for our later proofs, we mention the following
three homological formulas which are related to derived functors or
total derived functors. The first one is taken from \cite[Theorem
3.2.1, Theorem 3.2.13, Remark 3.2.27]{EJ}.

\begin{Lem}\label{tor-ext}
Let $R$ and $S$ be  rings. Suppose that $M$ is an $S$-$R$-bimodule
and $I$ is an injective  $S$-module.

$(1)$ If $N$ is an $R$-module, then
$$
\Hom_S(\Tor_i^R(M, \,N),\, I)\simeq \Ext^i_R(N,\, \Hom_S(M,\,
I))\;\, \mbox{for all}\;\, i\geq 0.$$

$(2)$ If $L$ is an $R\opp$-module which has a finitely generated
projective resolution in $R\opp\Modcat$, then
$$
\Hom_S(\Ext_R^i(L,\,M),\, I)\simeq \Tor^R_i(L, \,\Hom_S(M, I))\;\,
\mbox{for all}\;\, i\geq 0.$$
\end{Lem}

\smallskip
The next formula is proved in \cite[Section 2.1]{CX3}.

\begin{Lem}\label{complex}
Let $R$ and $S$ be rings. Suppose that $\cpx{X}$ is a bounded
complex of $R$-$S$-bimodules. If $\cpx{X}\in \Cb{\pmodcat{R}}$, then
there is a natural isomorphism of functors:
$$\Hom_R(\cpx{X}, R)\cpx{\otimes}_R-\,\lraf{\simeq}\cpx{\Hom}_R(\cpx{X},-):\C{R}\to\C{S}.$$
In particular,
$$\Hom_R(\cpx{X}, R)\otimesL_R-\,\lraf{\simeq}\rHom{_R}(\cpx{X},-):\D{R}\to\D{S}.$$
\end{Lem}

\smallskip
The last formula is useful for us to calculate the cohomology groups
of tensor products of complexes.

\begin{Lem}\label{Formula} Let $n$ be an integer, and let $S$ be a ring and $M$ an $S\opp$-module. Suppose that
$\cpx{Y}:=(Y^i)_{i\in\mathbb{Z}}$ is a complex in $\C S$ such that
$Y^i=0$ for all $i\geq n+1$, and $\Tor^S_j(M, Y^i)=0$ for all
$i\in\mathbb{Z}$ and $j\geq 1$. Let $m\in\mathbb{Z}$ with $m < n$.
If $\,\Tor^S_t\big(M, H^{m+t}(\cpx{Y})\big)=0=\Tor^S_{t-1}\big(M,
H^{m+t}(\cpx{Y})\big)$ for $0\leq t \leq n-m-1$, then
$H^m(M\otimes_S\cpx{Y})\simeq \Tor^S_{n-m}\big(M,
H^n(\cpx{Y})\big)$.
\end{Lem}

{\it Proof.} Suppose that $\cpx{Y}$ is  the following form:
$$
\cdots\lra Y^{m-1} \lraf{d^{m-1}} Y^{m} \lraf{d^m} Y^{m+1}\lra
\cdots\lra Y^{n-1} \lraf{d^{n-1}} Y^n\lra 0\lra \cdots$$ For
$i\in\mathbb{Z}$, define $C_i:=\Coker(d^{i-1})=Y^i/\Img(d^{i-1})$
and $I_i:=\Img(d^i)$. Then we have two short exact sequences of
$S$-modules for each $i\in \mathbb{Z}$:
$$
(a)\quad 0\lra H^i(\cpx{Y})\lra C_i\lraf{\pi_i} I_i\lra 0
\quad\mbox{and}\quad (b) \quad 0\lra I_i
\stackrel{\lambda_i}{\hookrightarrow} Y^{i+1}\lra C_{i+1}\lra 0.$$
Clearly, $H^i(\cpx{Y})=\Ker(\pi_i\lambda_i)$, and  $d^{i}: Y^i\to
Y^{i+1}$ is just the composite of the canonical surjection $Y^i \to
C_i$ with $\pi_i\lambda_i: C_i \to Y^{i+1}$.

$(1)$ We claim that if $M\otimes_S H^i(\cpx{Y})=0$, then
$H^i(M\otimes_S\cpx{Y})\simeq \Tor^S_1(M, C_{i+1})$.

In fact, since $M\otimes_S-: S\Modcat\to \mathbb{Z}\Modcat$ is right
exact, the sequence
$$M\otimes_S Y^{i-1}\lraf{1\otimes d^{i-1}} M\otimes_S
Y^i\lra M\otimes_S C_i\lra 0$$ is exact, that is, $\Coker(1\otimes
d^{i-1})\simeq  M\otimes_S C_i$. This implies that
$H^i(M\otimes_S\cpx{Y})\simeq \Ker(1\otimes \pi_i\lambda_i)$ where
$$1\otimes \pi_i\lambda_i=(1\otimes \pi_i) (1\otimes \lambda_i):
M\otimes_S C_i\to M\otimes_S Y^{i+1}, $$ which is the composite of
$1\otimes \pi_i: M\otimes_S C_i\lra M\otimes_S I_i$  with $1\otimes
\lambda_i: M\otimes_S I_i\lra M\otimes_S Y^{i+1}$.

Assume that $M\otimes_S H^i(\cpx{Y})=0$. Then $1\otimes\pi_i$ is an
isomorphism and $\Ker(1\otimes \pi_i\lambda_i)\simeq \Ker(1\otimes
\lambda_i)$. Now, we apply $M\otimes_S-$ to the sequence $(b)$, and
get the following exact sequence:
$$
\Tor^S_1(M, Y^{i+1})\lra \Tor^S_{1}(M, C_{i+1}) \lra M\otimes_S I_i
\lraf{1\otimes \lambda_i}  M\otimes_S Y^{i+1}
$$
Since $\Tor^S_1(M, Y^{i+1})=0$ by assumption, we obtain
$\Tor^S_{1}(M, C_{i+1})\simeq \Ker(1\otimes \lambda_i).$ It follows
that $$H^i(M\otimes_S\cpx{Y})\simeq \Ker(1\otimes
\pi_i\lambda_i)\simeq \Ker(1\otimes \lambda_i)\simeq \Tor^S_1(M,
C_{i+1}).$$ This finishes the claim (1).

$(2)$ We show that, for any $j\geq 1$, if $\Tor^S_j(M,
H^i(\cpx{Y}))=0=\Tor^S_{j-1}(M, H^i(\cpx{Y}))$, then
$$\Tor^S_j(M,\,C_i)\lraf{\simeq} \Tor^S_{j+1}(M, C_{i+1}).$$

This follows from applying $M\otimes_S-$ to the exact sequences
($a$) and ($b$), respectively, together with our assumptions on
$\cpx{Y}$.

$(3)$ Let $m\in\mathbb{Z}$ with $m\leq n-1$. Suppose that
$$\Tor^S_t\big(M, H^{m+t}(\cpx{Y})\big)=0=\Tor^S_{t-1}\big(M,
H^{m+t}(\cpx{Y})\big)\quad \mbox{for}\quad 0\leq t \leq n-m-1. $$
Then, by taking $t=0$, we have $M\otimes_S H^m(\cpx{Y})=0$.  Thanks
to (1), we have $H^m(M\otimes_S\cpx{Y})\simeq \Tor^S_1(M, C_{m+1})$.
Since $Y^i=0$ for $i\geq n+1$, it follows that $H^n(\cpx{Y})= C_n$.
This implies that if $n-m=1$, then $H^m(M\otimes_S\cpx{Y})\simeq
\Tor^S_{n-m}\big(M, H^n(\cpx{Y})\big)$.

Now, suppose $n-m\geq 2$. For $1\leq t\leq n-m-1$, we see from $(2)$
that $\Tor^S_t(M,\,C_{m+t})\lraf{\simeq} \Tor^S_{t+1}(M,
C_{m+t+1}).$ Thus
$$\Tor^S_1(M,\,C_{m+1})\simeq \Tor^S_{2}(M,
C_{m+2})\simeq\cdots\simeq \Tor^S_{n-m-1}(M,\,C_{n-1})\simeq
\Tor^S_{n-m}(M, C_{n}).
$$
Consequently, $H^m(M\otimes_S\cpx{Y})\simeq \Tor^S_1(M,
C_{m+1})\simeq\Tor^S_{n-m}(M, C_{n})= \Tor^S_{n-m}(M,
H^n(\cpx{Y}))$. This finishes the proof of Lemma \ref{Formula}.
$\square$

\subsection{Relative Mittag-Leffler modules} \label{sect2.3}
\medskip
Now, we recall the definition of relative Mittag-Leffler modules
(see \cite{Go}, \cite{HH}).

\begin{Def}\label{ML}{\rm
A right $R$-module $M$ is said to be \emph{$R$-Mittag-Leffler} if
the canonical map
$$\rho{_I}:\; M\otimes_RR^I\lra M^I,\; m\otimes(r_i)_{i\in I}\mapsto
(mr_i)_{i\in I} \, \mbox{ for } \; m\in M, \; r_i\in R,$$ is
injective for any nonempty set $I$.

A right $R$-module $M$ is said to be \emph{strongly
$R$-Mittag-Leffler} if the $m$-th syzygy  of $M$ is
$R$-Mittag-Leffler for every $m\geq 0$.}
\end{Def}

By \cite[Theorem 1]{Go}, a right $R$-module $M$ is
$R$-Mittag-Leffler if and only if, for any finitely generated
submodule $X$ of $M_R$, the inclusion $X\to M$ factorizes through a
finitely presented right $R$-module. This implies that if $M$ is
finitely presented, then it is $R$-Mittag-Leffler. Actually, for
such a module $M$, the above map $\rho_I$ is always bijective (see
\cite[Theorem 3.2.22]{EJ}). Further, if the ring $R$ is right
noetherian, then each right $R$-module is $R$-Mittag-Leffler since
each finitely generated right $R$-module is finitely presented.

In the next lemma, we shall collect some basic properties of
Mittag-Leffler modules for later use.

\begin{Lem}\label{MLP}
Let $R$ be a ring and $M$ a right $R$-module. Then the following
statements are true.

$(1)$ If $M$ is $R$-Mittag-Leffler, then so is each module in
$\Add(M_R)$. In particular, each projective right $R$-module is
$R$-Mittag-Leffler.

$(2)$ The first syzygy  of $M$ in $R{^{\rm op}}\Modcat$ is
$R$-Mittag-Leffler  if and only if  $\Tor^R_1(M, R^I)=0$ for every
nonempty set $I$.

$(3)$ $M$ is strongly $R$-Mittag-Leffler if and only if $M$ is
$R$-Mittag-Leffler and $\Tor^R_i(M, R^I)=0$ for each $i\geq 1$ and
every nonempty set $I$.

$(4)$ If $M$ is finitely generated, then $M$ is strongly
$R$-Mittag-Leffler if and only if $M$ has a finitely generated
projective resolution.
\end{Lem}

{\it Proof.}  $(1)$ follows from the fact that tensor functors
commute with direct sums.

$(2)$ Note that the first syzygy $\Omega_R(M)$ of $M$ depends on the
choice of projective presentations of $M_R$. However, the
``$R$-Mittag-Leffler'' property of $\Omega_R(M)$ is independent of
the choice of projective presentations of $M_R$. This is due to (1)
and Schanuel's Lemma in homological algebra.

So, we choose an exact sequence $$0\lra K_1\lraf{f} P_1\lra M\lra
0$$ of right $R$-modules with $P_1$ projective, and shall show that
$K_1$ is $R$-Mittag-Leffler if and only if $\Tor^R_1(M, R^I)=0$ for
any nonempty set $I$.

Obviously, we can construct the following exact commutative diagram:
$$\xymatrix{
0\ar[r] & \Tor^R_1(M, R^I)\ar[r] & K_1\otimes_RR^I\ar[r]^-{f\otimes
1}\ar[d]^-{\rho_2} & P_1\otimes_RR^I\ar[r]\ar[d]^-{\rho_1} & M\otimes_RR^I\ar[r]\ar[d] & 0\\
 & 0\ar[r] &  K_1{^I}\ar[r]^{f^I} & P_1{^I} \ar[r] & M{^I} \ar[r] & 0 }
$$
where $\rho_i$, $1\le i\le 2$, are the canonical maps (see
Definition \ref{ML}). Since the projective module $P_1$ is
$R$-Mittag-Leffler by $(1)$, the map $\rho_1$ is injective. This
means that $\rho_2$ is injective if and only if so is $f\otimes 1$.
Clearly, the former is equivalent to that $K_1$ is
$R$-Mittag-Leffler, while the latter is equivalent to that
$\Tor^R_1(M, R^I)=0$. This finishes the proof of $(2)$.

$(3)$ For each $i\geq 0$, let $\Omega^i_R(M)$ stand for the  $i$-th
syzygy  of $M$ in $R{^{\rm op}}\Modcat$. Then, for each nonempty set
$I$, we always have
$$\Tor^R_{i+1} (M, R^I)\simeq \Tor^R_1 (\Omega^{i}_R(M),\, R^I).$$
Now $(3)$ follows immediately from $(2)$.

$(4)$ The sufficient condition is clear. Now suppose that $M$ is
strongly $R$-Mittag-Leffler. We need only to show that the first
syzygy of $M$ is finitely generated, that is, $M$ is finitely
presented. However, this follows from the fact that the inclusion
map $M\hookrightarrow M$ factorizes through a finitely presented
right $R$-module. $\square$

\medskip
A special class of strongly Mittag-Leffler modules is the class of
tilting modules. The following result can be concluded from
\cite[Corollary 9.8]{HH}, which will play an important role in our
proof of the main result.

\begin{Lem}\label{TML}
If $M$ is a tilting right $R$-module, then $M$ is strongly
$R$-Mittag-Leffler.
\end{Lem}

As a corollary of Lemmas \ref{TML} and \ref{MLP} (4), we obtain the
following result which is a generalization of \cite[Corollary
4.7]{CX1}.

\begin{Koro}\label{2.9}
Let $M$ be a tilting right $R$-module. If $M$ is finitely generated,
then $M$ is classical.
\end{Koro}

{\bf Proof.} Suppose that $M_R$ is finitely generated. Then we can
get an exact sequence $(T3)'$ from $(T3)$ by using the argument in
\cite[Corollary 4.7]{CX1} repeatedly. This shows that $M_R$ is
actually a good tilting module. Since $M$ is strongly
$R$-Mittag-Leffler, it follows from Lemma \ref{MLP} (4) that $M$
admits a finitely generated projective resolution. Clearly, such a
resolution can be chosen to be of finite length since $M$ has finite
projective dimension. This implies that $M_R$ is classical.
$\square$

\section{Homological subcategories of derived module categories}\label{sect3}

In this section, we shall give the definitions of bireflective and
homological subcategories of derived module categories. In
particular, we shall establish some applicable criterions for
bireflective subcategories to be homological.

Let $R$ and $S$ be arbitrary rings.

Let $\lambda: R\ra S$ be a homomorphism of rings. We denote by
$\lambda_*:S\Modcat\to R\Modcat$ the restriction functor induced by
$\lambda$, and by $D(\lambda_*):\D{S}\to\D{R}$ the derived functor
of the exact functor $\lambda_*$. Recall that $\lambda$ is a
\emph{ring epimorphism} if $\lambda_*:S\Modcat\to R\Modcat$ is fully
faithful. This is equivalent to saying that the multiplication map
$S\otimes_RS\ra S$ is an isomorphism in $S\Modcat$.

Two ring epimorphisms $\lambda: R\to S$ and $\lambda': R\to S'$ are
said to be \emph{equivalent} if there is an isomorphism $\psi: S\to
S'$ of rings such that $\lambda'=\lambda\psi$. Note that there is a
bijection between the equivalence classes of ring epimorphisms
staring from $R$ and bireflective full subcategories of $R\Modcat$,
and that there is a bijection between bireflective full
subcategories of $R\Modcat$ and the abelian full subcategories of
$R\Modcat$ which are closed under arbitrary direct sums and direct
products (see, for example, \cite[Lemma 2.1]{CX1}).

Recall that a ring epimorphism $\lambda: R\ra S$ is called
\emph{homological} if $\;\Tor^R_i(S, S)=0$ for all $i>0$. This is
equivalent to that the functor $D(\lambda_*):\D{S}\to\D{R}$ is fully
faithful, or that $S\otimesL_RS\simeq S$ in $\D{S}$. It is known
that $D(\lambda_*)$ has a left adjoint $S\otimesL_R-$ and a right
adjoint $\rHom_R(S, -)$.

\label{def-homo} Let $\mathcal{Y}$ be a full triangulated
subcategory of $\D{R}$. We say that $\mathcal{Y}$ is
\emph{bireflective} if the inclusion $\mathcal{Y}\to \D{R}$ admits
both a left adjoint and a right adjoint.

Combining \cite[Chapter I, Proposition 2.3]{BI} with \cite[Section
2.3]{CX1}, we know that a full triangulated  subcategory
$\mathcal{Y}$ of $\D{R}$ is bireflective if and only if there exists
a recollement of triangulated categories of the form
$$
\xymatrix@C=1.2cm{\mathcal{Y}\ar[r]^-{{i_*}}
&\D{R}\ar[r]\ar@/^1.2pc/[l]\ar@/_1.2pc/[l]
&\mathcal{X}\ar@/^1.2pc/[l]\ar@/_1.2pc/[l] }\vspace{0.2cm}$$ where
$i_*$ is the inclusion functor. Here, by a recollement of
triangulated categories (see \cite{BBD}) we mean that there are six
triangle functors between triangulated categories in the following
diagram:
$$\xymatrix{\mathcal{Y}\ar^-{i_*=i_!}[r]&\D{R}\ar^-{j^!=j^*}[r]
\ar^-{i^!}@/^1.2pc/[l]\ar_-{i^*}@/_1.6pc/[l]
&\mathcal{X}\ar^-{j_*}@/^1.2pc/[l]\ar_-{j_!}@/_1.6pc/[l]}$$ such
that

$(1)$ $(i^*,i_*),(i_!,i^!),(j_!,j^!)$ and $(j^*,j_*)$ are adjoint
pairs,

$(2)$ $i_*,j_*$ and $j_!$ are fully faithful functors,

$(3)$ $i^!j_*=0$ (and thus also $j^! i_!=0$ and $i^*j_!=0$), and

$(4)$ for each object $X\in\D{R}$, there are two canonical
distinguished triangles in $\D{R}$:
$$
i_!i^!(X)\lra X\lra j_*j^*(X)\lra i_!i^!(X)[1],\qquad j_!j^!(X)\lra
X\lra i_*i^*(X)\lra j_!j^!(X)[1],
$$ where $i_!i^!(X)\ra X$ and $j_!j^!(X)\ra
X$ are counit adjunction morphisms, and where $X\ra j_*j^*(X)$ and
$X\ra i_*i^*(X)$ are unit adjunction morphisms.

Note that $\mathcal{X}$ is always equivalent to the full subcategory
$\Ker\big(\Hom_{\D R}(-, \mathcal{Y})\big)$ of $\D{R}$ as
triangulated categories ( for example, see \cite[Lemma 2.6]{CX1}).
But here we do not require that the triangulated category
$\mathcal{X}$ must be a subcategory of $\D{R}$ in general. For more
examples of recollements related to homological ring epimorphisms,
we refer the reader to \cite{CX2}.

\medskip
Clearly, if $\mathcal{Y}$ is homological (see Definition in Section
\ref{sect1}), then it is bireflective. Let us now consider the
converse of this statement.

\smallskip
From now on, we assume that $\mathcal{Y}$ is a \textbf{bireflective
subcategory} of $\D{R}$, and define $\mathscr{E}:=\mathcal{Y}\cap
R\Modcat$.

It is easy to see that $\mathcal{Y}$ is closed under isomorphisms,
arbitrary direct sums and direct products in $\D{R}$. This implies
that $\mathscr{E}$ also has the above properties in $R\Modcat$.
Moreover, $\mathscr{E}$ always admits the ``$2$ out of $3$"
property: For an arbitrary short exact sequence in $R\Modcat$, if
any two of its three terms belong to $\mathscr{E}$, then the third
one belongs to $\mathscr{E}$. By \cite[Lemma 2.1]{CX1},
$\mathscr{E}$ is an abelian subcategory of $R\Modcat$ if and only if
$\mathscr{E}$ is closed under kernels (respectively, cokernels) in
$R\Modcat$. This is also equivalent to saying that there exists a
unique ring epimorphism $\lambda:R\to S$ (up to equivalence) such
that $\mathscr{E}$ is equal to $\Img(\lambda_*)$.

If $\mathcal{Y}$ is homological via a homological ring epimorphism
$\lambda: R\to S$, then $\mathcal{Y}=\Img\big(D(\lambda_*)\big)$ and
$\mathscr{E}=\Img(\lambda_*)$. In this case, $\mathscr{E}$ must be a
full, abelian subcategory of $R\Modcat$.

In general, for a bireflective subcategory $\mathcal Y$ in $\D{R}$,
the category $\mathscr{E}$ may not be abelian. This means that
bireflective subcategories in $\D{R}$ may not be homological.
Alternatively, we can reach this point by looking at differential
graded rings: By the proof of \cite[Chapter IV, Proposition
1.1]{BI}, the complex $i^*(R)$ is a compact generator of
$\mathcal{Y}$. In particular, we have
$\mathcal{Y}={\rm{Tria}}(i^*(R))$. It follows from \cite[Chapter 5,
Theorem 8.5]{HHK} that there exists a dg (differential graded) ring
such that its dg derived category is equivalent to $\mathcal{Y}$ as
triangulated categories. In general, this dg ring may have
non-trivial cohomologies in  other degrees besides the degree $0$.
In other words, the category $\mathcal{Y}$ may not be realized by
the derived module category of an ordinary ring.

Let  $i_*:\mathcal{Y}\to \D{R}$ be the inclusion functor with $i^*:
\D{R}\to \mathcal{Y}$ as its left adjoint. Define $\Lambda:=\End_{\D
R}(i^*(R))$. Then, associated with $\mathcal{Y}$, there is a ring
homomorphism defined by
$$ \delta:\; R\lra\Lambda,\quad r\mapsto i^*(\cdot
r) \, \mbox{ for } \, r\in R,$$ where $\cdot r: R\to R$ is the right
multiplication by $r$ map. This ring homomorphism induces a functor
$$\delta_*: \Lambda\Modcat\lra R\Modcat,$$
called the restriction functor.

\medskip
The following result is motivated by \cite[Section 6 and Section
7]{nr1}.

\begin{Lem} \label{ringhom}
The following statements hold true.

$(1)$ For each $\cpx{Y}\in\mathcal{Y}$, we have $H^n(\cpx{Y})\in
\Img(\delta_*)$ for all $n\in\mathbb{Z}$.  In particular,
$H^n(i^*(R))$ is an $R$-$\Lambda$-bimodule for all $n\in
\mathbb{Z}$.

$(2)$ Let $\eta_R :R\to i_*i^*(R)$ be the unit adjunction morphism
with respect to the adjoint pair $(i^*, i_*)$. Then $\Lambda\simeq
H^0(i^*(R))$ as $R$-$\Lambda$-bimodules, and there exists a
commutative diagram of $R$-modules:
$$\xymatrix{R \ar[r]^-{\delta}\ar[rd]_-{H^0(\eta_R)} & \Lambda\ar[d]^-{\simeq}\\
& H^0(i^*(R))}$$

$(3)$ If $H^0(i^*(R))\in\mathcal{Y}$, then $H^n(i^*(R))=0$ for all
$n\geq 1$, the homomorphism $\delta$ is a ring epimorphism and
$$\mathcal{Y}=\{\cpx{Y}\in\D{R}\mid H^m(\cpx{Y})\in \Img(\delta_*)
\mbox{\;for\, all \,} m\in\mathbb{Z}\}.$$

\end{Lem}

\medskip
{\it Proof.} The proof of Lemma \ref{ringhom} is derived from
\cite[Section 6 and Section 7]{nr1}, where $\mathcal{Y}$ is related
to a set of two-term complexes in $\C{\pmodcat R}$.

By our convention, the full subcategory $\Img(\delta_*)$ of
$R\Modcat$ is required to be closed under isomorphisms in
$R\Modcat$.

Let $\eta_R :R\to i_*i^*(R)=i^*(R)$ be the unit adjunction morphism.

$(1)$ Let $\cpx{Y}\in\mathcal{Y}$. Then we obtain the following
isomorphisms for each $n\in\mathbb{Z}$:
$$
\Hom_{\D{R}}(i^*(R), \cpx{Y}[n])\lraf{\simeq} \Hom_{\D R}(R,
i_*(\cpx{Y}) [n])=\Hom_{\D{R}}(R, \cpx{Y}[n])\simeq H^n(\cpx{Y}),
$$
where the first isomorphism is given by $\Hom_{\D{R}}(\eta_R,
\cpx{Y}[n])$, which is actually an isomorphism of $R$-modules. Since
$\Hom_{\D R}(i^*(R), \cpx{Y}[n])$ is a left $\Lambda$-module, we
clearly have $H^n(\cpx{Y})\in \Img(\delta_*)$. If $\cpx{Y}= i^*(R)$,
then one can check that the composite of the following isomorphisms
$$ (\ast)\quad
\Hom_{\D{R}}(i^*(R), i^*(R)[n])\simeq \Hom_{\D R}(R,
i_*i^*(R)[n]))=\Hom_{\D{R}}(R, i^*(R)[n])\simeq H^n(i^*(R))$$ is an
isomorphism of $R$-$\Lambda$-bimodules. This implies that
$H^n(i^*(R))$ is an $R$-$\Lambda$-bimodule.

$(2)$ In $(\ast)$, we take $n=0$. This gives the first part of
$(2)$. For the second part of $(2)$, we note that there exists the
following commutative diagram of $R$-modules:
$$\xymatrix{\Hom_R(R, R) \ar[r]^-{i^*}\ar[rd]_-{\Hom_{\D R}(R, \eta_R)\qquad}
& \Hom_{\D R}(i^*(R), i^*(R))\ar[d]^-{\simeq}\\
& \Hom_{\D R}(R, i_*i^*(R))}$$ which implies the diagram in $(2)$ if
we identify $\Hom_R(R, R)$,\, $\Hom_{\D R}(R, i_*i^*(R))$ and
$\Hom_{\D R}(R, \eta_R)$  with $R$,\, $H^0(i^*(R))$ and
$H^0(\eta_R)$, respectively.

$(3)$ Define $$\mathcal{Y}':=\{\cpx{Y}\in\D{R}\mid H^m(\cpx{Y})\in
\Img(\delta_*) \mbox{\;for\, all \,} m\in\mathbb{Z}\}.$$ It follows
from $(1)$ that $\mathcal{Y}\subseteq \mathcal{Y}'$.

Suppose $H^0(i^*(R))\in\mathcal{Y}$. We shall prove that
$\mathcal{Y}'\subseteq \mathcal{Y}$, and so
$\mathcal{Y}=\mathcal{Y}'$.

In fact, from $(2)$ we see that $\Lambda\simeq H^0(i^*(R))$ as
$R$-modules, and so $_R\Lambda\in\mathcal{Y}$. Note that the derived
functor $D(\delta_*): \D{\Lambda}\to \D{R}$ admits a right adjoint,
and therefore it commutes with arbitrary direct sums. Since
$\mathcal{Y}$ is a full triangulated subcategory of $\D{R}$ closed
under arbitrary direct sums in $\D{R}$,  it follows from
$\D{\Lambda}={\rm{Tria}}(_\Lambda\Lambda)$ and
$_R\Lambda\in\mathcal{Y}$ that $\Img\big( D(\delta_*)
\big)\subseteq\mathcal{Y}.$ In particular, $\Img(\delta_*)\subseteq
\mathcal{Y}$.

To prove $\mathcal{Y}'\subseteq \mathcal{Y}$, we shall use the
following statements $(a)$-$(d)$ mentioned in \cite[Lemma 4.6]{HKL}.
For the definitions of homotopy limits and homotopy colimits in
triangulated categories, we refer to \cite[Section 2]{bn}.

$(a)$ By canonical truncations, one can show that each bounded
complex over $R$ can be generated by its cohomologies, that is, if
$\cpx{M}\in \Cb R$, then $\cpx{M} $ belongs to the smallest full
triangulated subcategory of $\D{R}$ containing $H^n(\cpx{M})$ with
all $n\in\mathbb{Z}$.

$(b)$ Any bounded-above complex over $R$ can be expressed as the
homotopy limit of its bounded ``quotient'' complexes, which are
obtained from canonical truncations.

$(c)$ Any bounded-below complex over $R$ can be expressed as the
homotopy colimit of its bounded ``sub'' complexes, which are
obtained from canonical truncations.

$(d)$ Any complex is generated by a bounded-above complex and a
bounded-below complex obtained by canonical truncations.

Recall that $\mathcal{Y}$ is a full triangulated subcategory of
$\D{R}$ closed under arbitrary direct sums and direct products in
$\D{R}$. Therefore it is closed under taking homotopy limits and
homotopy colimits in $\D{R}$. Now, by the fact
$\Img(\delta_*)\subseteq \mathcal{Y}$ and the above statements
$(a)$-$(d)$, we can show that $\mathcal{Y}'\subseteq \mathcal{Y}$.
Thus $\mathcal{Y}=\mathcal{Y'}.$

Next, we shall show that $H^n(i^*(R))=0$ for all $n\geq 1$. The idea
of the proof given here is essentially taken from \cite[Lemma
6.4]{nr1}.

On the one hand, from the adjoint pair $(i^*, i_*)$, we can obtain a
triangle in $\D{R}$:
$$
\cpx{X}\lra R \lraf{\eta_R} i^*(R)\lra \cpx{X}[1].
$$
It is cleat that the unit $\eta_R$ induces an isomorphism $\Hom_{\D
R}(i^*(R), \cpx{Y}[n])\simeq \Hom_{\D R}(R, \cpx{Y}[n])$ for each
$\cpx{Y}\in\mathcal{Y}$ and $n\in\mathbb{Z}$. This implies that
$\Hom_{\D R}(\cpx{X}, \cpx{Y}[n])=0$ for $\cpx{Y}\in \mathcal{Y}$
and $n\in \mathbb{Z}$.

On the other hand, by the canonical truncation at degree $0$, we
obtain a distinguished triangle of the following form in $\D R$:
$$
 i^*(R)^{\leq 0} \lraf{\alpha}  i^*(R) \lraf{\beta} i^*(R)^{\geq 1} \lra i^*(R)^{\leq 0}[1]
$$
such that $H^{s}\big(\,i^*(R)^{\leq 0}\big )\simeq \left\{\begin{array}{ll} 0   & \;\mbox{if}\; s\geq 1,\\
 H^s(i^*(R)) & \;\mbox{if}\;s\leq 0,\end{array} \right.$ and
$H^{t}\big(\,i^*(R)^{\geq 1}\big )\simeq \left\{\begin{array}{ll} 0            & \;\mbox{if}\; t\leq 0,\\
 H^t(i^*(R)) & \;\mbox{if}\;t\geq 1.\end{array} \right.$

\medskip
\noindent It follows that $\eta_R\beta=0$ and that there exists a
homomorphism $\gamma: R\to i^*(R)^{\leq 0}$ such that $\gamma
\,\alpha =\eta_R$. Since $i^*(R)\in\mathcal{Y}=\mathcal{Y}'$, we
know that $i^*(R)^{\leq 0}\in\mathcal{Y}$ and $\Hom_{\D R}(\cpx{X},
i^*(R)^{\leq 0})=0$. Consequently, there exists a homomorphism
$\theta: i^*(R)\to i^*(R)^{\leq 0}$ such that $\gamma=\eta_R\,
\theta$. So, we have the following diagram in $\D{R}$:
$$\xymatrix{& & i^*(R)^{\leq 0}\ar[d]_-{\alpha}&\\
\cpx{X}\ar[r] &R\ar@{-->}[ur]^-{\gamma} \ar@{-->}[rd]^{0}
\ar[r]^-{\eta_R}
& i^*(R)\ar@/_1.4pc/[u]_-{\theta}\ar[r]\ar[d]_-{\beta}& \cpx{X}[1] \\
&    &    i^*(R)^{\geq 1}\ar[d]& \\
&    &    i^*(R)^{\leq 0}[1] &}
$$
Further, one can check that $\eta_R \theta\, \alpha=\gamma
\,\alpha=\eta_R$. Since $\eta_R :R\to i_*i^*(R)=i^*(R)$ is a unit
morphism, we infer that $\theta \alpha= Id_{i^*(R)}$, and so
$$H^n(\theta \alpha)=H^n(\theta) H^n(\alpha)=Id{_{H^n(i^*(R))}}\; \mbox{for any }\; n\in\mathbb{Z}.$$
This means that $H^n(\theta): H^n(i^*(R))\ra
H^{n}\big(\,i^*(R)^{\leq 0}\big )$ is injective. Observe that
$H^{n}\big(\,i^*(R)^{\leq 0}\big )=0$ for $n\geq 1$. Hence
$H^n(i^*(R))=0$ for $n\geq 1$.

Finally, we shall prove that $\delta: R\to \Lambda$ is a ring
epimorphism.

Clearly, the $\delta$ is a ring epimorphism if and only if for every
$\Lambda$-module $M$, the induced map $\Hom_R(\delta, M):
\Hom_R(\Lambda, M)\lra \Hom_R(R, M)$ is an isomorphism. Observe that
$\Hom_R(\delta,M)$ is always surjective. To see that this map is
also injective, we shall use the commutative diagram in $(2)$ and
show that the induced map
$$\Hom_R\big(H^0(\eta_R), M\big): \Hom_R\big( H^0(i^*(R)),
M\big)\lra \Hom_R(R, M)$$ is injective. That is, we have to prove
that if $f_i: H^0(i^*(R))\to M$, with $i=1,2$, are two homomorphisms
in $R\Modcat$ such that $H^0(\eta_R) f_1= H^0(\eta_R) f_2$, then
$f_1=f_2$.

Now, we describe the map $H^0(\eta_R)$. Recall that $H^n(i^*(R))=0$
for all $n\geq 1$. Without loss of generality, we may assume that
$i^*(R)$ is of the following form $\big($up to isomorphism in $\D
R\big)$:
$$\cdots \lra
V^{-n}\lraf{d^{-n}} V^{-n+1}\lra\cdots \lra V^{-1}\lraf{d^{-1}} V^0
\lra 0\lra \cdots$$ From the canonical truncation, we can obtain the
following distinguished triangle in $\D{R}$:
$$
 \cpx{V}\,^{\leq -1}\lra i^*(R)\lraf{\pi} H^0(i^*(R))\lra \cpx{V}\,^{\leq -1}[1]
$$
where $\cpx{V}\,^{\leq -1}$ is of the form:
$$\cdots \lra V^{-n}\lra
V^{-n+1}\lra\cdots \lra V^{-2}\lra \Ker({d^{-1}}) \lra 0\lra
\cdots$$and  $\pi$ is the chain map induced by the canonical
surjection $V^0\to H^0(i^*(R))=\Coker(d^{-1})$. Applying $H^0(-)=
\Hom_{\D R}(R,-)$ to the above triangle, we see that
$H^0(\eta_R)=\eta_R\,\pi$ in $\D{R}$ and that $H^0(\pi)$ is an
isomorphism of $R$-modules.

Suppose that  $H^0(\eta_R) f_1= H^0(\eta_R) f_2: R\to M$ with $f_i:
H^0(i^*(R))\to M$ for $i=1,2$. Then $\eta_R\pi f_1=\eta_R\pi f_2$.
From the proof of $(2)$, we have $\Img(\delta_*)\subseteq
\mathcal{Y}$. Thus $_RM\in \mathcal{Y}$ since $M$ is an
$\Lambda$-module.  Note that the unit $\eta_R :R\to
i_*i^*(R)=i^*(R)$ induces an isomorphism $\Hom_{\D R}(i^*(R),
M)\simeq \Hom_{\D R}(R, M)$. Thus $\pi\, f_1=\pi\,f_2$ and
$H^0(\pi)f_1=H^0(\pi)f_2$. It follows from the isomorphism of
$H^0(\pi)$ that $f_1=f_2$. This means that $\Hom_R\big(H^0(\eta_R),
M\big)$ is injective, and thus $\delta$ is a ring epimorphism. This
finishes the proof of $(3)$. $\square$

\medskip
In the following, we shall systematically discuss when bireflective
subcategories of derived categories are homological. Note that some
partial answers have been given in the literature, for example, see
\cite[Theorem 0.7 and Proposition 5.6]{nr1}, \cite[Proposition
1.7]{HKL} and \cite[Proposition 3.6]{CX1}. Let us first mention the
following criterions.

\begin{Lem} \label{hep}
Let $\mathcal{Y}$ be a bireflective subcategory of $\D{R}$, and let
$i^*: \D{R}\to \mathcal{Y}$ be a left adjoint of the inclusion
$\mathcal{Y}\hookrightarrow \D{R}$. Then the following are
equivalent:

$(1)$ $\mathcal{Y}$ is homological.

$(2)$ $H^m(i^*(R))=0$ for any $m\neq 0$.

$(3)$  $H^0(i^*(R))\in\mathcal{Y}$ and  $H^m(i^*(R))=0$ for any
$m<0$.

$(4)$ $H^0(i^*(R))\in\mathcal{Y},$ and the associated ring
homomorphism $\delta: R\to \End_{\D R}(i^*(R))$ is a homological
ring epimorphism.

$(5)$ There exists a ring epimorphism $\lambda:R\to S$ such that
$_RS\in\mathcal{Y}$ and $i^*(R)$ is isomorphic in $\D R$ to a
complex $\cpx{Z}:=(Z^n)_{n\in\mathbb{Z}}$ with $Z^i\in S\Modcat$ for
$i=0, 1$.

$(6)$ $\mathscr{E}:=\mathcal{Y}\cap R\Modcat$ is an abelian
subcategory of $R\Modcat$ such that $i^*(R)$ is isomorphic in
$\D{R}$ to a complex $\cpx{Z}:=(Z^n)_{n\in\mathbb{Z}}$  with $Z^i\in
\mathscr{E}$ for $i=0, 1$.

In particular, if one of the above conditions is fulfilled, then
$\mathcal{Y}$ can be realized as the derived category of $\End_{\D
R}(i^*(R))$ via $\delta$.
\end{Lem}

{\it Proof.} It follows from the proof of \cite[Proposition
1.7]{HKL} that $(1)$ and $(2)$ are equivalent, and that $(2)$
implies both $(3)$ and $(4)$. By Lemma \ref{ringhom} (3), we know
that $(3)$ implies $(2)$.

Now, we show that $(4)$ implies $(1)$. In fact, since
$H^0(i^*(R))\in \mathcal{Y}$, it follows from Lemma \ref{ringhom}
(3) that
$$\mathcal{Y}=\{\cpx{Y}\in\D{R}\mid H^m(\cpx{Y})\in \Img(\delta_*)
\mbox{\;for\, all \,} m\in\mathbb{Z}\},$$ where $\delta: R\
\to\Lambda:= \End_{\D R}(i^*(R))$ is the associated ring
homomorphism. By assumption, $\delta$ is a homological ring
epimorphism, and therefore the derived functor
$D(\delta_*):\D{\Lambda}\to\D{R}$ is fully faithful. Furthermore, we
know from \cite[Lemma 4.6]{HKL} that
$$
\Img\big(D(\delta_*)\big)=\{\cpx{Y}\in\D{R}\mid H^m(\cpx{Y})\in
\Img(\delta_*) \mbox{\;for\, all \,} m\in\mathbb{Z}\}.
$$
Thus $\mathcal{Y}=\Img\big(D(\delta_*)\big)\subseteq \D{R}$, that
is, $\mathcal{Y}$ is homological by definition. Hence $(4)$ implies
$(1)$.

Consequently, we have proved that $(1)$-$(4)$ in Lemma \ref{hep} are
equivalent.

Note that $(5)$ and $(6)$ are equivalent because $\mathscr{E}$ is an
abelian subcategory of $R\Modcat$ if and only if there is a ring
epimorphism $\lambda:R\to S$ such that $\mathscr{E}=\Img(\lambda_*)$
(see \cite[Lemma 2.1]{CX1}).

In the following, we shall prove that $(1)$ implies $(5)$ and that
$(5)$ implies $(2)$.

Suppose that $\mathcal{Y}$ is homological, that is, there exists a
homological ring epimorphism $\lambda: R\to S$ such that the functor
$D(\lambda_*): \D{S}\to \D{R}$ induces a triangle equivalence from
$\D{S}$ to $\mathcal{Y}$. Thus $\mathcal{Y}=\Img(D(\lambda_*))$.
Since $i^*(R)\in\mathcal{Y}$, we have $i^*(R)\in\Img(D(\lambda_*))$.
It follows that there exists a complex
$\cpx{Z}:=(Z^n)_{n\in\mathbb{Z}}\in\C{S}$ such that
$i^*(R)\simeq\cpx{Z}$ in $\D{R}$. This shows $(5)$.

It remains to show that $(5)$ implies $(2)$. The idea of the
following proof arises from the proof of \cite[Proposition
3.6]{CX1}.

Let $\lambda:R\to S$ be a ring epimorphism satisfying the
assumptions in (5). We may identify $\Img(\lambda_*)$ with
$S\Modcat$ since $\lambda_*:S\Modcat\to R\Modcat$ is fully faithful.
Let $\cpx{Z}$ be a complex in $\C{R}$ such that $\cpx{Z}\simeq
i^*(R)$ in $\D{R}$. We may assume that $\cpx{Z}:=(Z^n,
d^n)_{n\in\mathbb{Z}}$ such that $Z^n\in S\Modcat$ for $n=0, 1$, and
define $\varphi=\Hom_{\D {R}}(\lambda, \cpx{Z}):\Hom_{\D R}(S,
\cpx{Z})\lra \Hom_{\D R}(R,\cpx{Z})$. We claim that the map
$\varphi$ is surjective.

In fact, there is a commutative diagram:
$$\xymatrix{\Hom{_{\K{R}}}(S, \cpx{Z})\ar[r]^-{q_1}\ar[d]^-{\varphi'}&
\Hom_{\D{R}}(S, \cpx{Z})\ar[d]^-{\varphi}\\
\Hom_{\K{R}}(R, \cpx{Z})\ar[r]^-{q_2}&\Hom_{\D{R}}(R, \cpx{Z}),}$$
where $\varphi'=\Hom_{\K{R}}(\lambda, \cpx{Z})$, and where $q_1$ and
$q_2$ are induced by the localization functor $q:\K{R}\to\D{R}$.
Clearly, the $q_2$ is a bijection. To prove that $\varphi $ is
surjective, it is sufficient to show that $\varphi'$ is surjective.

Let $\bar{\cpx{f}}:=\overline{(f^i)}\in \Hom_{\K{R}}(R,\cpx{Z})$
with $(f^i)_{i\in {\mathbb Z}}$ a chain map from $R$ to $\cpx{Z}$.
Then $f^i=0$ for any $i\neq 0$ and $f^0 d^0=0$. Since $Z^0$ is an
$S$-module, we can define $g: S\to Z^0$ by $s\mapsto s\,(1)f^0$ for
$s\in S$. One can check that $g$ is a homomorphism of $R$-modules
with $f^0=\lambda g$, as is shown in the following visual diagram:
$$\xymatrix{
            &             &R\ar[r]^-{\lambda}\ar[d]_-{f^0}& S\ar@{-->}[ld]_-{g}&&&\\
\cdots\ar[r]&Z^{-1}\ar[r]^-{d^{-1}} &Z^0\ar[r]^-{d^0} &
Z^1\ar[r]^-{d^1} & Z^2\ar[r]&\cdots}
$$
Since $\lambda:R\to S$ is a ring epimorphism and since $Z^1$ is an
$S$-module, the induced map $\Hom_R(\lambda, Z^1): \Hom_R(S, Z^1)\to
\Hom_R(R, Z^1)$ is a bijection. Thus, from this bijection together
with $\lambda g d^0=f^0d^0=0$, it follows that $gd^0=0$. Now, we can
define a morphism $\bar{\cpx{g}}:=\overline{(g^i)}\in
\Hom_{\K{R}}(S,\cpx{Z})$, where $(g^i)_{i\in {\mathbb Z}}$ is the
chain map with $g^0=g$ and $g^i=0$ for any $i\neq 0$. Thus
$\bar{\cpx{f}}=\lambda\bar{\cpx{g}}$. This shows that $\varphi'$ is
surjective. Consequently, the map $\varphi$ is surjective, and the
induced map
$$\Hom_{\D {R}}(\lambda, i^*(R)):\Hom_{\D R}(S, i^*(R))\to
\Hom_{\D R}(R, i^*(R))$$ is surjective since $\cpx{Z}\simeq i^*(R)$
in $\D{R}$.

Finally, we shall prove that $i^*(R)\simeq S$ in $\D{R}$. In
particular, this will give rise to $H^m(i^*(R))\simeq H^m(S)=0$ for
any $m\neq 0$, and therefore show $(2)$. So, it suffices to prove
that $i^*(R)\simeq S$ in $\D{R}$.

Indeed, let $i_*: \mathcal{Y}\to \D{R}$ be the inclusion, and let
$\eta_R :R\to i_*i^*(R)$ be the unit with respect to the adjoint
pair $(i^*, i_*)$. Clearly, $i^*(R)= i_*i^*(R)$ in $\D{R}$. Since we
have proved that $\Hom_{\D {R}}(\lambda, i^*(R))$ is surjective,
there exists a homomorphism $v: S\to i_*i^*(R)$ in $\D{R}$ such that
$\eta_R=\lambda\,v$. Furthermore, since $_RS$ belongs to
$\mathcal{Y}$ by assumption, we see that $\Hom_{\D R}(\eta_R,S):
\Hom_{\D R}(i^*(R),S)\ra \Hom_{\D R}(R,S)$ is an isomorphism. Thus
there exists a homomorphism $u: i_*i^*(R)\to S$ in $\D{R}$ such that
$\lambda=\eta_R\,u$. This yields the following commutative diagram
in $\D{R}$:
$$
\xymatrix{
R\ar@{=}[r]\ar[d]_-{\eta_R}&R\ar@{=}[r]\ar[d]_-{\lambda}&R \ar[d]^-{\eta_R}\\
i_*i^*(R)\ar@{-->}[r]^-{u}&S \ar@{-->}[r]^-{v}&\;i_*i^*(R), }
$$
which shows that $\eta_R=\eta_R uv$ and $\lambda=\lambda v u$. Since
$\Hom_{\D R}(\eta_R, i^*(R))$: $\Hom_{\D R}(i^*(R),i^*(R))\ra
\Hom_{\D R}(R,i_*i^*(R))$ is an isomorphism, we clearly have
$uv=1{_{i_*i^*(R)}}$. Note that $\Hom_R(\lambda,S):\Hom_R(S,
S)\to\Hom_R(R, S)$ is bijective since $\lambda:R\to S$ is a ring
epimorphism. It follows from $\lambda=\lambda vu$ that $vu=1_{S}$.
Thus the map $u$ is an isomorphism in $\D{R}$, and
$i^*(R)=i_*i^*(R)\simeq S$ in $\D{R}$. This shows that $(5)$ implies
$(2)$.

Hence all the statements in Lemma \ref{hep} are equivalent. This
finishes the proof. $\square$.

\medskip
Now, we mention a special bireflective subcategory of $\D{R}$, which
is constructed from complexes of finitely generated projective
$R$-modules. For the proof, we refer to \cite[Chapter III, Theorem
2.3; Chapter IV, Proposition 1.1]{BI}. See also \cite[Lemma
2.8]{CX1}.

\begin{Lem}\label{ref}
Let $\Sigma$ be a set of complexes in $\Cb{\pmodcat{R}}$. Define
$\mathcal{Y}:=\Ker\big(\,\Hom_{\D R}({\rm{Tria}}{(\Sigma)},
-)\big)$. Then $\mathcal{Y}$ is bireflective and equal to the full
subcategory of $\D {R}$ consisting of complexes $\cpx{Y}$ in $\D{R}$
such that $\Hom_{\D R}(\cpx{P}, \,\cpx{Y}[n])$ = $0$ for every
$\cpx{P}\in\Sigma$ and $n\in\mathbb{Z}$.
\end{Lem}

To develop properties of the bireflective subcategories of $\D{R}$
in Lemma \ref{ref}, we shall define the so-called generalized
localizations, which is motivated by a discussion with Silvana
Bazzoni in 2012. In fact, this notion was first discussed in
\cite{Kr} under the name ``homological localizations" for a set of
complexes in $\Cb{\pmodcat{R}}$, and is related to both the
telescope conjecture and algebraic $K$-theory. The reason for not
choosing the adjective word ``homological" in this note is that we
have reserved this word for ring epimorphisms.

\begin{Def} \label{genloc}{\rm
Let $R$ be a ring, and let $\Sigma$ be a set of complexes in $\C R$.
A homomorphism $\lambda_{\Sigma}: R\to R_{\Sigma}$ of rings is
called a \emph{generalized localization} of $R$ at $\Sigma$ provided
that

$(1)$ $\lambda_{\Sigma}$ is $\Sigma$-exact, that is, if $\cpx{P}$
belongs to $\Sigma$, then $R_{\Sigma}\otimes_R\cpx{P}$ is exact as a
complex over $R_{\Sigma}$, and

$(2)$ $\lambda_{\Sigma}$ is universally $\Sigma$-exact, that is, if
$S$ is a ring  together with a $\Sigma$-exact homomorphism
$\varphi:R\to S$, then there exists a unique ring homomorphism
$\psi:R_{\Sigma}\to S$ such that $\varphi=\lambda_{\Sigma}\psi$. }
\end{Def}

If $\Sigma$ consists only of two-term complexes in
$\Cb{\pmodcat{R}}$, then the generalized localization of $R$ at
$\Sigma$ is the \emph{universal localization} of $R$ at $\Sigma$ in
the sense of Cohn (see \cite{cohenbook1}). It was proved in
\cite{cohenbook1} that universal localizations always exist.
However, generalized localizations may not exist in general. For a
counterexample, we refer the reader to \cite[Example 15.2]{Kr}.

We remark that, in Definition \ref{genloc} (1), if $\Sigma$ consists
of complexes in $ \Cb{\pmodcat{R}}$, then, for each
$\cpx{P}:=(P^i)_{i\in\mathbb{Z}}\in\Sigma$, the complex
$R_{\Sigma}\otimes_R\cpx{P}$ is actually split exact as a complex
over $R_{\Sigma}$ since $R_\Sigma\otimes_RP^i$ is a projective
$R_\Sigma$-module for each $i$. Further, by Definition \ref{genloc}
(2), if $\lambda_i: R\to R_i$ is a generalized localization of $R$
at $\Sigma$ for $i=1,2$, then $\lambda_1$ and $\lambda_2$ are
equivalent, that is, there exists a ring isomorphism $\rho: R_1\to
R_2$ such that $\lambda_2=\lambda_1\rho$.

Suppose that $\mathcal {U}$ is a set of $R$-modules each of which
possesses a finitely generated projective resolution of finite
length. For each $U\in\mathcal{U}$, we choose such a projective
resolution $_pU$ of finite length, and set $\Sigma:=\{_pU\mid
U\in\mathcal{U}\}\subseteq \Cb{\pmodcat{R}}$, and let
$R_\mathcal{U}$ be the generalized localization of $R$ at $\Sigma$.
If $_PU '$ is another choice of finitely generated projective
resolution of finite length for $U$ , then the generalized
localization of $R$ at $\Sigma':=\{_pU '\mid U\in {\mathcal U}\}$ is
isomorphic to $R_\mathcal{U}$, that is, $R_\mathcal{U}$ does not
depend on the choice of projective resolutions of $U$. Thus, we may
say that $R_{\mathcal U}$ is the \emph{generalized  localization of
$R$ at $\mathcal U$.}

Generalized  localizations have the following simple properties
(compare with \cite[Theorem 3.1 and Lemma 3.2]{CX1}).

\begin{Lem}\label{genpro}
Let $R$ be a ring and let $\Sigma$ be a set of complexes in
$\Cb{\pmodcat{R}}$. Assume that the generalized localization
$\lambda_{\Sigma}: R\to R_{\Sigma}$ of $R$ at $\Sigma$ exists. Then
the following hold.

$(1)$ For any homomorphism $\varphi:R_\Sigma\to S$ of rings, the
ring homomorphism $\lambda_{\Sigma}\varphi: R\to S$ is
$\Sigma$-exact.

$(2)$ The ring homomorphism $\lambda_{\Sigma}$ is a ring
epimorphism.

$(3)$ Define $\Sigma^*:=\{\Hom_R(\cpx{P}, R)\mid
\cpx{P}\in\Sigma\}.$ Then $\lambda_{\Sigma}$ is also the generalized
localization of $R$ at the set $\Sigma^*$. In particular,
$R_{\Sigma^*}\simeq R_{\Sigma}$ as rings.
\end{Lem}

{\it Proof.} $(1)$ For each $\cpx{P}\in\Sigma$, we have the
following isomorphisms of complexes of $S$-modules:
$$ S\otimes_R\cpx{P}\simeq
(S\otimes_{R_\Sigma}R_\Sigma)\otimes_R\cpx{P}\simeq
S\otimes_{R_\Sigma}(R_\Sigma\otimes_R\cpx{P}).$$ Since
$R_\Sigma\otimes_R\cpx{P}$ is split exact in $\C{R_\Sigma}$, we see
that $S\otimes_R\cpx{P}$ is  also split exact in $\C{S}$. This means
that the ring homomorphism $\lambda_{\Sigma}\varphi$ is
$\Sigma$-exact.

$(2)$ Assume that $\varphi_i: R_\Sigma\to S$ is a ring homomorphism
for $i=1,2,$ such that
$\lambda_{\Sigma}\varphi_1=\lambda_{\Sigma}\varphi_2$. It follows
from $(1)$ that $\lambda_{\Sigma}\varphi_i$ is $\Sigma$-exact. By
the property $(2)$ in Definition \ref{genloc}, we obtain
$\varphi_1=\varphi_2$. This implies that $\lambda_{\Sigma}$ is a
ring epimorphism.

$(3)$ Note that $\cpx{P}$ is in $\Cb{\pmodcat{R}}$. It follows from
Lemma \ref{complex} that, for any homomorphism $R\to S$ of rings,
there are the following isomorphisms of complexes:
$$\Hom_R(\cpx{P}, R)\otimes_RS\simeq
\Hom_R(\cpx{P}, S)\simeq\Hom_R(\cpx{P}, \Hom_S({_S}S{_R},\,S))\simeq
\Hom_S(S\otimes_R\cpx{P}, S).$$ This implies that the complex
$\Hom_R(\cpx{P}, R)\otimes_RS$ is (split) exact in $\C{S\opp}$ if
and only if so is the complex $S\otimes_R\cpx{P}$ in $\C{S}$. Now,
$(3)$ follows immediately from the definition of generalized
localizations. $\square$

\medskip
In the following, we shall establish a relation between bireflective
subcategories of $\D{R}$ and generalized localizations. In
particular, the statements $(3)$ and $(4)$ in Lemma \ref{rg} below
will be useful for discussions in the next section and the proof of
Theorem \ref{main-result}.

\begin{Lem}\label{rg}
Let $\Sigma$ be a set of complexes in $\Cb{\pmodcat{R}}$, and let
$j_!:{\rm{Tria}}{(\Sigma)}\to \D{R}$ be the inclusion. Define
$\mathcal{Y}:=\Ker\big(\,\Hom_{\D R}({\rm{Tria}}{(\Sigma)},
-)\big)$. Then the following are true.

$(1)$ There exists a recollement of triangulated categories:
$$\xymatrix@C=1.2cm{\mathcal{Y}\ar[r]^-{{i_*}}
&\D{R}\ar[r]\ar@/^1.2pc/[l]\ar_-{i^*}@/_1.2pc/[l]
&{\rm{Tria}}{(\Sigma)} \ar@/^1.2pc/[l]\ar@/_1.2pc/[l]_{j_!\;}}$$

\smallskip
\noindent where $(i^*, i_*)$ is a pair of adjoint functors with
$i_*$ the inclusion.

\smallskip
$(2)$ The associated ring homomorphism $\delta: R\to
\Lambda:=\End_{\D{R}}(i^*(R))$ induced by $i^*$ admits the following
property: For any $\Sigma$-exact ring homomorphism $\varphi: R\to
S$, there exists a ring homomorphism $\psi: \Lambda \to S$ such that
$\varphi=\delta\psi$.

$(3)$ If $H^0(i^*(R))\in\mathcal{Y}$, then $\delta$ is a generalized
localization of $R$ at $\Sigma$. In particular, if the subcategory
$\mathcal{Y}$ of $\D{R}$ is homological, then $\delta$ is a
generalized localization of $R$ at $\Sigma$.

$(4)$ Define $\Sigma^*:=\{\Hom_R(\cpx{P}, R)\in
\Cb{\pmodcat{R\opp}}\mid \cpx{P}\in\Sigma\}$ and
$\mathcal{Y}':=\Ker\big(\,\Hom_{\D {R\opp}}({\rm{Tria}}{(\Sigma^*)},
-)\big).$ Then $\mathcal{Y}$ is homological in $\D{R}$ if and only
if so is $\mathcal{Y}'$ in $\D{{R\opp}}$.
\end{Lem}

{\it Proof.} $(1)$ can be concluded from  \cite[Lemma 2.6 and Lemma
2.8]{CX1}.

$(2)$ The proof here is motivated by \cite[Lemma 7.3]{nr1}. Let
$\varphi: R\to S$ be a $\Sigma$-exact ring homomorphism. Since
$S\otimes_R\cpx{P}$ is exact in $\C{S}$ for $\cpx{P}\in\Sigma$, we
have $S\otimesL_R\cpx{P}=S\otimes_R\cpx{P}\simeq 0$ in $\D{S}$.
Further, the functor $S\otimesL_R-: \D{R}\to\D{S}$ commutes with
arbitrary direct sums, so $S\otimesL_R\cpx{X}\simeq 0$ for each
$\cpx{X}\in {\rm Tria}(\Sigma)$.

Let $\D{R}/{\rm Tria}(\Sigma)$ denote the Verdier quotient of
$\D{R}$ by the full triangulated subcategory ${\rm Tria}(\Sigma)$.
It follows from the recollement in $(1)$ that $i^*$ induces a
triangle equivalence:
$$\D{R}/{\rm Tria}(\Sigma)\lraf{\simeq}\mathcal{Y}.$$
Since $S\otimesL_R-$ sends ${\rm Tria}(\Sigma)$ to zero, there
exists a triangle functor $F:\mathcal{Y}\to \D{S}$ together with a
natural isomorphism of triangle functors:
$$\Phi:\;S\otimesL_R-\lraf{\simeq}
F\,i^*:\,\D{R}\lra \D{S}.$$ This clearly induces the following
canonical ring homomorphisms:
$$\Lambda:=\End_{\D{R}}(i^*(R)) \lraf{F}
\End_{\D{S}}\big(F(i^*(R))\big)\simeq \End_{\D S}(S\otimesL_RR)
\simeq \End_{\D S}(S)\simeq S
$$
where the first isomorphism is induced by the natural isomorphism
$\Phi_R: S\otimesL_RR\lra F(i^*(R))$ in $\D{S}$. Now, we define
$\psi: \Lambda \to S$ to be the composite of the above ring
homomorphisms. Then it is easy to check that $\varphi=\delta\psi$.
Consequently, the $\delta$ has the property mentioned  in $(2)$.

$(3)$ Assume that $H^0(i^*(R))\in\mathcal{Y}$. By Lemma
\ref{ringhom} (3), the map $\delta$ is a ring epimorphism. Combining
this with $(2)$, we know that $\delta$ satisfies the condition $(2)$
in Definition \ref{genloc}. To see that $\delta $ is the generalized
localization of $R$ at $\Sigma$, we have to show that $\delta$
satisfies the condition $(1)$ in Definition \ref{genloc}, that is,
$\delta$ is $\Sigma$-exact.

In fact, by Lemma \ref{ringhom} (2), we have $\Lambda\simeq
H^0(i^*(R))$ as $R$-modules. This gives rise to
$_R\Lambda\in\mathcal{Y}$. Note that $\Hom_{\D R}(\cpx{X},
\cpx{Y})=0$ for $\cpx{X}\in {\rm{Tria}}{(\Sigma)}$ and
$\cpx{Y}\in\mathcal{Y}$. In particular, we have $\Hom_{\D
R}(\cpx{P}, \Lambda[n])=0$ for any $\cpx{P}\in\Sigma$ and
$n\in\mathbb{Z}$. It follows that $H^n(\Hom_R(\cpx{P},
\Lambda))\simeq \Hom_{\K R}(\cpx{P}, \Lambda[n])\simeq \Hom_{\D
R}(\cpx{P}, \Lambda[n])=0$, and therefore the complex
$\Hom_R(\cpx{P}, \Lambda)$ is exact. Since $\cpx{P}\in
\Cb{\pmodcat{R}}$, we have $\Hom_R(\cpx{P},
\Lambda)\in\Cb{\pmodcat{\Lambda\opp}}$. This implies that
$\Hom_R(\cpx{P}, \Lambda)$ is split exact, and therefore the complex
$\Hom_{\Lambda\opp}(\Hom_R(\cpx{P}, \Lambda), \Lambda)$ over
$\Lambda$ is split exact. Now, we claim that the latter complex is
isomorphic to the complex $\Lambda\otimes_R\cpx{P}$ in
$\C{\Lambda}$. Actually, this follows from the following general
fact in homological algebra:

For any finitely generated projective $R$-module $P$, there exists a
natural isomorphism of $\Lambda$-modules:
$$\Lambda\otimes_RP \lra \Hom_{\Lambda\opp}(\Hom_R(P, \Lambda),
\Lambda), \quad x\otimes p\mapsto [f\mapsto x\,(p)f]$$ for
$x\in\Lambda$, $p\in P$ and $f\in\Hom_R(P, \Lambda)$. Consequently,
the complex $\Lambda\otimes_R\cpx{P}$ is exact in $\C{\Lambda}$, and
thus $\delta$ is $\Sigma$-exact. Hence $\delta$ is a generalized
localization of $R$ at $\Sigma$.

Clearly, the second part of Lemma \ref{rg} (3) follows from the
equivalences of $(1)$ and $(4)$ in Lemma \ref{hep}.

$(4)$ We shall only prove the necessity of $(4)$ since the
sufficiency of $(4)$ can be proved similarly.

Suppose that $\mathcal{Y}$ is homological in $\D{R}$. It follows
from Lemma \ref{hep} (4) and  Lemma \ref{rg} (3) that the ring
homomorphism $\delta: R\to \Lambda$ is not only a homological ring
epimorphism, but also a generalized  localization of $R$ at
$\Sigma$. Moreover, by Lemma \ref{genpro} (3), the map $\delta$ is
also a generalized localization of $R$ at $\Sigma^*$.

Note that $\mathcal{Y}'$ is a bireflective subcategory of
$\D{R\opp}$ by Lemma \ref{ref}. Now, let ${\bf L}$ be a left adjoint
of the inclusion $\mathcal{Y}'\to \D{R\opp}$. To show that
$\mathcal{Y}'$ is homological in $\D{R\opp}$, we employ  the
equivalences of $(1)$ and $(4)$ in Lemma \ref{hep}, and prove that

$(a)$ $H^0({\bf L}(R))\in\mathcal{Y}'$ and

$(b)$ the ring homomorphism $\delta': R\lra \Lambda':=\End_{\D
{R\opp}}({\bf L}(R))$ induced by $\bf{L}$ is homological.

Clearly, under the assumption $(a)$, we see from $(3)$ that
$\delta'$ is a generalized  localization of $R$ at $\Sigma^*$. Since
$\delta$ is also a generalized  localization of $R$ at $\Sigma^*$,
there exists a ring isomorphism $\rho: \Lambda' \lra \Lambda$ such
that $\delta=\delta' \rho$. Note that $\delta$ is homological. It
follows that $\delta'$ is homological.

It remains to show $(a)$. In fact, since $H^0({\bf L}(R))\simeq
\Lambda'$ as right $R$-modules by Lemma \ref{ringhom} (2), it is
sufficient to prove that the right $R$-module $\Lambda'$ belongs to
$\mathcal{Y}'$. However, by $(1)$ and Lemma \ref{ref}, we have
$$\mathcal{Y}'=\{\cpx{Y}\in\D{R\opp} \mid \Hom_{\D
{R\opp}}\big(\Hom_R(\cpx{P}, R),\, \cpx{Y}[n]\big)=0\; \mbox{for}\;
\cpx{P}\in\Sigma\; \mbox{and}\; n\in\mathbb{Z}\},$$ and by the
isomorphism $\rho$ and $\delta=\delta' \rho$, we get $\Lambda'
\simeq \Lambda$ as right $R$-modules. Consequently, to show
$\Lambda'_{R}\in {\mathcal Y}'$, it is enough to show that
$\Lambda_R$ belongs to $\mathcal{Y}'$, that is, we have to prove
that $\Hom_{\D {R\opp}}\big(\Hom_R(\cpx{P}, R),\, \Lambda[n]
\big)=0$ for any $\cpx{P}\in\Sigma$ and $n\in\mathbb{Z}$.

Let $\cpx{P}\in \Sigma$, and set $\cpx{P}{^*}:=\Hom_R(\cpx{P}, R)$.
Since $\cpx{P}$ is a complex in $\Cb{\pmodcat{R}}$, we see from
Lemma \ref{complex} that $\Hom_{R\opp}(\cpx{P}{^*}, \Lambda)\simeq
\Lambda\otimes_R\cpx{P}$ as complexes in $\C \Lambda$, and therefore
there exist the following isomorphisms:
$$\Hom_{\D {R\opp}}\big(\cpx{P}{^*},\, \Lambda[n]
\big)\simeq \Hom_{\K {R\opp}}\big(\cpx{P}{^*},\, \Lambda[n]
\big)\simeq H^n(\Hom_{R\opp}(\cpx{P}{^*}, \Lambda))\simeq
H^n(\Lambda\otimes_R\cpx{P}).$$  Since $\delta: R\to \Lambda$ is a
generalized  localization of $R$ at $\Sigma$, the complex
$\Lambda\otimes_R\cpx{P}$ is exact in $\C{\Lambda}$, that is,
$H^n(\Lambda\otimes_R\cpx{P})=0$ for any $n\in\mathbb{Z}$. Thus
$\Hom_{\D {R\opp}}\big(\cpx{P}{^*},\, \Lambda[n] \big)=0$ for
$n\in\mathbb{Z}$. Thus $\Lambda_R\in \mathcal{Y}'$, and the proof of
the necessity of $(4)$ is completed. $\square$

\medskip
As an application of Lemma \ref{rg} (3), we have the following
result which says that generalized localizations can be constructed
from homological ring epimorphisms.
\begin{Koro}
Let $\lambda:R\to S$ be a homological ring epimorphism. Suppose that
$_RS$ has a finitely generated projective resolution of finite
length. Let $\cpx{P}$ be a complex in $\Cb{\pmodcat{R}}$, which is
isomorphic in $\D R$ to the mapping cone of $\lambda$. Then
$\lambda$ is a generalized localization of $R$ at $\cpx{P}$.
\end{Koro}

{\it Proof.} Since $\lambda$ is homological and $\cpx{P}$ is
isomorphic to the mapping cone of $\lambda$ in $\D{R}$, it follows
from \cite[Section 4]{NS} that there is a recollement of
triangulated categories:
$$
\xymatrix@C=1.2cm{\D{S}\ar[r]^-{D(\lambda_*)}
&\D{R}\ar[r]\ar@/^1.2pc/[l]\ar_-{S\otimesL_R-}@/_1.2pc/[l]
&{\rm{Tria}}{(\cpx{P})} \ar@/^1.2pc/[l]\ar@/_1.2pc/[l]_{j_!\;}}$$

\smallskip
\noindent where $j_!$ is the inclusion. This shows that
$\mathcal{Y}:=\Ker\big(\Hom_{\D R}({\rm Tria}(\cpx{P}),-)\big)$ is
equivalent to $\D S$. Thus $\mathcal Y$ is homological. Note that
$S\otimesL_RR\simeq S$ and $\End_R(_RS)\simeq S$. By Lemma \ref{rg}
(3), we know that $\lambda$ is a generalized localization of $R$ at
$\cpx{P}$. $\square$

\section{Ringel modules}
\label{sect4}

This section is devoted to preparations for proofs of our main
results in this paper. First, we introduce a special class of
modules, called Ringel modules, which can be constructed from both
good tilting  and cotilting modules, and then discuss certain
bireflective subcategories (of derived module categories) arising
from Ringel modules. Finally, we shall describe when these
subcategories are homological. In particular, we shall establish a
key proposition, Proposition \ref{realization}, which will be
applied in later sections.

\smallskip
Throughout this section, let $R$ be an arbitrary ring, $M$ an
$R$-module and $S$ the endomorphism ring of $_RM$. Then $M$ becomes
naturally an $R$-$S$-bimodule. Further, let $n$ be an arbitrary but
fixed natural number.

\begin{Def}\label{rm}
{\rm The $R$-module $M$ is called an $n$-\emph{Ringel module}
provided that the following three conditions are fulfilled:

$(R1)$ there exists an exact sequence
$$
0\lra P_n \lra \cdots \lra P_1\lra P_0\lra M\lra 0$$ of $R$-modules
such that $P_i\in \add(_RR)$ for all $0\leq i\leq n$,

$(R2)$ $\Ext^j_R(M, M)=0$ for all $j\geq 1$, and

$(R3)$ there exists an exact sequence
$$
0\lra {}_RR \lra M_0 \lraf{\nu} M_1\lra \cdots \lra M_n\lra 0$$ of
$R$-modules such that $M_i\in\Prod(_RM)$ for all $0\leq i\leq n$.

\smallskip

An $n$-Ringel $R$-module $M$ is said to be \emph{perfect} if the
ring $S$ is right noetherian; and \emph{good} if

$(R4)$ the right $S$-module $M$ is strongly $S$-Mittag-Leffler (see
Definition \ref{ML}). }\end{Def}

Classical tilting modules are good Ringel modules. Conversely, for a
Ringel module $M$, if each $M_i$ in $(R3)$ is isomorphic to a direct
summand of finite direct products of copies of $M$, then $M$ becomes
a classical tilting module (see Introduction).

If a Ringel $R$-module $M$ has the property $\Prod(_RM) = \Add(_RM)$
(for example, $M_S$ is of finite length), then $_RM$ is a tilting
module. In this case, $_RM$ is even classical (see Corollary
\ref{2.9}).

Moreover, if the ring $S$ is right noetherian (see the statements
following Definition \ref{ML}), then any right $S$-module is
$S$-Mittag-Leffler. Thus each perfect Ringel $R$-module must be
good.

It is worth noting that good tilting (or cotilting) modules may not
be Ringel modules because it may not be finitely generated. For
example, the infinitely generated $\mathbb{Z}$-module
$\mathbb{Q}\oplus \mathbb{Q}/\mathbb{Z}$ is a good tilting module,
but not a Ringel module. Clearly, the good $1$-cotilting
$\mathbb{Z}$-module $\Hom_{\mathbb{Z}}(\mathbb{Q}\oplus
\mathbb{Q}/\mathbb{Z},\mathbb{Q}/{\mathbb Z})$ is not a Ringel
module.

Assume that $_RM$ satisfies $(R1)$. Then $M$ is isomorphic in
$\D{R}$ to the following complex of finitely generated projective
$R$-modules:
$$
\cdots \lra 0\lra P_n \lra \cdots \lra P_1\lra P_0\lra 0\lra\cdots
$$
It follows from Lemma \ref{ref} that
$\mathcal{Y}:=\{\cpx{Y}\in\D{R}\mid \Hom_{\D R}(M,\, \cpx{Y}[m])=0
\mbox{\;for\, all \,} m\in\mathbb{Z}\}$ is a bireflective
subcategory of $\D{R}$.

Now, assume that $M$ satisfies both $(R1)$ and $(R2)$. Then the
functors
$${\bf G}:={}_RM\otimesL_S-: \;\D{S}\lra\D{R}\quad \mbox{and}\quad  {\bf
H}:=\rHom_R(M,-):\; \D{R}\lra\D{S}$$ induce a triangle equivalence:
$\D{S}\lraf{\simeq}{\rm{Tria}}{(_RM)}$ (see \cite[Chapter 5,
Corollary 8.4, Theorem 8.5]{HHK}). Moreover, $\mathcal{Y}=\Ker({\bf
H })$ since $H^m\big(\rHom_R(M, \cpx{Y})\big)\simeq \Hom_{\D R}(M,
\cpx{Y}[m])$ for each $\cpx{Y}\in\D{R}$ and $m\in\mathbb{Z}$.

Thus, by Lemma \ref{rg} (1) and (3) as well as Lemma \ref{hep}, we
have the following useful result for constructing recollements of
derived module categories.

\begin{Lem}\label{rt}
Suppose that the $R$-module $M$ satisfies  $(R1)$ and $(R2)$. Then
there exists a recollement of triangulated categories:
$$(\ast)\quad\;
\xymatrix@C=1.2cm{\mathcal{Y}\ar[r]^-{{i_*}} &\D{R}\ar[r]^-{{\bf
H}}\ar@/^1.2pc/[l]\ar_-{i^*}@/_1.2pc/[l] &\D{S}
\ar@/^1.2pc/[l]\ar@/_1.2pc/[l]_{{\bf G}} }$$

\medskip
\noindent where $(i^*, i_*)$ is a pair of adjoint functors with
$i_*$ the inclusion.

If, in addition, the category $\mathcal{Y}$ is homological in
$\D{R}$, then the generalized  localization $\lambda: R\to R_M$ of
$R$ at $M$ exists and is homological, which induces a recollement of
derived module categories:
$$\xymatrix@C=1.2cm{\D{R_M}\ar[r]^-{D(\lambda_*)}
&\D{R}\ar[r]^-{{\bf H}}\ar@/^1.2pc/[l]\ar@/_1.2pc/[l] &\D{S}
\ar@/^1.2pc/[l]\ar@/_1.2pc/[l]_{{\bf G}}}$$
\end{Lem}

\medskip
In the following, we shall consider when the category $\mathcal{Y}$
is homological. In general, this category is not homological since
the category
$$\mathscr{E}:=\mathcal{Y}\cap R\Modcat=\{Y\in
R\Modcat \mid \Ext^m_R(M, Y)=0\mbox{\;for\, all \,} m\geq0\}$$ may
not be an abelian subcategory of $R\Modcat$. So, we need to impose
some additional conditions on the module $M$.

By Lemma \ref{hep}, whether $\mathcal{Y}$ is homological is
completely determined by the cohomology groups of $i_*i^*(R)$. So,
to calculate these cohomology groups efficiently, we shall
concentrate on good Ringel modules.

From now on, we assume that $_RM$ is a {\bf good} $n$-Ringel module,
and define $\cpx{M}$ to be the complex
$$\cdots \lra 0\lra M_0
\lraf{\nu} M_1\lra \cdots \lra M_n\lra 0\lra \cdots$$ arising from
$(R3)$ in Definition \ref{rm}, where $M_i$ is in degree $i$ for
$0\leq i\leq n$.

First of all, we establish the following result.

\begin{Lem}\label{prep}
The following statements are true.

\smallskip
$(1)$ For each $X\in\Prod(_RM)$, the evaluation map $\theta_{X}:
M\otimes_S\Hom_R(M, X)\lra X$ is injective and $\Coker(\theta
_X)\in\mathscr{E}$.

$(2)$
$$H^j(i_*i^*(R))\simeq
\left\{\begin{array}{ll} 0                        & \;\mbox{if}\; j<0,\\
H^{j+1}\big({_R}M\otimes_S\Hom_R(M,\,\cpx{M})\big)&
\;\mbox{if}\;j>0.\end{array} \right.$$

$(3)$ For $n=0$, the complex $i_*i^*(R)$ is isomorphic in $\D{R}$ to
the stalk complex $\Coker(\theta_{M_0})$. For $n\ge 1$, the complex
$i_*i^*(R)$ is isomorphic in $\D{R}$ to a complex of the form
$$
0\lra E^0 \lra E^1 \lra \cdots \lra E^{n-1}\lra 0
$$
with $E^m\in\mathscr{E}$ for  $0\leq m\leq n-1$.
\end{Lem}

\medskip
{\it Proof.} Recall that $M$ is an $R$-$S$-bimodule with
$S=\End_R(M)$. So we have a pair of adjoint functors:
$$_RM\otimes_S-: S\Modcat \lra R\Modcat \quad \mbox{and}\quad\Hom_R(M, -): R\Modcat \lra S\Modcat. $$
This can be naturally extended to a pair of adjoint triangle
functors between homotopy categories:
$$_RM\otimes_S-: \K{S}\lra \K{R} \quad \mbox{and}\quad\Hom_R(M, -): \K{R} \lra \K{S}. $$
By passing to derived categories, we obtain the derived functors
${\bf G}$ and ${\bf{H}}$, respectively. Further, let $$\theta:\,
M\otimes_S\Hom_R(M, -)\lra Id_{R\mbox{-}{\rm Mod}}\quad \mbox{and}
\quad \varepsilon: {\bf GH } \lra Id_{\D R}$$ be the counit
adjunctions with respect to $\big(M\otimes_S-, \Hom_R(M, -)\big)$
and  $({\bf G}, {\bf H})$, respectively.

Note that, for each $\cpx{X}\in\D{R}$, it follows from the
recollement $(\ast)$ in Lemma \ref{rt} that there exists a canonical
distinguished triangle in $\D{R}$:
$${\bf GH}(\cpx{X})\lraf{\varepsilon_{\cpx{X}}} \cpx{X} \lra i_*i^*(\cpx{X})\lra{\bf
GH}(\cpx{X})[1].$$

$(1)$ Let $X\in\Prod(_RM)$. To verify that $\theta_X$ is injective,
it is sufficient to show that
$$\theta_{M^I}: M\otimes_S\Hom_R(M, M^I)\lra  M^I$$
is injective for any nonempty set $I$.  Since $\Hom_R(M, M^I)\simeq
\Hom_R(M, M)^I$, the injection of $\theta_{M^I}$ is equivalent to
saying that the canonical map $\rho{_I}:\; M\otimes_SS^I\lra M^I$,
defined in Definition \ref{ML}, is injective. This holds exactly if
$M$ is $S$-Mittag-Leffler. However, the axiom $(R4)$ ensures that
$M$ is $S$-Mittag-Leffler. Thus $\theta_{X}: M\otimes_S\Hom_R(M,
X)\lra X$ is injective.

To prove $\Coker(\theta_X)\in\mathscr{E}:={\mathcal Y}\cap
R\Modcat$, we demonstrate that there is the following commutative
diagram in $\D{R}$:
$$(a)\quad
\xymatrix { {\bf GH}(X) \ar[r]^-{\varepsilon_{X}}\ar[d]_-{\simeq} &
X\ar@{=}[d]\ar[r]
& i_*i^*(X) \ar[r]\ar[d]^-{\simeq} & {\bf GH}(X)[1]\ar[d]^-{\simeq}\\
M\otimes_S\Hom_R(M, X)\ar[r]^-{\theta_{X}} & X\ar[r] &
\Coker(\theta_{X}) \ar[r] & M\otimes_S\Hom_R(M, X)[1]}
$$
With the help of this diagram and the recollement $(\ast)$ in Lemma
\ref{rt}, we have $i_*i^*(X)\in\mathcal{Y}$, and therefore
$$i_*i^*(X)\simeq\Coker(\theta_{X})\in \mathcal{Y}\cap
R\Modcat=\mathscr{E}.$$ This will finish the proof of $(1)$. So we
shall prove the existence of the above diagram $(a)$.

In fact, we shall first show that there exists a commutative diagram
$(b)$ in $\D{R}$:
$$
(b)\quad \xymatrix
{{\bf GH} (X) \ar[r]^-{\varepsilon_{X}}\ar[d]_-{\simeq} & X\ar@{=}[d]\\
M\otimes_S\Hom_R(M, X)\ar[r]^-{\theta_{X}} & X}
$$

This can be seen as follows: In Corollary \ref{counit}, we take $F:=
{}_RM\otimes_S-$ and $G:=\Hom_R(M,-)$. Then ${\bf G}=\mathbb{L}F$
and ${\bf H}=\mathbb{R}G$. To prove the existence of $(b)$, it
suffices to prove $X\in \mathcal{R}{_G}$ and $G(X)\in
\mathcal{L}{_F}$. For the definitions of $\mathcal{R}{_G}$ and
$\mathcal{L}{_F}$, we refer to Lemma \ref{homo}.

Observe that $X\in \mathcal{R}{_G}$ if and only if $\Ext^j_R(M,
X)=0$ for any $j>0$. Since $X\in\Prod(_RM)$, it suffices to show
that $\Ext^j_R(M, M^I)=0$ for any $j>0$ and any set $I$. This
follows from $\Ext^j_R(M, M^I)\simeq\Ext^j_R(M, M)^I=0$ by the axiom
$(R2)$. Thus $X\in \mathcal{R}{_G}$.

Note that $G(X)\in \mathcal{L}{_F}$ if and only if $\Tor^S_j(M,
G(X))=0$ for any $j>0$. Since $X\in\Prod(_RM)$ and $G$ commutes with
arbitrary direct products in $R\Modcat$,  we have $G(X)\in
\Prod(_SS)$. This means that, to prove $G(X)\in \mathcal{L}{_F}$, it
is sufficient to check $\Tor^S_j(M, S^I)=0$ for any $j>0$ and any
set $I$. However, since $M$ is a good Ringel module, the right
$S$-module $M$ is strongly $S$-Mittag-Leffler by the axiom $(R4)$,
and therefore $\Tor^S_j (M, S^I)=0$  by Lemma \ref{MLP} (3). This
shows  $G(X)\in \mathcal{L}{_F}$.

Hence, by Corollary \ref{counit}, the diagram $(b)$ does exist. Now,
by the recollement $(\ast)$ in Lemma \ref{rt}, we can extend
$\varepsilon_{X}$ to a canonical triangle in $\D{R}$: ${\bf
GH}(X)\lraf{\varepsilon_{X}} X \lra i_*i^*(X)\lra{\bf GH}(X)[1].$
Since each short exact sequence in $R\Modcat$ induces a canonical
triangle in $\D{R}$: $$M\otimes_S\Hom_R(M, X)\lraf{\theta_X} X \lra
\Coker(\theta_X)\lra M\otimes_S\Hom_R(M, X)[1],$$ the diagram $(a)$
follows from the commutative diagram $(b)$.

$(2)$ Since $M$ is a Ringel $R$-module, it follows from $(R3)$ that
there is a quasi-isomorphism $R\to \cpx{M}$ in $\K{R}$.
Consequently, we can form the following commutative diagram in
$\D{R}$:
$$
(c)\quad \xymatrix
{{\bf GH}(R) \ar[r]^-{\varepsilon_{R}}\ar[d]_-{\simeq} & R\ar[d]^-{\simeq}\\
{\bf GH}(\cpx{M}) \ar[r]^-{\varepsilon_{\cpx{M}}} & \cpx{M}}
$$

Next, using Corollary \ref{counit} again, we shall show that there
exists a commutative diagram in $\D{R}$:
$$
(d)\quad \xymatrix
{{\bf GH} (\cpx{M}) \ar[r]^-{\varepsilon_{\cpx{M}}}\ar[d]_-{\simeq} & \cpx{M}\ar@{=}[d]\\
M\otimes_S\Hom_R(M, \cpx{M})\ar[r]^-{\theta_{\cpx{M}}} & \cpx{M}}
$$
By Corollary \ref{counit}, we need only to show that $\cpx{M}\in
\mathcal{R}{_G}$ and $G(\cpx{M})\in \mathcal{L}{_F}$.

On the one hand, by the axiom $(R3)$ of Definition \ref{rm},
$\cpx{M}$ is a bounded complex such that each term of it belongs to
$\Prod(M)$. On the other hand, by Lemma \ref{homo}, the categories
$\mathcal{R}{_G}$ and $\mathcal{L}{_F}$ are triangulated
subcategories of $\K{R}$ and $\K{S}$, respectively. Thus, to prove
that $\cpx{M}\in \mathcal{R}{_G}$ and $G(\cpx{M})\in
\mathcal{L}{_F}$, it is enough to prove that $X\in \mathcal{R}{_G}$
and $G(X)\in \mathcal{L}{_F}$ for any $X\in\Prod(_RM)$. Clearly, the
latter has been shown in $(1)$. Thus $(d)$ follows directly from
Corollary \ref{counit}.

Note that $\theta_X: M\otimes_S\Hom_R(M,X)\lra X$ is injective by
$(1)$. Since $M_i\in\Prod(_RM)$ by the axiom $(R3)$, each map
$\theta_{M_i}$ is injective for $ 0\leq i\leq n$. This clearly
induces a complex $\Coker(\theta_{\cpx{M}})$ of the form:
$$
0\lra \Coker(\theta_{M_0}) \lraf{\partial_0}
\Coker(\theta_{M_1})\lraf{\partial_1} \cdots\lra
\Coker(\theta_{M_{n-1}})\lraf{\partial_{n-1}}
\Coker(\theta_{M_n})\lra 0\;\mbox{ in }\; \C R
$$ such that there is an exact sequence of complexes over $R$:
$$0 \lra  M\otimes_S\Hom_R(M, \cpx{M})\lraf{\theta_{\cpx{M}}}\cpx{M}\lra
\Coker\big(\theta_{\cpx{M}}\big) \lra 0.$$ Since each exact sequence
of complexes over $R$ can be naturally extended to a canonical
triangle in $\D{R}$, we obtain a triangle in $\D{R}$:
$$(e)\quad M\otimes_S\Hom_R(M, \cpx{M})\lraf{\theta_{\cpx{M}}}\cpx{M}\lra
\Coker\big(\theta_{\cpx{M}}\big) \lra M\otimes_S\Hom_R(M,
\cpx{M})[1].
$$
Certainly, we also have a canonical triangle in $\D{R}$ from the
recollement $(\ast)$ in Lemma \ref{rt}:
$$
(f)\quad {\bf GH}(R)\lra R \lra i_*i^*(R)\lra{\bf GH}(R)[1].
$$
So, combining $(c)$, $(d)$, $(e)$ with $(f)$, one can easily
construct the following commutative diagram in $\D{R}$:
$$ \xymatrix { {\bf GH}(R)
\ar[r]^-{\varepsilon_{R}}\ar[d]_-{\simeq} & R\ar[d]^-{\simeq}\ar[r]
& i_*i^*(R) \ar[r]\ar[d]^-{\simeq} & {\bf GH}(R)[1]\ar[d]^-{\simeq}\\
M\otimes_S\Hom_R(M, \cpx{M})\ar[r]^-{\theta_{\cpx{M}}} &
\cpx{M}\ar[r] & \Coker(\theta_{\cpx{M}}) \ar[r] &
M\otimes_S\Hom_R(M, \cpx{M})[1]}
$$
In particular, we have $i_*i^*(R)\simeq \Coker(\theta_{\cpx{M}})$ in
$\D{R}$, and therefore
$$H^j(i_*i^*(R))\simeq H^j\big(\Coker(\theta_{\cpx{M}})\big)\;
\mbox{for any}\; j\in \mathbb{Z}.$$ This implies that
$H^j(i_*i^*(R))=0$ for $j<0$ or $j>n$.

Now, combining $(e)$ with $R\simeq \cpx{M}$ in $\D{R}$, we obtain a
triangle in $\D{R}$:
$$ M\otimes_S\Hom_R(M, \cpx{M})\lra R\lra
\Coker\big(\theta_{\cpx{M}}\big) \lra M\otimes_S\Hom_R(M,
\cpx{M})[1].
$$
Applying the cohomology functor $H^j$ to this triangle, one can
check that
$$H^j(i_*i^*(R))\simeq H^j\big(\Coker(\theta_{\cpx{M}})\big)
\simeq H^{j+1}( M\otimes_S\Hom_R(M, \cpx{M}))\; \mbox{for any}\;
j>0.$$ Thus $(2)$ follows.

$(3)$ For $n=0$, the conclusion follows from $i_*i^*(R)\simeq
\Coker(\theta_{\cpx{M}})$ trivially. So, we may assume $n\ge 1$. By
the final part of the proof of $(2)$, we know that
$$ i_*i^*(R)\simeq \Coker(\theta_{\cpx{M}})\; \mbox{ in}\;
\D{R}\quad \mbox{and}\quad
H^n\big(\Coker(\theta_{\cpx{M}})\big)\simeq H^{n+1}(
M\otimes_S\Hom_R(M, \cpx{M})).$$ Since the $(n+1)$-term of the
complex $M\otimes_S\Hom_R(M, \cpx{M})$ is zero, we see that
$H^n\big(\Coker(\theta_{\cpx{M}})\big)=0$. This implies that the
$(n-1)$-th differential $\partial_{n-1}$ of the complex
$\Coker(\theta_{\cpx{M}})$ is surjective. It follows that
$\Coker(\theta_{\cpx{M}})$ is isomorphic in $\D{R}$ to the following
complex:
$$ (\dag)\quad
0\lra \Coker(\theta_{M_0}) \lraf{\partial_0}
\Coker(\theta_{M_1})\lraf{\partial_1} \cdots\lra
\Coker(\theta_{M_{n-2}})\lraf{\partial{_{n-2}}}\;\Ker(\partial_{n-1})\lra
0.
$$
Since $M_m\in\Prod(_RM)$ for $0\leq m\leq n$ by the axiom $(R3)$, we
see from $(1)$ that $\Coker(\theta_{M_m})\in \mathscr{E}$. Note that
$\mathscr{E}$ is always closed under kernels of surjective
homomorphisms in $R\Modcat$. Thus
$\Ker(\partial_{n-1})\in\mathscr{E}$. This means that $(\dag)$ is a
bounded complex with all of its terms in $\mathscr{E}$.

Consequently, the complex $i_*i^*(R)$ is isomorphic in $\D{R}$ to
the complex $(\dag)$ with the required form in Lemma \ref{prep} (3).
This finishes the proof. $\square$

\medskip
{\it Remark.} By the proof of Lemma \ref{prep} (2), we see that the
complex ${_R}M\otimes_S\Hom_R(M,\,\cpx{M})$ is isomorphic in $\D{R}$
to both ${_R}M\otimesL_S\Hom_R(M,\,\cpx{M})$ and ${\bf GH}(R)$. This
implies that, up to isomorphism, the cohomology groups
$H^j\big({_R}M\otimes_S\Hom_R(M,\,\cpx{M})\big)$, for
$j\in\mathbb{Z}$, are independent of the choice of the complex
$\cpx{M}$ which arises in the axiom $(R3)$ of Definition \ref{rm}.

\medskip
With the help of Lemma \ref{hep} and Lemma \ref{prep}, we can prove
the following key proposition.

\begin{Prop}\label{realization} The following statements are equivalent:

$(1)$ The full triangulated subcategory $\mathcal{Y}$ of $\D{R}$ is
homological.

$(2)$ The category $\mathscr{E}$ is an  abelian subcategory of
$R\Modcat$.

$(3)$ $H^j\big({_R}M\otimes_S\Hom_R(M,\,\cpx{M})\big)=0$ for any
$j\geq 2$.

$(4)$ The kernel of the  homomorphism $\partial_0:
\Coker(\theta_{M_0})\lra \Coker(\theta_{M_1})$ induced from $\nu$
belongs to $\mathscr{E}$.
\end{Prop}

{\it Proof.} The equivalences of $(1)$ and $(2)$ follow from those
of $(1)$ and $(6)$ in Lemma \ref{hep} together with Lemma \ref{prep}
(3), while the equivalences of $(1)$ and $(3)$ follow from those of
$(1)$ and $(2)$ in Lemma \ref{hep} together with Lemma \ref{prep}
(2). Now we prove that $(1)$ and $(4)$ are equivalent. By Lemma
\ref{prep} (2) and the equivalence of $(1)$ and $(3)$ in Lemma
\ref{hep}, we see that (1) is equivalent to $H^0(i_*i^*(R))\in
\mathcal{Y}$. By the proof of Lemma \ref{prep} (2), we infer that
$H^0(i_*i^*(R))\simeq H^0(\Coker(\theta_{\cpx{M}}))\simeq
\Ker(\partial_0)$. Thus, $(1)$ is equivalent to $\Ker(\partial_0)\in
\mathcal{Y}\cap \Modcat{R}=\mathscr{E}$. $\square$

\medskip
As a consequence of Proposition \ref{realization}, we have the
following handy characterizations.

\begin{Koro}\label{appl}
Assume that the projective dimension of $_RM$ is equal to $n$. Then
the following are true.

$(1)$ If $n\leq 1$, then $\mathcal{Y}$ is always homological.

$(2)$ If $n=2$, then $\mathcal{Y}$ is homological if and only if
$M\otimes_S\Ext_R^2(M, R)=0$.

$(3)$ Suppose that $n\geq 3$ and $\Tor_i^S(M,\,\Ext_R^j(M,
R))=0\quad \mbox{for}\quad 2\leq j\leq n-1\quad \mbox{and}\quad
0\leq i\leq j-2.$ Then $\mathcal{Y}$ is homological if and only if
$$
\Tor_k^S(M, \,\Ext_R^n(M, R))=0\quad \mbox{for}\quad 0\leq k\leq
n-2.
$$
\end{Koro}

\medskip
{\it Proof.} The key point in the proof is to check when the $j$-th
cohomology group $H^j\big({_R}M\otimes_S\Hom_R(M,\,\cpx{M})\big)$
vanishes for $j\geq 2$. Note that
$H^j\big(M\otimes_S\Hom_R(M,\,\cpx{M})\big)=0$ for all $j>n$.

For $n\le 1$, the conclusion in Corollary \ref{appl} is clear. So,
we suppose $n\geq 2$. By the axiom $(R2)$, we have $\Ext^j_R(M,
M)=0$ for all $j\geq 1$. It follows that $\Ext^j_R(M, M^I)\simeq
\Ext^j_R(M, M)^I=0$ for any nonempty set $I$, and therefore
$\Ext^j_R(M, X)=0$ for any $X\in\Prod(M)$.

By the axiom $(R3)$, there exists an exact sequence in $R\Modcat$:
$$
0\lra R \lra M_0 \lra M_1\lra \cdots \lra M_n\lra 0$$ such that
$M_i\in\Prod(M)$ for $0\leq i\leq n$. Since $\Ext^j_R(M, X)=0$ for
any $X\in\Prod(M)$ and $j\ge 1$, we know that  the following complex
$\Hom_R(M, \cpx{M}):$
$$0\lra
\Hom_R(M, M_0) \lra \Hom_R(M, M_1) \lra\Hom_R(M, M_2)\lra  \cdots
\lra \Hom_R(M, M_n)\lra 0$$ satisfies that $H^j\big(\Hom_R(M,
\cpx{M})\big)\simeq \Ext^j_R(M, R)$ for each $j\geq 1$.

$(2)$ Let $n=2$. Consider the complex
$M\otimes_S\Hom_R(M,\,\cpx{M}):$
$$0\lra M\otimes_S\Hom_R(M, M_0) \lra M\otimes_S\Hom_R(M, M_1) \lra
M\otimes_S\Hom_R(M, M_2)\lra 0.$$ Since the functor $_RM\otimes_S-:
S\Modcat\to R\Modcat$ is right exact, we have
$$H^2\big(M\otimes_S\Hom_R(M,\,\cpx{M})\big)\simeq M\otimes_S
H^2\big(\Hom_R(M,\,\cpx{M})\big)\simeq M\otimes_S\Ext_R^2(M, R).$$
Now, the statement $(2)$ follows from the equivalences of $(1)$ and
$(3)$ in Proposition \ref{realization}.

$(3)$ Under the assumption of $(3)$, we claim that
$$H^m\big(M\otimes_S\Hom_R(M,\,\cpx{M})\big)\simeq \Tor_{n-m}^S(M,
\,\Ext_R^n(M, R))\quad \mbox{for}\quad 2\leq m\leq n.$$
Consequently, the statement $(3)$ will follow from the equivalences
of $(1)$ and $(3)$ in Proposition \ref{realization}.

In the following, we shall apply Lemma \ref{Formula} to prove this
claim. Define $\cpx{Y}:=\Hom_R(M,\,\cpx{M})$. This is a complex over
$S$ with $Y^i=\Hom_R(M, M_i)$ for $0\leq i\leq n$ and $Y^i=0$ for
$i\geq n+1$. Moreover, since the right $S$-module $M$ is strongly
$S$-Mittag-Leffler by the axiom $(R4)$, it follows from the proof of
Lemma \ref{prep} (1) that
$$\Tor^S_k\big(M, \Hom_R(M, X)\big)=0\; \mbox{ for all}\; k\geq
1\;\mbox{and}\; X\in\Prod(M).$$ This implies that $\Tor^S_k(M,
Y^i)=0$ for all $i\in \mathbb{Z}$ and $k\geq 1$.

Recall that $H^j(\cpx{Y})\simeq \Ext^j_R(M, R)$ for all $j\geq 1$.
By assumption, we obtain
$$\Tor_i^S(M,\,H^j(\cpx{Y}))=0\quad \mbox{for}\quad 2\leq j\leq n-1\quad \mbox{and}\quad
0\leq i\leq j-2.$$ Clearly, this implies that, for each $2\leq m\leq
n-1$, we have
$$\,\Tor^S_t\big(M,
H^{m+t}(\cpx{Y})\big)=0=\Tor^S_{t-1}\big(M, H^{m+t}(\cpx{Y})\big)
\;\mbox{for}\;\, 0\leq t \leq n-m-1.$$ It follows from Lemma
\ref{Formula} that $H^m(M\otimes_S\cpx{Y})\simeq \Tor^S_{n-m}\big(M,
H^n(\cpx{Y})\big)\simeq\Tor_{n-m}^S(M, \,\Ext_R^n(M, R))$.

To finish the proof of the claim, it remains to prove
$H^n(M\otimes_S\cpx{Y})\simeq M\otimes_S\Ext_R^n(M, R)$. However,
since the functor $M\otimes_S-$ is right exact and since $Y^i=0$ for
$i\geq n+1$, we see that $H^n(M\otimes_S\cpx{Y})\simeq M\otimes_S
H^n(\cpx{Y})\simeq M\otimes_S \Ext_R^n(M, R)$. This finishes the
proof of the above-mentioned claim. Thus $(3)$ holds.  $\square$

\medskip
As another consequence of Proposition \ref{realization}, we mention
the following result which is not used in this note, but of its own
interest.

\begin{Koro}\label{appl'}
$(1)$ If $M_0\in\Add(_RM)$, then $_RM$ is a classical tilting
module.

$(2)$ If $M_1\in\Add(_RM)$, then $\mathcal Y$ is homological in $\D
R$.
\end{Koro}

{\it Proof.} (1) Suppose $M_0\in\Add(_RM).$ We claim that
$\Coker(\theta_{M_0})=0$. In fact, since $_RM$ is finitely generated
by the axiom $(R1)$, the functor $\Hom_R(M, -): R\Modcat\to
S\Modcat$ commutes with arbitrary direct sums. It follows that the
evaluation map $\theta_{X}: M\otimes_S\Hom_R(M, X)\lra X$ is an
isomorphism for each $X\in\Add(_RM)$. Since $M_0\in\Add(_RM)$, the
map $\theta_{M_0}: M\otimes_S\Hom_R(M, M_0)\lra M_0$ is an
isomorphism, and therefore $\Coker(\theta_{M_0})=0$. Combining this
with the proof of Proposition \ref{realization}, we have
$H^0(i_*i^*(R))\simeq \Ker(\partial_0)=0$. Note that $\End_{\D
R}(i^*(R))\simeq H^0(i^*(R))=H^0(i_*i^*(R))$ as $R$-modules by Lemma
\ref{ringhom} (2). This implies that $\End_{\D R}(i^*(R))=0$ and so
$\mathcal{Y}=0$ by Lemma \ref{ringhom} (1). Now, it follows from
Lemma \ref{rt} that $\rHom_R(M,-):\; \D{R}\lra\D{S}$ is a triangle
equivalence. Consequently, $_RM$ is a classical tilting module by
\cite[Chapter 5, Theorem 4.1]{HHK}.

(2) It follows from the proof of (1) that $\Coker(\theta_{M_1})=0$.
Thus (2) follows from Proposition \ref{realization} and Lemma
\ref{prep} (1). $\square$

\section{Application to tilting modules:
Proofs of Theorem \ref{main-result} and Corollary \ref{cor}}\label{sect5}

In this section, we first develop some properties of (good) tilting
modules, and then give a method to construct good Ringel modules.
With these preparations in hand, we finally apply Proposition
\ref{realization} to prove Theorem \ref{main-result} and Corollary
\ref{cor}.

\smallskip
Throughout this section, $A$ will be a ring and $n$ a natural
number. In addition, we assume that $T$ is a {\bf{good}}
 $n$-tilting $A$-module with $(T1), (T2)$ and $(T3)'$. Let
$B:=\End_A(T)$.

\medskip
First of all, we shall mention a few basic properties of good
tilting modules in the following lemma. For proofs, we refer to
\cite[Chapter 11, Lemma 2.7]{HHK}, \cite[Proposition 1.4, Lemma
1.5]{Bz2} and \cite[Proposition 3.5]{Bz1}.

\begin{Lem} \label{tilt1}
The following hold true for the tilting module $_AT$.

$(1)$ The torsion class  $T^{\bot}:=\{X\in A\Modcat \mid \Ext^i
_A(T,X)=0 \mbox{\;for\, all \,} i\geq 1\}$ in $A\Modcat$ is closed
under arbitrary direct sums in $A\Modcat$.

$(2)$ The right $B$-module $T$ has a finitely generated projective
resolution of length at most $n$:
$$
0\lra \Hom_A(T_n, T)\lra \cdots \lra\Hom_A(T_1, T)\lra\Hom_A(T_0,
T)\lra T_B\lra 0$$ with  $T_i\in\add(_AT)$ for all $0\leq i\leq n$.

$(3)$ The map $A^{\opp} \to \End{_{B^{\opp}}}(T) $, defined by
$a\mapsto [t\mapsto at]$ for $a\in A$ and $t\in T$, is an
isomorphism of rings. Moreover, $\Ext^i_{B^{\opp}}(T,T)=0$ for all $
i\geq 1$.

$(4)$ If $T_n=0$ in the axiom $(T3)'$, then $_AT$ is an
$(n-1)$-tilting module.
\end{Lem}

Let us introduce some notation which will be used throughout this
section.

Define
$$
G:={_A}T\otimesL_B- : \,\D{B}\to \D{A},\quad H:=\rHom_A(T,-):\,
\D{A}\to \D{B},
$$
$$\cpx{Q}:=\;\cdots\lra 0\lra \Hom_A(T, T_0)\lra\Hom_A(T, T_1)\lra
\cdots \lra\Hom_A(T, T_n)\lra 0\lra \cdots$$ where $\Hom_A(T, T_i)$
is of degree $i$ for $0\leq i\leq n$, and
$\cpx{Q}{^*}:=\Hom_B(\cpx{Q}, B)\in \C{\pmodcat{B\opp}}$. Clearly,
$\cpx{Q}{^*}$ is isomorphic in $\Cb{\pmodcat{B\opp}}$ to the complex
$$
\;\cdots\lra 0\lra \Hom_A(T_n, T)\lra \cdots \lra\Hom_A(T_1,
T)\lra\Hom_A(T_0, T)\lra 0\lra \cdots $$

\medskip

The following result is due to Bazzoni \cite[Theorem\,2.2]{Bz2},
which says that, in general, $\D{A}$ is not equivalent to $\D B$,
but a full subcategory of $\D{B}$.

\begin{Lem}\label{Bz}
The functor $H: \D{A}\to \D{B}$ is fully faithful, and
$\Img(H)=\Ker(\Hom_{\D B}(\Ker(G), -))$.
\end{Lem}
\smallskip

The next result supplies a way to understand good tilting modules
$T$ by some special objects or by subcategories of derived module
categories. In particular, the category $\Ker(G)$ is a bireflective
subcategory of $\D{B}$.

\begin{Lem}\label{tilt2} For the tilting $A$-module $T$, we have
the following:

$(1)$ $H(A)\simeq \cpx{Q}$ in $\D{B}$ and
$\Hom_{\D{B}}(\cpx{Q},\cpx{Q}[m])=0$ for any $m\ne 0$.

$(2)$ $\Ker(G)=\{\cpx{Y}\in\D{B}\mid \Hom_{\D B}(\cpx{Q},\,
\cpx{Y}[i])=0 \mbox{\;for\, all \,} i\in\mathbb{Z}\}.$

$(3)$ Let $j_!:{\rm{Tria}}{(\cpx{Q})}\lra \D{B}$ and $i_*:
\Ker(G)\lra \D{B}$ be the inclusions. Then there exists a
recollement of triangulated categories together with a triangle
equivalence:
$$(\star)\quad
\xymatrix@C=1.2cm{\Ker(G) \ar[r]^-{{i_*}}
&\D{B}\ar[r]^{j^!}\ar@/^1.2pc/[l]\ar_-{i^*}@/_1.2pc/[l]
&{\rm{Tria}}{(\cpx{Q})}\ar[r]^-{G\,j_*}_-{\simeq}
\ar@/^1.2pc/[l]_{j_*}\ar@/_1.2pc/[l]_{j_!\;} & \D{A} }$$

\bigskip
\noindent such that $G\,j_*\,j^!$ is naturally isomorphic to $G$.
\end{Lem}

{\it Proof.} We remark that Lemma \ref{tilt2} is implied in
\cite{Bz2}. For convenience of the reader, we give a proof here.

$(1)$ By the axiom $(T3)'$, the stalk complex $A$ is
quasi-isomorphic in $\C{A}$ to the complex $\cpx{T}$ of the form:
$$ \cdots \lra 0 \lra T_0 \lra
T_1\lra \cdots \lra T_n\lra 0\lra \cdots$$ where $T_i\in \add(T)$ is
in degree $i$ for $0\leq i\leq n$.  Further, by the axiom $(T2)$, we
have $T_i\in T^{\bot}:=\{X\in A\Modcat \mid \Ext^i _A(T,X)=0
\mbox{\;for\, all \,} i\geq 1\}$. It follows from Lemma \ref{homo}
(1) that $H(A)\simeq H(\cpx{T})\simeq \Hom_A(T, \cpx{T})=\cpx{Q}$ in
$\D{B}$. Since the functor $H$ is fully faithful by Lemma \ref{Bz},
we obtain
$$\Hom_{\D B}(\cpx{Q}, \cpx{Q}[m])\simeq \Hom_{\D B}(H(A),
H(A)[m])\simeq \Hom_{\D A}(A, A[m])\simeq \Ext^m_A(A, A)=0$$ for any
$m\neq 0$. This shows $(1)$.

$(2)$ Since $\cpx{Q}\in \Cb{\pmodcat{B}}$ and since $\cpx{Q}{^*}$ is
quasi-isomorphic to $T_B$ by Lemma \ref{tilt1} (2), we have the
following natural isomorphisms of triangle functors:
$$\rHom{_B}(\cpx{Q},-)\lraf{\simeq} \cpx{Q}{^*}\otimesL_B-\lraf{\simeq} {_\mathbb{Z}}T\otimesL_B-: \D{B}\lra\D{\mathbb
Z},$$ where the first isomorphism follows from Lemma \ref{complex}.
Note that  $H^m(\rHom_{\D B}(\cpx{Q}, \cpx{Y}))\simeq \Hom_{\D
B}(\cpx{Q}, \cpx{Y}[m])$ for $m\in\mathbb{Z}$ and $\cpx{Y}\in\D{B}$.
This shows $(2)$.

$(3)$  Since $\cpx{Q}\in \Cb{\pmodcat{B}}$, we know from $(2)$ and
Lemma \ref{rg} (1) that there exists a recollement of triangulated
categories:
$$(\star\star)\quad
\xymatrix@C=1.2cm{\Ker(G)\ar[r]^-{{i_*}}
&\D{B}\ar[r]^{j^!}\ar@/^1.2pc/[l]_-{i^!}\ar_-{i^*}@/_1.2pc/[l]
&{\rm{Tria}}{(\cpx{Q})}
\ar@/^1.2pc/[l]_{j_*}\ar@/_1.2pc/[l]_{j_!\;}}$$

\noindent On the one hand, by the correspondence of recollements and
TTF (torsion, torsion-free) triples (see, for example, \cite[Section
2.3]{CX1}), we infer from $(\star\star)$ that
$\Img(j_*)=\Ker(\Hom_{\D B}(\Ker(G), -))$ and that the functor $j_*:
{\rm{Tria}}{(\cpx{Q})}\lra \Img(j_*)$ is a triangle equivalence with
the restriction of $j^!$ to $\Img(j_*)$ as  its quasi-inverse. On
the other hand, it follows from Lemma \ref{Bz} that
$\Img(H)=\Ker(\Hom_{\D B}(\Ker(G), -))$ and the functor $H:\D{A}\
\lra \Img(H)$ is a triangle equivalence with the restriction of $G$
to $\Img(H)$ as its quasi-inverse. Consequently, we see that
$\Img(j_*)=\Img(H)$ and the composition $G\,j_*:
{\rm{Tria}}{(\cpx{Q})}\lra \D{A} $ of $j_*$ with $G$ is also a
triangle equivalence.

It remains to check
$$G\lraf{\simeq} G\,j_* j^!: \;\D{B}\lra \D{A}.$$ In fact, for
any $\cpx{X}\in\D{B}$, by the recollement $(\star\star)$, there
exists a canonical triangle in $\D{B}:$
$$i_*i^!(\cpx{X})\lra \cpx{X}\lra j_*j^!(\cpx{X})\lra
i_*i^!(\cpx{X})[1].$$  Since $\Img(i_*i^!)=\Img(i_*)=\Ker(G)$, we
know that $G(\cpx{X}) \lraf{\simeq} Gj_*j^!(\cpx{X})$ in $\D{B}$.
This proves $(3)$. $\square$

\medskip
Next, we shall investigate when the subcategory $\Ker(G)$ of $\D{B}$
is homological. The following result conveys that this discussion
can be proceeded along the right $B$-module $T$.

\begin{Lem}\label{tilt4}
The category $\Ker(G)$ is a homological subcategory of $\D{B}$ if
and only if  $\,\Ker\big(\,\rHom_{B\opp}(T, -)\big)$ is a
homological subcategory of $\D{B\opp}$.
\end{Lem}

{\it Proof.} In Lemma \ref{rg}, we take $R:=B$ and
$\Sigma:=\{\cpx{Q}\}$. Then $\Sigma^*=\{\cpx{Q}{^*}\}$ where
$\cpx{Q}{^*}:=\Hom_B(\cpx{Q}, B)$. Since $\cpx{Q}{^*}$ is
quasi-isomorphic to $T_B$ by Lemma \ref{tilt1} (2), we infer that
$\cpx{Q}{^*}\lraf{\simeq} T_B$ in $\D{B\opp}$ and that there exists
a natural isomorphism of triangle functors:
$$\rHom{_{B\opp}}(T,\,-)\lraf{\simeq}
\rHom_{B\opp}(\cpx{Q}{^*},\,-): \D{B\opp}\lra \D{\mathbb{Z}}.$$ This
implies that $$\Ker\big(\rHom_{B\opp}(T,\,
-)\big)=\Ker\big(\rHom_{B\opp}(\cpx{Q}{^*},\, -)\big)=\{\cpx{Y}\mid
\Hom_{\D {B\opp}}(\cpx{Q}{^*}, \,\cpx{Y}[m])=0\; \mbox{for}\;
m\in\mathbb{Z}\}. $$  Thus Lemma \ref{tilt4} follows from Lemmas
\ref{ref} and \ref{rg} (4). $\square$

\medskip
Next, we point out that each good tilting module naturally
corresponds to a good Ringel module. This guarantees that we can
apply Proposition \ref{realization} to show Theorem
\ref{main-result}.

\begin{Lem}\label{tilt3}
The right $B$-module $T_B$ is a good $n$-Ringel module.
\end{Lem}

{\it Proof.}  By Lemma \ref{tilt1} (2), the axiom $(R1)$ holds for
$T_B$, and the projective dimension of $T_B$ is at most $n$.
Moreover, by Lemma \ref{tilt1} (3), the axiom $(R2)$ also holds for
$T_B$. Now, we check the axiom $(R3)$ for $T_B$.

In fact, according to the axiom $(T1)$, the module $_AT$ admits a
projective resolution of $A$-modules:
$$
0\lra P_n \lra \cdots \lra P_1\lra P_0\lra T\lra 0$$ with
$P_i\in\Add(_AA)$ for $0\leq i\leq n$. Since $\Ext^j_A(T, T)=0$ for
each $j\geq 1$ by the axiom $(T2)$, it follows that the sequence
$$
0\lra B \lra\Hom_A(P_0, T)\lra \Hom_A(P_1, T)\lra \cdots \lra
\Hom_A(P_n, T) \lra 0$$ of right $B$-modules is exact. Note that
$\Hom_A(P_i, T)\in \Prod(T_B)$ due to $P_i\in\Add(_AA)$. This means
that the axiom $(R3)$ holds for $T_B$. Thus the right $B$-module
$T_B$ is an $n$-Ringel module.

It remains to prove that $T_B$ is good, that is, $T_B$ satisfies the
axiom $(R4)$.

Actually, by Lemma \ref{tilt1} (3), the map $A^{\opp}\lra
\End{_{B^{\opp}}}(T) $, defined by $a\mapsto [t\mapsto at]$ for
$a\in A$ and $t\in T$, is an isomorphism of rings. Further, it
follows from Lemma \ref{TML} that the right $A\opp$-module $T$ is
strongly $A\opp$-Mittag-Leffler. Hence, the right
$\End{_{B^{\opp}}}(T)$-module $T$ is strongly
$\End{_{B^{\opp}}}(T)$-Mittag-Leffler. Thus, by definition, the
$n$-Ringel $B\opp$-module $T$ is good. $\square$

\smallskip
{\it Remark.}  If $_AT$ is infinitely generated, then the right
$B$-module $T$ is not a tilting module. In fact, it follows from
Lemma \ref{tilt1} (2) that $T_B$ is finitely generated. Suppose
contrarily that $T_B$ is a tilting right $B$-module. Then, by
Corollary \ref{2.9}, the right $B$-module $T_B$ is classical, and
therefore $_AT$ is classical by Lemma \ref{tilt1} (2)-(3). This is a
contradiction.

\medskip
Now, with the previous preparations, we are in the position to prove
Theorem \ref{main-result}.

\medskip
{\bf Proof of Theorem \ref{main-result}.} We shall use Proposition
\ref{realization} to show the equivalences in Theorem
\ref{main-result}.

Recall that we denote by $\cpx{P}$ the complex which is the deleted
projective resolution of $_AT$:
$$
\cdots\lra 0\lra P_n \lra \cdots \lra P_1\lraf{\sigma} P_0\lra 0\lra
\cdots
$$ appearing in the axiom $(T1)$. Here, $P_i$ is in degree $-i$ for $0\le i\le
n$.

By Lemma \ref{tilt3}, we know that $T$ is a good $n$-Ringel
$B\opp$-module and that the exact sequence in the axiom $(R3)$ can
be chosen as
$$
0\lra B_B \lra\Hom_A(P_0, T)\lra \Hom_A(P_1, T)\lra \cdots \lra
\Hom_A(P_n, T) \lra 0.$$ In particular, the complex $\cpx{M}$ in
Proposition \ref{realization} can be chosen to be the following
complex:
$$
\Hom_A(\cpx{P}, T):\; \cdots \lra 0\lra \Hom_A(P_0, T)\lra
\Hom_A(P_1, T)\lra \cdots \lra \Hom_A(P_n, T) \lra 0\lra \cdots$$
Now, in Proposition \ref{realization}, we take $R:=B\opp\,$,
$S:=A\opp$ and $M:= {}_RT_S$.  Further, let $${\bf
H}=\rHom_{B\opp}(T, -):\D{B\opp}\lra \D{A\opp}.$$ It follows from
Lemma \ref{tilt4} that $\Ker(G)$ is homological in $\D R$ if and
only if so is $\Ker({\bf H})$ in $\D{B\opp}$. In other words, the
statement $(1)$ in Theorem \ref{main-result} is equivalent to the
following statement:

$(1')$ The category  $\Ker({\bf H})$ is a homological subcategory of
$\D{B\opp}$.

In the following, we shall show that $(1')$ is equivalent to $(2),
(3)$ and $(4)$, respectively.

We first show that $(1')$ and $(2)$ are equivalent. In fact, it
follows form Proposition \ref{realization} that $(1')$ is equivalent
to

$(2')$ The category $\mathscr{E}:=\{Y\in {B\opp}\Modcat \mid
\Ext^m_{B\opp}(T, Y)=0\mbox{\;for\, all \,} m\geq0\}$ is an abelian
subcategory of $B\opp\Modcat$.

\noindent So, we will show that $(2')$ is equivalent to $(2)$. For
this aim, we set $\mathscr{A}:=\{X\in B\Modcat \mid \Tor_m
^B(T,X)=0\mbox{\;for\, all \,} m\geq0\}$, and establish a connection
between $\mathscr{A}$ and $\mathscr{E}$. Let $(-)^\vee$ be the
 dual functor $\Hom_\mathbb{Z}(-, \mathbb{Q}/\mathbb{Z}):
\mathbb{Z}\Modcat\lra \mathbb{Z}\Modcat$.

Now, we claim that $(-)^\vee$ induces two exact functors:
$$
(-)^\vee: \mathscr{A}\lra \mathscr{E}\quad \mbox{and}\quad (-)^\vee:
\mathscr{E}\lra \mathscr{A}
$$
such that $X\in\mathscr{A}$ if and only if $X^\vee\in\mathscr{E}$,
and that $Y\in\mathscr{E}$ if and only if $Y^\vee\in\mathscr{A}$,
where $X\in B\Modcat$ and $Y\in B\opp\Modcat$.

In fact, it is known that $\mathbb{Q}/\mathbb{Z}$ is an injective
cogenerator for $\mathbb{Z}\Modcat$, and that $(-)^\vee$ admits the
following properties:

$(a)$ For each $M\in\mathbb{Z}\Modcat$, if $M^\vee=0$, then $M=0$.

$(b)$ A sequence $0\to X_1\to X_2\to X_3\to 0$ of
$\mathbb{Z}$-modules is exact if and only if $0\to (X_3)^\vee \to
(X_2)^\vee \to (X_1)^\vee \to 0$ is exact.

On the one hand, for each $X\in B\Modcat$, it follows from Lemma
\ref{tor-ext} (1) that
$$
(\Tor_m^B(T, \,X))^\vee \simeq \Ext^m_{B\opp}(T,\, X^\vee)\;\,
\mbox{for all}\;\, m\geq 0.
$$
This implies that $X\in\mathscr{A}$ if and only if
$X^\vee\in\mathscr{E}$. This is due to $(a)$.

On the other hand, since $T_B$ has a finitely generated projective
resolution in $B\opp\Modcat$ by Lemma \ref{tilt1} (2), it follows
from Lemma \ref{tor-ext} (2) that
$$
(\Ext_{B\opp}^m(T,\,Y))^\vee \simeq \Tor^{B}_m(T, \,Y^\vee)\;\,
\mbox{ for all}\;\, m\geq 0 \mbox{  and for any } Y\in
B\opp\Modcat.$$ This means that $Y\in\mathscr{E}$ if and only if
$Y^\vee\in\mathscr{A}$, again due to $(a)$. This finishes the proof
of the claim.

Recall that $\mathscr{A}$ always admits the ``$2$ out of $3$"
property: For an arbitrary short exact sequence in $B\Modcat$, if
any two of its three terms belong to $\mathscr{A}$, then so does the
third. Moreover, $\mathscr{A}$ is an abelian subcategory of
$B\Modcat$ if and only if $\mathscr{A}$ is closed under kernels
(respectively, cokernels) in $B\Modcat$. Clearly, similar statements
hold for the subcategory $\mathscr{E}$ of $B\opp\Modcat$.

By the above-proved claim,  one can easily show that $\mathscr{A}$
is closed under kernels in $B\Modcat$ if and only if $\mathscr{E}$
is closed under cokernels in $B\opp\Modcat$. It follows that
$\mathscr{A}$ is an abelian subcategory of $B\Modcat$ if and only if
$\mathscr{E}$ is an abelian subcategory of $B\opp\Modcat$. Thus
$(2')$ is equivalent to $(2)$, and therefore  $(1')$ and $(2)$ are
equivalent.

Next, we shall verify that $(1')$ and $(3)$ are equivalent.
Actually, it follows form Proposition \ref{realization} that $(1')$
is also equivalent to the following statement:

\smallskip
$(3')$ $H^j\big(\Hom_{B\opp}(T,\,\cpx{M})\otimes_AT\big)=0$ for all
$j\geq 2$, where $\Hom_{B\opp}(T,\,\cpx{M}):=\Hom_{B\opp}(T,
\Hom_A(\cpx{P}, T))$ is the complex of the form: {\small{$$0\lra
\Hom_{B\opp}(T, \Hom_A(P_0, T))\lra \Hom_{B\opp}(T, \Hom_A(P_1,
T))\lra \cdots \lra \Hom_{B\opp}(T, \Hom_A(P_n, T)) \lra 0,
$$}}
with $\Hom_{B\opp}(T, \Hom_A(P_i, T))$ in degree $i$ for $0\leq
i\leq n$.

So it suffices to verify that $(3')$ and $(3)$ are equivalent.
Clearly, for this purpose, it is enough to show that
$\Hom_A(\cpx{P}, A) \simeq \Hom_{B\opp}(T, \Hom_A(\cpx{P}, T))$ as
complexes over $A\opp$.

Note that there exists a natural isomorphism of additive functors:
$$
\Hom_{B\opp}(T, \Hom_A(-, T))\lraf{\simeq} \Hom_{B\opp}(\Hom_A(A,
T),\,\Hom_A(-, T)): A\Modcat\to A\opp\Modcat.
$$
Moreover, the functor $\Phi:=\Hom_A(-, T)$ yields a natural
transformation:
$$
\Hom_A(-, A) \lra\Hom_{B\opp}(\Phi(A),\,\Phi(-)\,): A\Modcat\to
A\opp\Modcat.
$$
Now we shall show that this transformation is even a natural
isomorphism. Clearly, it is sufficient to prove that $$\Phi:
\Hom_A(X, A) \lraf{\simeq} \Hom_{B\opp}(\Phi(A),\,\Phi(X))$$  for
any projective $A$-module $X$. In the following, we will show that
this holds even for any $A$-module $X$.

In fact, since $T$ is a good tilting $A$-module, it follows from the
axiom $(T3)'$ that there exists an exact sequence $0\lra A \lra T_0
\lra T_1$ with $T_i\in\add(T)$ for $i=0, 1$. By Lemma \ref{tilt1}
(2), we obtain another exact sequence $\Phi(T_1)\lra\Phi(T_0)\lra
\Phi(A)\lra 0$ of $B\opp$-modules. This gives rise to the following
exact commutative diagram:
$$\xymatrix{
0\ar[r] & \Hom_A(X, A)\ar[r]\ar[d]^-{\Phi} & \Hom_A(X, T_0)
\ar[r]\ar[d]^-{\simeq}
&\Hom_A(X, T_1)\ar[d]^-{\simeq}\\
0\ar[r] &\Hom_{B\opp}(\Phi(A),\,\Phi(X)) \ar[r]&
\Hom_{B\opp}(\Phi(T_0),\,\Phi(X))\ar[r]
&\Hom_{B\opp}(\Phi(T_1),\,\Phi(X)) }$$ where the isomorphisms in the
second and third columns are due to $T_0\in\add(T)$ and
$T_1\in\add(T)$, respectively. Consequently, the $\Phi: \Hom_A(X, A)
\lra \Hom_{B\opp}(\Phi(A),\,\Phi(X))$ in the first column is an
isomorphism. This implies that
$$
\Hom_A(-, A)
\lraf{\simeq}\Hom_{B\opp}(\Phi(A),\,\Phi(-)\,)\lraf{\simeq}
\Hom_{B\opp}(T, \Hom_A(-, T)): A\Modcat\to A\opp\Modcat.
$$
Thus $\Hom_A(\cpx{P}, A)\simeq \Hom_{B\opp}(T, \Hom_A(\cpx{P}, T))$
as complexes over $A\opp$. Thus $(3')$ is equivalent to $(3)$.

It remains to show that $(1')$ is equivalent to $(4)$.

For each right $B$-module $Y$, let
$\theta_Y:\Hom_{B\opp}(_AT_B,\,Y)\otimes_AT \lra Y$ be the
evaluation map. Then it follows from the equivalence of $(1)$ and
$(4)$ in Proposition \ref{realization} that $(1')$ is equivalent to
the following statement:

$(4')$ The kernel of the  homomorphism $\partial_0:
\Coker\big(\theta_{\,\Phi(P_0)})\lra
\Coker\big(\theta_{\,\Phi(P_1)}\big)$ induced from the homomorphism
$\Phi(\sigma):\Phi(P_0)\lra \Phi(P_1)$ belongs to $\mathscr{E}$.

Now, we claim that $K\simeq \Ker(\partial_0)$ as right $B$-modules
(see the definition of $K$ in Theorem \ref{main-result} (4)). This
will show that $(1')$ and $(4)$ are equivalent.

To check the above isomorphism,  we first define the following map
for each $A$-module $X$:
$$\zeta_X: \Hom_A(X, A)\otimes_AT \lra \Hom_A(X, T), \;f\otimes
t\mapsto [x\mapsto (x)f\,t]$$ for $f\in\Hom_A(X, A), \,t\in T$ and
$x\in X$. This yields a natural transformation $\zeta: \Hom_A(-,
A)\otimes_AT\lra \Hom_A(-, T)$ from $A\Modcat$ to $B\opp\Modcat$.
Clearly, by definition, we have $\varphi_i=\zeta_{P_i}$ for $i=0,
1$.

Recall that, under the identification of $\Phi(A)$ with $T$ as
$A$-$B$-bimodules, the functor $\Phi$ induces an isomorphism
$\Hom_A(X, A)\lraf{\simeq} \Hom_{B\opp}(T,\,\Phi(X))$ of
$A\opp$-modules. In this sense, one can easily construct the
following commutative diagram:
$$\xymatrix{\Hom_A(X, A)\otimes_AT
\ar[d]^-{\Phi\otimes 1}_-{\simeq}\ar[rr]^-{\zeta_X}&& \Hom_A(X, T)\ar@{=}[d]\\
\Hom_{B\opp}(T,\,\Phi(X))\otimes_AT\ar[rr]^-{\theta_{\,\Phi(X)}}&&\Phi(X)}$$
 This
implies that $\Coker(\zeta_X)$ is naturally isomorphic to
$\Coker\big(\theta_{\,\Phi(X)}\big)$ as $B\opp$-modules. Since
$\varphi_i=\zeta_{P_i}$ for $i=0, 1$, we show that $K\simeq
\Ker(\partial_0)$ as $B\opp$-modules.

Hence, we have proved that the statements $(1)$-$(4)$ in Theorem
\ref{main-result} are equivalent.

Now, suppose $n=2$. Then the complex $\cpx{P}$ is of the following
form:
$$\cdots \lra 0\lra P_2 \lra P_1 \lra P_0 \lra
0\lra \cdots
$$
which is a deleted projective resolution of $_AT$. Since $(1)$ and
$(3)$ in Theorem \ref{main-result} are equivalent, we see that $(1)$
holds if and only if $H^2\big(\Hom_A(\cpx{P},
A)\otimes_AT_B\big)=0$. However, since the tensor functor
$-\otimes_AT_B: A\opp\Modcat\lra B\opp\Modcat$ is always right
exact, we have
$$H^2\big(\Hom_A(\cpx{P}, A)\otimes_AT_B\big)\simeq
H^2\big(\Hom_A(\cpx{P}, A)\big)\otimes_AT\simeq\Ext^2_A(T,
A)\otimes_AT.$$ This finishes the proof of Theorem
\ref{main-result}. $\square$

\medskip
{\it Remarks.}  $(1)$ If the category $\Ker(_AT\otimesL_B-)$ in
Theorem \ref{main-result} is homological in $\D B$, then it follows
from Lemma \ref{rt} (see also Lemma \ref{tilt2} (3)) that the
generalized localization $\lambda: B\to B_T$ of $B$ at the module
$T_B$ exists and is homological, which gives rise to a recollement
of derived module categories:
$$\xymatrix@C=1.2cm{\D{B_T}\ar[r]^-{D(\lambda_*)}
&\D{B}\ar[r]^-{ {}_AT\otimesL_B-}\ar@/^1.2pc/[l]\ar@/_1.3pc/[l]
&\D{A} \ar@/^1.2pc/[l]\ar@/_1.3pc/[l]}$$

\smallskip
$(2)$ Combining the remark following Lemma \ref{prep} with the proof
of Theorem \ref{main-result}, we infer that the complex
$\Hom_A(\cpx{P}, A)\otimes_AT_B$ in Theorem \ref{main-result} is
isomorphic in $\D{B\opp}$ to both $\Hom_A(\cpx{P}, A)\otimesL_AT_B$
and $\rHom_{B\opp}(T,\,B)\otimesL_AT$. This implies that, up to
isomorphism, the cohomology group $H^m\big(\Hom_A(\cpx{P},
A)\otimes_AT_B\big)$ in Theorem \ref{main-result} (3) is independent
of the choice of the projective resolutions of $_AT$ for all $m\in
\mathbb{Z}$.

$(3)$ By the proof of the equivalence of $(1)$ and $(4)$ in Theorem
\ref{main-result}, we know that $\Coker(\zeta_X)\simeq
\Coker\big(\theta_{\,\Phi(X)}\big)$ as $B\opp$-modules for $X\in
A\Modcat$. If $X\in\Add(_AA)$, then $\Phi(X)\in \Prod(T_B)$, and
therefore it follows from Lemma \ref{prep} (1) that
$\Coker\big(\theta_{\,\Phi(X)}\big)$ belongs to $\mathscr{E}:=\{Y\in
{B\opp}\Modcat \mid \Ext^m_{B\opp}(T, Y)=0\mbox{\;for\, all \,}
m\geq0\}$. Particularly, in Theorem \ref{main-result} (4), we always
have $\Coker(\varphi_i)\in\mathscr{E}$ for $i=1, 2$. Note that
$\mathscr{E}$ is closed under kernels of surjective homomorphisms in
$B\opp\Modcat$. Hence, if the homomorphism $\widetilde{\sigma}:
\Coker(\varphi_0)\lra\Coker(\varphi_1)$ induced from $\sigma: P_1\to
P_0$ is surjective, then the kernel $K$ of $\widetilde{\sigma}$ does
belong to $\mathscr{E}$, and therefore the category
$\Ker(T\otimesL_{B}-)$ is homological in $\D{B}$ by the equivalence
of $(1)$ and $(4)$ in Theorem \ref{main-result}.

\medskip
Clearly, the maps $\pi$ and $\omega$ in the definition of tilting
modules induce two canonical quasi-isomorphisms $\widetilde{\pi}:
\cpx{P}\lra T$ and $\widetilde{\omega}: A\lra \cpx{T}$ in $\C{A}$,
respectively. Consequently, both $\widetilde{\pi}$ and
$\widetilde{\omega}$ are isomorphisms in $\D{A}$.

\medskip
As a preparation for the proof of Corollary \ref{cor}, we shall
first establish the following lemma.

\begin{Lem}\label{Z-iso}
The complex $\Hom_A(\cpx{P}, A)$ is isomorphic in $\D{\mathbb{Z}}$
to the following complex:
$$
\Hom_A(T, \cpx{T}):\; \cdots \lra 0\lra \Hom_A(T, T_0)\lra \Hom_A(T,
T_1)\lra \cdots \lra \Hom_A(T, T_n) \lra 0\lra \cdots$$ In
particular, if $A$ is commutative, then $\Hom_A(\cpx{P},
A)\otimes_AT_B\simeq \Hom_A(T, \cpx{T})\otimesL_AT_B$ in
$\D{B\opp}.$

\end{Lem}

{\it Proof.} Since $\widetilde{\pi}$ and $\widetilde{\omega}$ are
chain maps in $\C{A}$, we can obtain  two chain maps in
$\C{\mathbb{Z}}$:
$$\xymatrix{\Hom_A(\cpx{P}, A)\ar[r]^-{(\widetilde{\omega})^*} & \cpx{\Hom}_A(\cpx{P},
\cpx{T}) & \ar[l]_-{(\widetilde{\pi})_*} \Hom_A(T, \cpx{T}).}$$ Now,
we claim that both chain maps are quasi-isomorphisms.

To check this claim, we apply the cohomology functor  $H^i(-)$ to
these chain maps for $i\in\mathbb{Z}$, and construct the following
commutative diagram:
$$\xymatrix{H^i(\Hom_A(\cpx{P}, A))\ar[d]^-{\simeq}\ar[r]^-{H^i((\widetilde{\omega})^*)} & H^i(\cpx{\Hom}_A(\cpx{P},
\cpx{T}))\ar[d]^-{\simeq} & \ar[l]_-{H^i((\widetilde{\pi})_*)}
H^i(\Hom_A(T,\,\cpx{T}))\ar[d]^-{\simeq}\\
\Hom_{\K A}(\cpx{P},
A[i])\ar[d]^-{q_1}\ar[r]^-{(\widetilde{\omega})^*} & \Hom_{\K
A}(\cpx{P}, \cpx{T}[i]) \ar[d]^-{q_2}& \ar[l]_-{(\widetilde{\pi})_*}
\Hom_{\K A}(T, \cpx{T}[i])\ar[d]^-{q_3}\\
\Hom_{\D A}(\cpx{P}, A[i])\ar[r]^-{(\widetilde{\omega})^*}_-{\simeq}
& \Hom_{\D A}(\cpx{P}, \cpx{T}[i])&
\ar[l]_-{(\widetilde{\pi})_*}^-{\simeq} \Hom_{\D A}(T,
\cpx{T}[i])}$$ where the maps $q_j$, for $1\leq j\leq 3$, are
induced by the localization functor $q: \K{A}\to \D{A}$, and where
the isomorphisms in the third row are due to the isomorphisms
$\widetilde{\omega}$ and $\widetilde{\pi}$ in $\D{A}$.

Since $\cpx{P}$ is a bounded complex of projective $A$-modules, both
$q_1$ and $q_2$ are bijective. This implies that
$H^i((\widetilde{\omega})^*)$ is also bijective, and therefore
$(\widetilde{\omega})^*$ is a quasi-isomorphism.

Note that $(\widetilde{\pi})_*$ is a quasi-isomorphism if and only
if  $H^i((\widetilde{\pi})_*)$ is bijective for each
$i\in\mathbb{Z}$. This is also equivalent to saying that $q_3$ is
bijective in the above diagram. Actually, to prove the bijection of
$q_3$, it is enough to show that, for $X\in\add(_AT)$ and
$i\in\mathbb{Z}$, the canonical map $\Hom_{\K A}(T, X[i])\lra
\Hom_{\D A}(T, X[i])$ induced by $q$ is bijective since $\cpx{T}$ is
a bounded complex with each term in $\add(_AT)$. However, this
follows directly from the axiom $(T2)$. Thus $(\widetilde{\pi})_*$
is a quasi-isomorphism.

Consequently, the complexes $\Hom_A(\cpx{P}, A)$ and $\Hom_A(T,
\cpx{T})$ are isomorphic in $\D{\mathbb{Z}}$.

Now, assume that $A$ is commutative. Then each
 $A$-module can be naturally regarded as a right $A$-module and even as an
$A$-$A$-bimodule. In particular, the complex $\cpx{T}$ can be
regarded as a complex of $A$-$A$-bimodules. In this sense,  both
$\widetilde{\pi}: \cpx{P}\lra T$ and $\widetilde{\omega}: A\lra
\cpx{T}$ are quasi-isomorphisms of complexes of $A$-$A$-bimodules.
Moreover, one can check that the chain maps $(\widetilde{\omega})^*$
and $(\widetilde{\pi})_*$ are quasi-isomorphisms in $\C{A\opp}$.
This implies that $\Hom_A(\cpx{P}, A)\simeq  \Hom_A(T, \cpx{T})$ in
$\D{A\opp}$. Note that $\Hom_A(\cpx{P}, A)\otimes_AT_B\simeq
\Hom_A(\cpx{P}, A)\otimesL_AT_B$ in $\D{B\opp}$ (see the above
remark (2)). As a result, we have $\Hom_A(\cpx{P},
A)\otimes_AT_B\simeq \Hom_A(T, \cpx{T})\otimesL_AT_B$ in
$\D{B\opp}.$ $\square$

\smallskip

{\bf Proof of Corollary \ref{cor}.} $(1)$ By the remark (3) at the
end of the proof of Theorem \ref{main-result}, we know that if the
homomorphism $\widetilde{\sigma}:
\Coker(\varphi_0)\lra\Coker(\varphi_1)$ induced from $\sigma: P_1\to
P_0$ (see Theorem \ref{main-result} (4)) is surjective, then
$\Ker(_AT\otimesL_B-)$ is homological in $\D{B}$.

Now, we verify this sufficient condition for the good tilting module
$_AT$ which satisfies the assumption in $(1)$.

In fact, by assumption, we can assume that $_AM$ has a projective
resolution: $0\lra P_1'\lraf{\sigma'} P_0'\lra {}_AM\lra 0$ with
$P_0', P_1'\in\Add(_AA)$, and that $_AN$ has a projective
presentation: $P_1''\lra P_0''\lraf{\sigma''} {}_AN\lra 0$ with
$P_0''\in\Add(_AA)$ and $P_1''\in\add(_AA)$. Since $_AT=M\oplus N$,
we can choose $\sigma=\left(\begin{array}{cc} \sigma' & 0\\
0 & \sigma''\end{array}\right): P_1'\oplus P_1''\lra P_0'\oplus
P_0''$. Recall that $\zeta: \Hom_A(-, A)\otimes_AT\lra \Hom_A(-, T)$
is a natural transformation from $A\Modcat$ to $B\opp\Modcat$ (see
the proof of Theorem \ref{main-result}). Certainly, if
$X\in\add(_AA)$, then $\zeta_X$ is an isomorphism, and so
$\Coker(\zeta_X)=0$.

Let $\widetilde{\sigma}':
\Coker(\zeta_{P_0'})\lra\Coker(\zeta_{P_1'})$ and
$\widetilde{\sigma}'':
\Coker(\zeta_{P_0''})\lra\Coker(\zeta_{P_1''})$ be the homomorphisms
induced from $\sigma'$ and $\sigma''$, respectively. By definition,
we have $\varphi_i=\zeta_{P_i}$ for $i=0, 1$, and
$$\widetilde{\sigma}=\left(\begin{array}{cc} \widetilde{\sigma}' & 0\\
0 & \widetilde{\sigma}''\end{array}\right):
\Coker(\zeta_{P_0'})\oplus \Coker(\zeta_{P_0''})\lra
\Coker(\zeta_{P_1'})\oplus \Coker(\zeta_{P_1''}).$$

Now, we show that $\widetilde{\sigma}$ is surjective, or
equivalently, both $\widetilde{\sigma}'$ and $\widetilde{\sigma}''$
are surjective. In fact, since $P_1''\in\add(_AA)$, we see that
$\Coker(\zeta_{P_1''})=0$. Thus $\widetilde{\sigma}''$ is
surjective. As $_AM$ is a direct summand of $_AT$ and of projective
dimension at most $1$, it follows from the axiom $(T2)$ that the map
$\Hom_A(\sigma', T): \Hom_A(P_0', T)\lra \Hom_A(P_1', T)$ is
surjective. This implies that $\widetilde{\sigma}'$ is a surjection.
Consequently, $\widetilde{\sigma}$ is surjective. Thus
$\Ker(_AT\otimesL_B-)$ is homological in $\D{B}$. This finishes the
proof of $(1)$.

$(2)$ Suppose that $\Ker(_AT\otimesL_B-)$ in Theorem
\ref{main-result} is homological. By Theorem \ref{main-result}, we
have $H^m\big(\Hom_A(\cpx{P}, A)\otimes_AT_B\big)=0$ for all $m\geq
2.$ In the sequel, we shall show that if $H^n\big(\Hom_A(\cpx{P},
A)\otimes_AT_B\big)=0$, then $T_n=0.$

In fact, since $A$ is commutative, it follows from the proof of
Lemma \ref{Z-iso} that $\Hom_A(\cpx{P}, A)\simeq \Hom_A(T, \cpx{T})$
in $\D{A\opp}$. Note that the tensor functor $-\otimes_AT_B:
A\opp\Modcat\lra B\opp\Modcat$ is right exact. This means that
$$0=H^n\big(\Hom_A(\cpx{P}, A)\otimes_AT_B\big)\simeq
H^n(\Hom_A(\cpx{P}, A))\otimes_AT\simeq H^n(\Hom_A(T,
\cpx{T}))\otimes_AT.$$ In particular, we have $H^n(\Hom_A(T_n,
\cpx{T}))\otimes_AT_n=0$, due to $T_n\in\add(_AT)$.

Recall that the complex $\Hom_A(T_n, \cpx{T})$ is of the form
$$
\cdots \lra 0\lra \Hom_A(T_n, T_0)\lra \cdots \lra \Hom_A(T_n,
T_{n-1})\lra \Hom_A(T_n, T_n) \lra 0\lra \cdots$$ As $\Hom_A(T_n,
T_{n-1})=0$ by our assumption in Corollary \ref{cor} (2), we obtain
$H^n(\Hom_A(T_n, \cpx{T}))=\Hom_A(T_n, T_n)$. Thus
$\End_A(T_n)\otimes_AT_n=0$. It follows from the surjective map
$$\End_A(T_n)\otimes_AT_n \lra T_n,\; f\otimes x\mapsto (x)f \;\,
\mbox{ for } f\in\End_A(T_n) \;\mbox{ and }\; x\in T_n$$ that
$T_n=0$. This finishes the proof of the above claim.

By our assumption, we have $\Hom_A(T_{i+1}, T_i)=0$ for $1\leq i
\leq n-1$.  Now, we can proceed by induction on $n$ to show that
$T_j=0$ for $2\leq j\leq n$. Thus, by Lemma \ref{tilt1} (4), $T$ is
a $1$-tilting module, that is, the projective dimension of $_AT$ is
at most $1$.

The sufficiency of Corollary \ref{cor} (2) follows from Theorem
\ref{main-result}, see also \cite[Theorem 1.1 (1)]{CX1}. This
finishes the proof of Corollary \ref{cor}. $\square$

\smallskip
Let us end this section by constructing an example of infinitely
generated $n$-tilting modules $T$ such that $\Ker(T\otimesL_B-)$ are
homological.

Let $A$ be an arbitrary ring with  a classical $n$-tilting
$A$-module $T'$.  Suppose $_AT' = M\oplus N$ with $M$ a nonzero
$A$-module of projective dimension at most $1$. Let $I$ be an
infinite set, and let $T:=M^{(I)}\oplus N$. Then $T$ is a good
$n$-tilting module. Since $T$ satisfies Corollary \ref{cor} (1), we
see that $\Ker(T\otimesL_B-)$ is homological in $\D B$.

\section{Applications to cotilting modules\label{sect6}}

Our main purpose in this section is to show Theorem \ref{coth} and
develop some conditions which can be used to decide if subcategories
induced from cotilting modules are homological or not. We also
provide an example to show that recollements provided by cotilting
modules depend upon the choice of injective cogenerators.

\subsection{Proof of Theorem \ref{coth}}\label{6.1}

In this section, we shall apply the results in Section \ref{sect4}
to deal with cotilting modules. First, we shall construct Ringel
modules from good cotilting modules, and then use Proposition
\ref{realization} to show the main result, Corollary
\ref{real-cotilt}, of this section, and finally give the proof of
Theorem \ref{coth}.

Suppose that $A$ is a ring and that $W$ is a fixed injective
cogenerator for $A\Modcat$. Recall that an $A$-module $W$ is called
a \emph{cogenerator} for $A\Modcat$ if, for any $A$-module $Y$,
there exists an injective homomorphism $Y\to W^I$ in $A\Modcat$ with
$I$ a set. This is also equivalent to saying that, for any non-zero
homomorphism $f: X\to Y$ in $A\Modcat$, there exists a homomorphism
$g\in\Hom_A(Y, W)$ such that $fg$ is non-zero.

Let us recall the definition of $n$-cotilting modules for $n$ a
natural number.

\begin{Def} {\rm \label{cotilting}
An  $A$-module $U$ is called an $n$-\emph{cotilting module} if the
following three conditions are satisfied:

$(C1)$ there exists an exact sequence
$$
0\lra U\lra I_0 \lraf{\delta} I_1\lra \cdots \lra I_n\lra 0$$ of
$A$-modules such that $I_i$ is an injective module for every $0\leq
i\leq n$;

$(C2)$ $\Ext^j_A(U^{I}, U)=0$ for each $j\geq 1$ and for every
nonempty set $I$; and

$(C3)$ there exists an exact sequence
$$
0\lra U_n \lra \cdots \lra U_1 \lra U_0\lra  W\lra 0$$ of
$A$-modules, such that $U_i\in\Prod(_AU)$ for all $0\leq i\leq n$.

\smallskip

An $n$-cotilting $A$-module $U$ is said to be \emph{good} if it
satisfies ($C1$), ($C2$) and

$(C3)'$ there is an exact sequence
$$0\lra U_n \lra \cdots \lra U_1 \lra U_0\lra  W\lra 0$$ of $A$-modules, such that  $U_i\in\add(_AU)$ for
all $0\leq i\leq n$.

We say that $U$ is a (good) cotilting $A$-module if $_AU$ is (good)
$n$-cotilting for some $n\in\mathbb{N}$.}
\end{Def}

We remark that if both $W_1$ and $W_2$ are injective cogenerators
for $A\Modcat$, then $\Prod(W_1)=\Prod(W_2)$. This implies that the
definition of cotilting modules is independent of the choice of
injective cogenerators for $A\Modcat$. However, the definition of
good cotilting modules relies on the choice of injective
cogenerators for $A\Modcat$.

As in the case of tilting modules, for a given $n$-cotilting
$A$-module $U$ with $(C1)$-$(C3)$, the $A$-module
$U':=\bigoplus_{i=0}^n U_i$ is a good $n$-cotilting module which is
equivalent to the given one in the sense that $\Prod(U)=\Prod(U')$.

\medskip
From now on, we assume that $U$ is a {\bf{good}} $n$-cotilting
$A$-module with $(C_1), (C_2)$ and $(C_3)'$, where the module $W$ in
$(C_3)'$ is referred to the fixed injective cogenerator for $A$-Mod.
In this event, we shall call $U$ a \emph{good $n$-cotilting
$A$-module with respect to $W$}.

\smallskip
Let $R:=\End_A(U)$, $M:=\Hom_A(U, W)$ and $\Lambda:=\End_A(W)$. Then
$M$ is an $R$-$\Lambda$-bimodule.

\smallskip
First of all, we collect some basic properties of good cotilting
modules in the following lemma.

\begin{Lem}\label{cotilt1}
The following hold for the cotilting module $U$.

$(1)$ The  $R$-module $M$ has a finitely generated projective
resolution of length at most $n$:
$$
0\lra \Hom_A(U, U_n)\lra \cdots \lra\Hom_A(U, U_1)\lra\Hom_A(U,
U_0)\lra M\lra 0$$ such that $U_m\in\add(_AU)$ for all $0\leq m\leq
n$.

$(2)$ The Hom-functor $\Hom_A(U,-): A\Modcat\ra R\Modcat$ induces an
isomorphism of rings: $\Lambda \simeq \End_R(M)$,  and $\Ext^i_R(M,
M)=0$ for all $ i\geq 1$.

$(3)$ The module $M$ is an $n$-Ringel $R$-module.
\end{Lem}
{\it Proof.} $(1)$ Applying the functor $\Hom_A(U,-)$ to the
sequence
$$0\lra U_n \lra \cdots \lra U_1 \lra U_0\lra  W\lra 0$$
in the axiom $(C_3)'$, we obtain the sequence in $(1)$ with all
$\Hom_A(U, U_i)\in \add(_RR)$. The exactness of this sequence
follows directly from the axiom $(C2)$. This also implies that the
projective dimension of $_RM$ is at most $n$.

$(2)$ Denote by $\Psi$ the Hom-functor $\Hom_A(U, -): A\Modcat\to
R\Modcat$. Then $\Psi(U)=R$, \,$\Psi(W)=M$ and, for every
$X\in\add(_AU)$, we have
$$\Hom_A(X,
W)\lraf{\simeq}\Hom_R(\Psi(X), \Psi(W)).$$

Clearly, if $n=0$, then $W=U_0$, $M=\Hom_A(U, U_0)$ as $R$-modules.
In this case, one can easily check $(2)$.

Suppose $n\geq 1$. By $(1)$, the $R$-module $M=\Psi(W)$ has a
finitely generated projective resolution
$$
0\lra\Psi(U_n)\lra \cdots \lra\Psi(U_1)\lra\Psi(U_0)\lra \Psi(W)\lra
0$$ with  $U_m\in\add(U)$ for all $0\leq m\leq n$. Applying the
functor $\Hom_A(-, W)$ to the resolution of $W$ in $(C3)'$, we can
construct the following commutative diagram: {\footnotesize
$$\xymatrix{
0\ar[r] & \Hom_A(W, W)\ar[r]\ar[d]_-{\Psi} &\Hom_A(U_0, W)
\ar[r]\ar[d]_-{\simeq} &\Hom_A(U_1, W) \ar[r]\ar[d]_-{\simeq}
&\cdots \ar[r] & \Hom_A(U_n, W)
\ar[r]\ar[d]_-{\simeq} &0\\
0\ar[r] & \Hom_R(\Psi(W), \Psi(W))\ar[r] &\Hom_R(\Psi(U_0), \Psi(W))
\ar[r]&\Hom_R(\Psi(U_1), \Psi(W)) \ar[r] & \cdots \ar[r] &
\Hom_R(\Psi(U_n), \Psi(W)) \ar[r]&0}
$$}

\noindent where the isomorphisms in the diagram are due to
$U_m\in\add(_AU)$ for $m\le n$. Since $_AW$ is injective, the first
row in the diagram is exact. Note that the following sequence
$$0\lra \Hom_R(\Psi(W),
\Psi(W))\lra \Hom_R(\Psi(U_0), \Psi(W)) \lra \Hom_R(\Psi(U_1),
\Psi(W))$$ is always exact since $\Psi(U_1)\lra\Psi(U_0)\lra
\Psi(W)\lra 0$ is exact in $R\Modcat$. This implies that the map
$\Psi: \End_A(W)\lra \End_R(\Psi(W))$ is an isomorphism of rings and
that the second row in the diagram is also exact. Thus $\Ext^i_R(M,
M)=\Ext^i_R(\Psi(W), \Psi(W))=0$ for all $ i\geq 1$.

$(3)$ We check the axioms $(R1)$-$(R3)$ in Definition \ref{rm} for
$M$. Clearly, the axioms $(R1)$ and $(R2)$ follow from $(1)$ and
$(2)$, respectively. It remains to show the axiom $(R3)$ for $M$. In
fact, by the axiom $(C1)$, there exists an exact sequence of
$A$-modules:
$$
0\lra U\lra I_0 \lra I_1\lra \cdots \lra I_n\lra 0$$ where  $I_i$ is
an injective module for $0\leq i\leq n$. Since $W$ is an injective
cogenerator for $A\Modcat$, we have $I_i\in\Prod(_AW)$. Moreover,
from the axiom $(C2)$, we see that $\Ext^j_A(U, U)=0$ for all $j\geq
1$. This implies that the following sequence
$$
0\lra R\lra \Hom_A(U, I_0) \lra \Hom_A(U, I_1)\lra \cdots \lra
\Hom_A(U, I_n)\lra 0$$ is exact. Since the functor $\Hom_A(U, -)$
commutes with arbitrary direct products, it follows from
$I_i\in\Prod(_AW)$ that $\Hom_A(U, I_i)\in\Prod\big({}_R\Hom_A(U,
W)\big)=\Prod(_RM)$. This shows that $_RM$ satisfies the axiom
$(R3)$. Therefore $M$ is an $n$-Ringel $R$-module. $\square$

\medskip
Observe that, by Lemma \ref{cotilt1} (2), the ring $\End_R(M)$ can
be naturally identified with $\Lambda$ (up to isomorphism of rings).
Now, we define
$${\bf G}:= {}_RM\otimesL_\Lambda-:
\;\D{\Lambda}\lra\D{R}\quad \mbox{and}\quad  {\bf
H}:=\rHom_R(M,-):\; \D{R}\lra\D{\Lambda}.$$ Since $_RM$  is a Ringel
$R$-module satisfying both $(R1)$ and $(R2)$ in Definition \ref{rm},
it follows from Lemma \ref{rt} that there exists a recollement of
triangulated categories:
$$
\xymatrix@C=1.2cm{\Ker({\bf H})\ar[r]^-{{i_*}} &\D{R}\ar[r]^-{{\bf
H}}\ar@/^1.2pc/[l]\ar_-{i^*}@/_1.2pc/[l] &\D{\Lambda}
\ar@/^1.2pc/[l]\ar@/_1.2pc/[l]_{{\bf G}} }$$

\medskip
\noindent where $(i^*, i_*)$ is a pair of adjoint functors with
$i_*$ the inclusion.

If $\Ker({\bf H})$ is homological, then it follows from Lemma
\ref{rt} that the generalized  localization $\lambda: R\to R_M$ of
$R$ at $M$ exists and induces a recollement of derived module
categories:
$$(\ddag)\quad
\xymatrix@C=1.2cm{\D{R_M}\ar[r]^-{D(\lambda_*)} &\D{R}\ar[r]^-{{\bf
H}}\ar@/^1.3pc/[l]\ar@/_1.2pc/[l] &\D{\Lambda}
\ar@/^1.2pc/[l]\ar@/_1.2pc/[l]_{{\bf G}}}$$

\medskip \noindent Thus we may construct recollements of derived
module categories from good cotilting modules. Here, a problem
arises naturally:

\medskip
{\bf Problem:} When is $\Ker({\bf H})$ homological in $\D{R}$?

\medskip
This seems to be a difficult problem because we cannot directly
apply Proposition \ref{realization} to the Ringel module $_RM$. The
reason is that we do not know whether $_RM$ is \emph{good}.
Actually, we do not know whether the right $\Lambda$-module $M$ is
strongly $\Lambda$-Mittag-Leffler. Certainly, if $\Lambda$ is right
noetherian, then $M$ is a perfect Ringel $R$-module (see Definition
\ref{rm}), and must be good.

Though we cannot solve this problem entirely, we do have some
partial solutions to the problem.

\begin{Koro}\label{real-cotilt}
Suppose that  $A$ is a ring together with an injective cogenerator
$W$ for $A\Modcat$. Let $U$ be a good $n$-cotilting $A$-module with
respect to $W$. Suppose that $\Lambda:=\End_A(W)$ is a right
noetherian ring. Then the following are equivalent:

$(a)$ $\Ker({\bf H})$ is homological in $\D{R}$.

$(b)$ $H^m\big({_R}\Hom_A(U,
W)\otimes_{\Lambda}\Hom_A(W,\,\cpx{I})\big)=0$ for all $m\geq 2$,
where $\cpx{I}$ is a deleted injective coresolution of $_AU$:
$$\cdots \lra 0\lra I_0 \lraf{\delta} I_1\lra \cdots \lra I_n\lra 0\lra
\cdots
$$
with $I_i$ in degree $i$ for all $0\leq i\leq n$.

$(c)$ The kernel $K$ of the homomorphism
$\Coker(\phi_0)\lra\Coker(\phi_1)$ induced from the map $\delta:
I_0\to I_1$ satisfies  $\,\Ext^m_{R}(M, K)=0$ for all $ m\geq 0$,
where $\phi_i:\Hom_A(U, W)\otimes_{\Lambda}\Hom_A(W,\,I_i)\lra
\Hom_A(U, I_i)$ is the composition map for $i=0, 1$.

\end{Koro}

 {\it Proof.} By the proof of Lemma \ref{cotilt1} (3),
the module $M:=\Hom_A(U, W)$ is an $n$-Ringel $R$-module. Moreover,
the sequence in the axiom $(R3)$ can be chosen as follows:
$$
0\lra R\lra \Hom_A(U, I_0) \lra \Hom_A(U, I_1)\lra \cdots \lra
\Hom_A(U, I_n)\lra 0.$$  In this case, the complex $\cpx{M}$ can be
defined as the following complex:
$$
\Hom_A(U, \cpx{I}):\; 0\lra \Hom_A(U, I_0) \lra \Hom_A(U, I_1)\lra
\cdots \lra \Hom_A(U, I_n)\lra 0.$$ Under the assumption that
$\Lambda$ is right noetherian, we know that $M$ is a good Ringel
$R$-module. So it follows from Proposition \ref{realization} that
$(a)$ is equivalent to the following:

$(b')$ $H^j\big({_R}M\otimes_{\Lambda}\Hom_R(M,\,\cpx{M})\big)=0$
for any $j\geq 2$, where $\cpx{M}:=\Hom_A(U, \cpx{I})$.

To prove that $(a)$ and $(b)$ in Corollary \ref{real-cotilt} are
equivalent, it is sufficient to show that $(b')$ and $(b)$ are
equivalent. For this purpose, we shall show that
$\Hom_R(M,\,\cpx{M})\simeq\Hom_A(W, \cpx{I})$ as complexes over
$\Lambda$.

Let $\Psi=\Hom_A(U, -): A\Modcat\to R\Modcat$. Then $\Psi(W)=M$ and
$\cpx{M}=\Psi(\cpx{I})$. Clearly, the functor $\Psi$ induces a
natural transformation
$$
\Hom_A(W,-) \lra \Hom_R(\Psi(W), \Psi(-)) : A\Modcat \lra
\Lambda\Modcat.
$$
This yields a chain map from $\Hom_A(W, \cpx{I})\lra \Hom_R(\Psi(W),
\Psi(\cpx{I}))=\Hom_R(M, \cpx{M})$ in $\C \Lambda$, that is,
$$\xymatrix{
0\ar[r] &\Hom_A(W, I_0) \ar[r]\ar[d] &\Hom_A(W, I_1)
\ar[r]\ar[d]&\cdots \ar[r] & \Hom_A(W, I_n)
\ar[r]\ar[d] &0\\
0\ar[r] &\Hom_R(\Psi(W), \Psi(I_0)) \ar[r]&\Hom_R(\Psi(W),
\Psi(I_1)) \ar[r] & \cdots \ar[r] & \Hom_R(\Psi(W), \Psi(I_n))
\ar[r]&0}
$$
Note that all $I_i$ are injective $A$-modules. To verify that this
chain map is an isomorphism of complexes, it is enough to show that
$\Psi$ induces an isomorphism of $\Lambda$-modules:
$$
\Hom_A(W,X) \lraf{\simeq} \Hom_R(\Psi(W), \Psi(X))
$$
for any injective $A$-module $X$. In the following, we shall prove
that this holds even for any $A$-module $X$.

Suppose $n=0$. By the axiom $(C3)'$, we know that $W = U_0$ as
$A$-modules with $U_0\in\add(_AU)$. It is clear that $\Hom_A(U_0, X)
\lraf{\simeq} \Hom_R(\Psi(U_0), \Psi(X))$ since $U_0\in\add(_AU)$.
Thus $\Hom_A(W,X) \lraf{\simeq} \Hom_R(\Psi(W), \Psi(X)).$

Now, suppose $n\geq 1$. By the axiom $(C3)'$ and Lemma \ref{cotilt1}
(1), there exists an exact sequence $U_1 \lra U_0\lra W\lra 0$ of
$A$-modules with $U_0, U_1\in\add(_AU)$ such that
$\Psi(U_1)\lra\Psi(U_0)\lra \Psi(W)\lra 0$ is also exact in
$R\Modcat$. From this sequence,  we may construct the following
exact commutative diagram:
$$\xymatrix{
0\ar[r] & \Hom_A(W, X)\ar[r]\ar[d]^-{\Psi} & \Hom_A(U_0, X)
\ar[r]\ar[d]^-{\simeq}&\Hom_A(U_1, X)\ar[d]^-{\simeq}\\
0\ar[r] & \Hom_R(\Psi(W), \Psi(X))  \ar[r]& \Hom_R(\Psi(U_0),
\Psi(X)) \ar[r] &\Hom_R(\Psi(U_1), \Psi(X)) }$$ where the last two
vertical maps are isomorphisms  since $U_0, U_1\in\add(_AU)$. This
means that $\Hom_A(W,X) \lraf{\simeq} \Hom_R(\Psi(W), \Psi(X))$ for
every $A$-module $X$.

Consequently, we see that $\Hom_A(W,
\cpx{I})\simeq\Hom_R(M,\,\cpx{M})$ as complexes over $\Lambda$. Thus
$(b')$ and $(b)$, and therefore, also $(a)$ and $(b)$, are
equivalent.

Note that if we identify $\Hom_R(M,\,\cpx{M})$ with $\Hom_A(W,
\cpx{I})$ as complexes over $\Lambda$, then the equivalence of $(a)$
and $(c)$ in Corollary \ref{real-cotilt} can be concluded from that
of $(1)$ and $(4)$ in Proposition \ref{realization}. Here, we leave
the details to the reader. $\square$

\medskip
As a consequence of  Corollary \ref{real-cotilt} (see also Corollary
\ref{appl}), we have the following result.

\begin{Koro}\label{cotilt2}
Let $U$ be a good $n$-cotilting $A$-module with respect to the
injective cogenerator $_AW$. Suppose that $\Lambda:=\End_A(W)$ is a
right noetherian ring.

$(1)$ If $_AU=M\oplus N$ such that $_AM$ has injective dimension at
most $1$ and that $_AN$ has an injective copresentation $0\lra
{}_AN\lra E_0\lra E_1$ with $E_1\in\add(_AW)$, then $\Ker({\bf H})$
is homological in $\D{R}$.

$(2)$ If $n=2$, then $\Ker({\bf H})$ is homological in $\D R$ if and
only if $\Hom_A(U, W)\otimes_{\Lambda}\Ext^2_A(W, U)=0.$
\end{Koro}

{\it Proof.} The idea of the proof of $(1)$ is very similar to that
of Corollary \ref{cor} (1). Here, we just give a sketch of the
proof.

Note that $\mathscr{E}:=\{Y\in R\Modcat \mid \Ext^m_R(M,
Y)=0\mbox{\;for\, all \,} m\geq0\}$ is closed under kernels of
surjective homomorphisms in $R\Modcat$, and that $\Coker(\phi_0)$
and $\Coker(\phi_1)$ (see Corollary \ref{real-cotilt} (c)) always
belong to $\mathscr{E}$ by Lemma \ref{prep} (1). Thus, according to
the equivalence of $(a)$ and $(c)$ in Corollary \ref{real-cotilt},
if we want to show $(1)$, then it suffices to verify that the
homomorphism $\widetilde{\delta}: \Coker(\phi_0)\lra\Coker(\phi_1)$
induced from $\delta: I_0\to I_1$ is surjective. Actually, this is
guaranteed by the assumption that the injective dimension of $_AM$
is at most $1$ and $E_1\in\add(_AW)$. For more details, we refer the
reader to the proof of Corollary \ref{cor} (1).

As to $(2)$, we keep the notation in the proof of Corollary
\ref{real-cotilt}. Suppose $n=2$. Then the complex $\cpx{I}$ in
Corollary \ref{real-cotilt} (b) has the following form
$$\cdots \lra 0\lra I_0 \lra I_1\lra I_2\lra 0\lra
\cdots.
$$
 By Corollary
\ref{real-cotilt}, the category $\Ker({\bf H})$ is homological if
and only if
$H^2\big({_R}M\otimes_{\Lambda}\Hom_A(W,\,\cpx{I})\big)=0$, where
$M:=\Hom_A(U, W)$. Note that the tensor functor
$_RM\otimes_{\Lambda}-: \Lambda\Modcat\lra R\Modcat$ is right exact.
Consequently, we have
$$H^2\big({_R}M\otimes_{\Lambda}\Hom_A(W,\,\cpx{I})\big)\simeq
M\otimes_{\Lambda} H^2(\Hom_A(W,\,\cpx{I}))\simeq
M\otimes_{\Lambda}\Ext^2_A(W, U).$$ This shows $(2)$. $\square$

\medskip
Finally, we point out a special case for which the ring $\Lambda$ in
Corollary \ref{real-cotilt} is right noetherian.

Let $k$ be a commutative Artin ring. Let $\rad(k)$ be the radical of
$k$ (that is, the intersection of all maximal ideals of $k$), and
let $J$ be the injective envelope of $k/\rad(k)$. We say that a
$k$-algebra $A$ is an \emph{Artin $k$-algebra}, or Artin algebra for
short, if $A$ is finitely generated as a $k$-module.

Suppose that $A$ is an Artin $k$-algebra. It is well known that the
functor $\Hom_k(-, J)$ is a duality between the category $A\modcat$
of finitely generated $A$-modules and that of finitely generated
$A\opp$-modules. In particular, the dual module $\Hom_k(A_A, J)$ of
the right $A$-module $A_A$ is an injective cogenerator for
$A\modcat$, or even for $A\Modcat$. In this case, we shall call
$\Hom_k(A_A, J)$ the \emph{ordinary injective cogenerator} for
$A\Modcat$.

Note that $\End_A(\Hom_k(A_A, J))\simeq \End_{A\opp}(A)\opp\simeq A$
as rings. So, if the module $W$ in Corollary \ref{real-cotilt} is
chosen to be the module $\Hom_k(A_A, J)$, then the ring
$\Lambda:=\End_A(W)$ is isomorphic to $A$. Since $A$ is an Artin
algebra, it is a left and right Artin ring, and certainly a right
noetherian ring. Thus $\Lambda$ is right noetherian and always
satisfies the assumption in Corollary \ref{real-cotilt}.

\medskip
{\bf Proof of Theorem \ref{coth}.} Recall that $_AW$ is the ordinary
injective cogenerator over the Artin algebra $A$. According to the
above-mentioned facts, the ring $\Lambda:=\End_A(W)$ is isomorphic
to $A$, and therefore right noetherian. Since $_AU$ is a good
$1$-cotilting module with respect to $W$, we know from Corollary
\ref{cotilt2} (1) that the category $\Ker({\bf H})$ is homological.
Now, Theorem \ref{coth} follows from the diagram $(\ddag)$ above
Corollary \ref{real-cotilt}. $\square$

\medskip
Let us end this section by a couple of remarks related to the
results in this section.

{\it Remarks.} $(1)$ If $A$ is a commutative ring and $W$ is an
injective cogenerator for $A\Modcat$, then the dual module
$\Hom_A(T, W)$ of a tilting $A$-module $T$ is always a cotilting
$A$-module. However, there exist cotilting modules over Pr\"ufer
domains, which are not equivalent to the dual modules of any tilting
modules (see \cite[Chapter 11, Section 4.16]{HHK}). This means that
the investigation of infinitely generated cotilting modules cannot
be carried out by using dual arguments of infinitely generated
tilting modules.

$(2)$ Corollary \ref{real-cotilt} provides actually a recollement of
$\D{\End_A(U)}$  with $\D{R_M}$ on the left-hand side and
$\D{\Lambda}$ on the right-hand side (see $(\ddag)$ for notation).
This recollement depends upon the choice of injective cogenerators
for $A\Modcat$. That is, for a fixed cotilting module $_AU$, if
different injective cogenerators $W$ for $A\Modcat$ are chosen in
the axiom $(C3)'$, then one may get completely different
recollements of $\D{\End_A(U)}$.

For example, let ${\mathbb Q}_{(p)}$, $\mathbb Q$, ${\mathbb Z}_p$
and ${\mathbb Q}_p$ denote the rings of $p$-integers, rational
numbers, $p$-adic integers and $p$-adic numbers, respectively.
Recall that ${\mathbb Q}_{(p)}$ is the localization of $\mathbb{Z}$
at the prime ideal $p\mathbb{Z}$. In particular, it is a local
Dedekind domain. Moreover, let $E(\mathbb{Z}/p\mathbb{Z})$ be the
injective envelope of $\mathbb{Z}/p\mathbb{Z}$, which is an
injective cogenerator for the category of ${\mathbb
Q}_{(p)}$-modules.

Now, we take $A:=\mathbb{Q}_{(p)}$, $T:=\mathbb{Q}\oplus
E(\mathbb{Z}/p\mathbb{Z})$ and
$U:=\Hom_A(T,E(\mathbb{Z}/p\mathbb{Z}))$.  Due to \cite[Section
7.1]{CX1}, we have

$(a)$ the module $T$ is a Bass $1$-tilting module over $A$, and
therefore $U$ is an $1$-cotilting $A$-module.

$(b)$ $\End_A(E(\mathbb{Z}/p\mathbb{Z}))\simeq {\mathbb Z}_p\,$ and
$\,\Hom_A(\mathbb{Q},
E(\mathbb{Z}/p\mathbb{Z}))\simeq\mathbb{Q}\otimes_A\End_A(E(\mathbb{Z}/p\mathbb{Z}))\simeq
\mathbb{Q}\otimes_A {\mathbb Z}_p\simeq {\mathbb Q}_p $. Thus
$U\simeq {\mathbb Z}_p\oplus {\mathbb Q}_p$ as $A$-modules.

$(c)$ By \cite[Lemma 6.5(3)]{CX1}, there exists an exact sequence of
${\mathbb Z}_p$-modules (and also $A$-modules):
$$
(*')\quad 0\lra {\mathbb Z}_p \lraf{\varphi} {\mathbb Q}_p\lra
E(\mathbb{Z}/p\mathbb{Z})\lra 0.
$$

Note that ${\mathbb Q}_p$ is an injective and flat $A$-module and
that $(*')$ is an injective coresolution of ${\mathbb Z}_p$ as an
$A$-module. This also implies that $W:={\mathbb Q}_p\oplus
E(\mathbb{Z}/p\mathbb{Z})$ is an injective cogenerator for
$A\Modcat$.

On the one hand, we may consider $U$ as a good $1$-cotilting
$A$-module with respect to $W$. Applying $\Hom_A(U, -)$ to the
sequence $(*')$, we get a projective resolution of $\Hom_A(U,
E(\mathbb{Z}/p\mathbb{Z}))$ as an $\End_A(U)$-module:
$$
 \quad 0\lra \Hom_A(U, {\mathbb Z}_p) \lraf{\varphi^*} \Hom_A(U, {\mathbb Q}_p)\lra
\Hom_A(U, E(\mathbb{Z}/p\mathbb{Z}))\lra 0.
$$
Since both ${\mathbb Q}_p$ and $E(\mathbb{Z}/p\mathbb{Z})$ belong to
$\add(_AW)$, one can use Lemma \ref{cotilt1} to show that $\Hom_A(U,
W)$ is a classical $1$-tilting $\End_A(U)$-module such that
$\End_{\End_A(U)}(\Hom_A(U, W))\simeq \End_A(W)$ as rings. It
follows that $\End_A(U)$ and  $\End_A(W)$ are derived equivalent. In
this case, we get a trivial recollement:
$\D{\End_A(U)}\lraf{\simeq}\D{\Lambda}$ with $\Lambda:=\End_A(W)$.
Note that this derived equivalence can also be seen from
\cite[Theorem 1.1]{hx2}.

On the other hand, we consider $U$ as a good $1$-cotilting
$A$-module with respect to $W':= E(\mathbb{Z}/p\mathbb{Z})$.
Clearly, the sequence $(*')$ can paly the role in the axiom $(C3)'$.
Since $\End_A(E(\mathbb{Z}/p\mathbb{Z})) \simeq\mathbb{Z}_p,$
we know from \cite[Corollary 2.5.16]{EJ} that
$\End_A(E(\mathbb{Z}/p\mathbb{Z}))$ is a noetherian ring. This
implies that $U$ satisfies the assumptions in Corollary
\ref{cotilt2} (1).

By \cite[Theorem 3.4.1]{EJ}, one can check that
$$
\End_A(\mathbb{Z}_p)\simeq \mathbb{Z}_p, \;\, \Hom_A(\mathbb{Q}_p,
\mathbb{Z}_p)=0=\Ext_A^1(\mathbb{Q}_p,
\mathbb{Z}_p)=\Hom_A(E(\mathbb{Z}/p\mathbb{Z}), \mathbb{Q}_p),$$ and
further that
$$
\End_A(U)\simeq \left(\begin{array}{lc} \mathbb{Z}_p & \End_A(\mathbb{Q}_p)\\
0 &  \End_A(\mathbb{Q}_p)\end{array}\right) \;\, \mbox{and} \;\,
\End_A(W)\simeq \left(\begin{array}{lc} \End_A(\mathbb{Q}_p)& \End_A(\mathbb{Q}_p)\\
\quad 0 & \mathbb{Z}_p\end{array}\right).
$$
Moreover, the universal localization of $\End_A(U)$ at the map
$\varphi^*$, or at the module $\Hom_A(U,
E(\mathbb{Z}/p\mathbb{Z}))$, is isomorphic to
$M_2(\End_A(\mathbb{Q}_p))$, the $2\times 2$ matrix ring over
$\End_A(\mathbb{Q}_p)$.

Now, we can construct the following non-trivial recollement of
derived module categories from the cotilting module $U$ with respect
to $W'= E(\mathbb{Z}/p\mathbb{Z})$:

$$\xymatrix@C=1.2cm{\D{\End_A(\mathbb{Q}_p)}\ar[r]
&\D{\End_A(U)}\ar[r]\ar@/^1.3pc/[l]\ar@/_1.2pc/[l] &\D{\mathbb{Z}_p}
\ar@/^1.2pc/[l]\ar@/_1.2pc/[l]}$$

\bigskip\noindent Thus, the recollement $(\ddag)$ above Corollary \ref{real-cotilt} constructed
from a cotilting module $U$ depends on injective cogenerator with
respect to which the $U$ is defined.

\subsection{Necessary conditions of homological subcategories from cotilting modules \label{6.2}}

We keep the notation in Section \ref{6.1}. For the cotilting module
$U$, we denote by
$$0\lra U_n \lraf{\partial_n} U_{n-1}\lra \cdots \lraf{\partial_2}
U_1 \lraf{\partial_1}U_0\lraf{\partial_0} W \lra 0
$$ the exact sequence in the axiom $(C_3)'$, and by
$\cpx{U}$ the following complex
$$\cdots \lra 0\lra U_n \lraf{\partial_n} U_{n-1}\lra \cdots \lraf{\partial_2}
U_1 \lraf{\partial_1} U_0 \lra 0\lra \cdots$$ with $U_i$ in degree
$-i$ for all $0\leq i\leq n$. Then $\partial_0$ induces a canonical
quasi-isomorphism $\widetilde{\partial}_0: \cpx{U}\ra W$ in $\C{A}$.
Recall that the complex $\cpx{I}$ in Corollary \ref{real-cotilt} (b)
also yields a canonical quasi-isomorphism $\xi: U\ra \cpx{I}$ in
$\C{A}$.

Furthermore, by the proof of the first part of Lemma \ref{Z-iso},
one can show that $\widetilde{\partial}_0$ and $\xi$ do induce the
following quasi-isomorphisms
$$(\ast)\quad
\xymatrix{\Hom_A(W, \cpx{I})\ar[r]^-{(\widetilde{\partial}_0)_*} &
\cpx{\Hom}_A(\cpx{U}, \cpx{I}) & \ar[l]_-{\xi^*} \Hom_A(\cpx{U},
U)}$$ in $\C{\mathbb{Z}}$. Here, we leave checking the details to
the reader.

Consequently, the morphism
$(\widetilde{\partial}_0)_*({\xi^*})^{-1}: \Hom_A(W, \cpx{I})\lra
\Hom_A(\cpx{U}, U)$ in $\D{\mathbb Z}$ is an isomorphism (compare
with Lemma \ref{Z-iso}). Due to the $A$-$\Lambda$-bimodule structure
of $W$, the former complex belongs to $\C{\Lambda}$. However, the
latter complex might not be a complex of $\Lambda$-modules since
$\cpx{U}$ is not necessarily a complex of $A$-$\Lambda$-bimodules in
general. This means that this isomorphism may not be extended to an
isomorphism in $\D{\Lambda}$. Nonetheless, for some special
cotilting modules, we do have this isomorphism in $\D{\Lambda}$. For
instance, in the case described in the following lemma.

\begin{Lem}\label{Z-iso2}
Suppose that $\Hom_A(U_i, U_{i+1})=0$ for $0\leq i < n$.

$(1)$ There exist a series of ring homomorphisms $\rho_j:
\Lambda\lra \End_A(U_j)$ for $0\leq j\leq n$, such that
$\,\widetilde{\partial}_0: \cpx{U}\lra W$ is a quasi-isomorphism in
$\C{A\otimes_\mathbb{Z}\Lambda\opp}$. In particular, the complexes
$\Hom_A(W, \cpx{I})$ and $\Hom_A(\cpx{U}, U)$ are isomorphic in
$\D{\Lambda}$.

$(2)$ If $\Ext^k_A(W, U_k)=\Ext_A^{k+1}(W, U_k)=0$ for all $0\leq k
< n$, then $\rho_n: \Lambda\lra \End_A(U_n)$ is an isomorphism.
\end{Lem}

{\it Proof.} $(1)$ Set $K_0:=W$, $K_n:=U_n$ and
$K_{m}:=\Ker(\partial_{m-1})$ for $1\leq m<n$. Then, for each $0\leq
i<n$, we have a short exact sequence $0\lra K_{i+1}\lra
U_i\lraf{\partial_i} K_i\lra 0$ of $A$-modules. In the following, we
shall define two ring homomorphisms $\varphi_i: \End_A(K_i)\lra
\End_A(U_i)$ and $\psi_i: \End_A(K_i)\lra \End_A(K_{i+1})$.

By Lemma \ref{cotilt1} (1), the sequence
$$0\lra \Hom_A(U, K_{i+1})\lra \Hom_A(U, U_i)\lraf{\partial_i^*}
\Hom_A(U, K_i)\lra 0$$ is exact. In particular, for $U_i\in\add(U)$,
the sequence
$$0\lra \Hom_A(U_i, K_{i+1})\lra \Hom_A(U_i, U_i)\lraf{\partial_i^*}
\Hom_A(U_i, K_i)\lra 0$$ is exact. Let $f\in\End_A(K_i)$. Then there
is a homomorphism $g\in\End_A(U_i)$ such that $\partial_i
f=g\,\partial_i$. We claim that such a $g$ is unique. Actually, if
there exists another $g'\in\End_A(U_i)$ such that $\partial_i
f=g'\,\partial_i$. Then $(g-g')\partial_i=0$, and so the map $g-g'$
factorizes through $K_{i+1}$. Note that each homomorphism $U_i\to
K_{i+1}$ also factorizes through $U_{i+1}$ via $\partial_{i+1}$.
This implies that $g-g':U_i\ra U_i$ factorizes through $U_{i+1}$.
However, since $\Hom_A(U_i, U_{i+1})=0$ by assumption, we have
$g=g'$. Hence, for a given $f$, such a $g$ is unique.

Now, we define $\varphi_i: f\mapsto g$ and $\psi_i: f\mapsto h$
where $h$ is the restriction of $g$ to $K_{i+1}$. This can be
illustrated by the following commutative diagram:
$$ \xymatrix{
0\ar[r]&\,K_{i+1}\ar@{-->}[d]^-{h}\ar[r]^-{\lambda_{i+1}}&U_i\ar@{-->}[d]^-{g}\ar[r]^-{\partial_i}
&K_i\ar[d]^-{f}\ar[r] &0\\
0\ar[r] &
K_{i+1}\ar[r]^-{\lambda_{i+1}}&U_i\ar[r]^-{\partial_i}&K_i\ar[r] &0}
$$
where $\lambda_{i+1}$ is the inclusion for $0\leq i\leq n-2$ and
$\lambda_{n}:=\partial_n$. Clearly, both $\varphi_i$ and $\psi_i$
are ring homomorphisms.

Recall that $\Lambda:=\End_A(W)=\End_A(K_0)$. Furthermore, for
$0\leq j\leq n$, we define $\rho_j: \Lambda\to \End_A(U_j)$ as
follows: If $j=0$, then $\rho_0:=\varphi_0\,$; if $j\geq 1$, then
$\rho_j$ is defined to be the composite of the following ring
homomorphisms:
$$\Lambda\lraf{\psi_0}\End_A(K_1) \lraf{\psi_1} \End_A(K_2)\lra \cdots
\lra \End_A(K_{j-1}) \lraf{\psi_{j-1}} \End_A(K_j)\lraf{\varphi_j}
\End_A(U_j)
$$
where $\varphi_n$ stands for the identity map. By definition, for
each $\lambda\in\Lambda$, there exists an exact commutative diagram
of $A$-modules:
$$ \xymatrix{0\ar[r] & U_n \ar[r]^-{\partial_n}\ar[d]_-{(\lambda)\rho_n} & U_{n-1}
\ar[r]\ar[d]^-{(\lambda)\rho_{n-1}}&\cdots \ar[r]^-{\partial_2} &
U_1 \ar[r]^-{\partial_1} \ar[d]^-{(\lambda)\rho_1} &
U_0\ar[r]^-{\partial_0} \ar[d]^-{(\lambda)\rho_0} & W \ar[r]\ar[d]^-{\lambda}  &0\\
0\ar[r] & U_n \ar[r]^-{\partial_n} & U_{n-1}\ar[r] &\cdots
\ar[r]^-{\partial_2} & U_1 \ar[r]^-{\partial_1} &
U_0\ar[r]^-{\partial_0} & W \ar[r] & 0}
$$
Note that $U_j$ is a natural $A$-$\End_A(U_j)$-bimodule and can be
regarded as an $A$-$\Lambda$-bimodule via $\rho_j$. It follows from
the above commutative diagram that $\partial_j$ is a homomorphism of
$A$-$\Lambda$-bimodules. This implies that $\widetilde{\partial}_0:
\cpx{U}\lra W$ can be viewed as a quasi-isomorphism in
$\C{A\otimes_\mathbb{Z}\Lambda^{op}}$. In this sense, the
quasi-isomorphisms in $(\ast)$ actually belong to $\C{\Lambda}$.
Thus $\Hom_A(W, \cpx{I})$ and $\Hom_A(\cpx{U}, U)$ are isomorphic in
$\D{\Lambda}$. This finishes $(1)$.

$(2)$ To show that $\rho_n$ is an isomorphism of rings, it suffices
to prove that $\psi_i$ is an isomorphism for $0\leq i\leq n-1$. Let
$i$ be such a fixed number. If $\Hom_A(K_i, U_i)=0$, then $\psi_i$
is injective. If the induced map $(\lambda_{i+1})_*:\Hom_A(U_i,
U_i)\lra \Hom_A(K_{i+1}, U_i)$ is surjective, then so is $\psi_i$.
Thus, by our assumptions in (2), to show that $\psi_i$ is an
isomorphism, it suffices to show that $\Hom_A(K_i, U_i)\simeq
\Ext^i_A(W, U_i)$ and that there exists an exact sequence of abelian
groups:
$$(\ast\ast)\quad
 \Hom_A(U_i,\,U_i)\lraf{(\lambda_{i+1})_*} \Hom_A(K_{i+1},
U_i)\lra \Ext_A^{i+1}(W, U_i)\lra 0.$$

In fact, since $U_s\in\add(_AU)$ for $0\leq s\leq n$, we have
$\Ext^r_A(U_s, X)=0$ for each $r\geq 1$ and $X\in\add(_AU)$ by the
axiom $(C2)$. Now, for $1\leq j\leq n$ and $X\in \add(_AU)$, one can
apply $\Hom_A(-, X)$ to the long exact sequence
$$0\lra K_{j}\lraf{\lambda_j} U_{j-1}\lra \cdots \lra U_1\lra U_0\lra W\lra 0,$$
and get an exact sequence $\Hom_A(U_{j-1}, X)\lraf{(\lambda_j)_*}
\Hom_A(K_j, X)\lra \Ext_A^j(W, X)\lra 0$ of abelian groups. If we
take $j:=i$ and $X:=U_i$, then $\Hom_A(K_i, U_i)\simeq\Ext_A^i(W,
U_i)$ since $\Hom_A(U_{i-1}, U_i)=0$ by assumption. If we take
$j:=i+1$ and $X:=U_i$, then we get the required sequence
$(\ast\ast)$. This finishes the proof of (2). $\square$

\medskip
The following result will be used for getting a counterexample which
demonstrates that, in general, the category $\Ker({\bf H})$ in
Corollary \ref{real-cotilt} may not be homological.

\begin{Koro}\label{new}
Keep all the assumptions in Corollary \ref{real-cotilt}. Further,
suppose that $n\geq 2$ and $U$ has injective dimension exactly equal
to $n$. If $\,\Hom_A(U_i, U_{i+1})=\Ext^i_A(W, U_i)$
$=\Ext_A^{i+1}(W, U_i)=0$ for all $0\leq i < n$, then the category
$\Ker({\bf H})$ is not a homological subcategory of $\D{R}$.
\end{Koro}

{\it Proof.} Suppose contrarily that $\Ker({\bf H})$ is homological
in $\D{R}$. Then, by Corollary \ref{real-cotilt}, we certainly have
$H^n\big({_R}\Hom_A(U,\,W)\otimes_{\Lambda}\Hom_A(W,\,\cpx{I})\big)=0$.
Furthermore, since $\Hom_A(U_i, U_{i+1})=0$ for all $0\leq i\leq
n-1$, we know from Lemma \ref{Z-iso2} (1) that $\Hom_A(W,
\cpx{I})\simeq \Hom_A(\cpx{U}, U)$ in $\D{\Lambda}$. Thus {\small
$$0=H^n\big(\Hom_A(U,\,W)\otimes_{\Lambda}\Hom_A(W,\,\cpx{I})\big)
\simeq\Hom_A(U,\,W)\otimes_{\Lambda}H^n\big(\Hom_A(W,\,\cpx{I})\big)\simeq
\Hom_A(U,\,W)\otimes_{\Lambda}H^n\big(\Hom_A(\cpx{U},\,U)\big).$$\vspace{-0.3cm}}

\noindent In particular, we have
$\Hom_A(U,\,W)\otimes_{\Lambda}H^n\big(\Hom_A(\cpx{U},\,U_n)\big)=0$,
due to $U_n\in\add(_AU)$. Recall that the complex $\Hom_A(\cpx{U},
U_n)$ is of the form {\small $$0\lra \Hom_A(U_0, U)
\lraf{(\partial_1)_*} \Hom_A(U_1, U_n)\lra \cdots\lra
\Hom_A(U_{n-2}, U_n) \lraf{(\partial_{n-1})_*} \Hom_A(U_{n-1}, U_n)
\lraf{(\partial_n)_*} \Hom_A(U_n, U_n) \lra 0$$\vspace{-0.5cm}}

\noindent with $\Hom_A(U_n, U_n)$ in degree $n$. Since
$\Hom_A(U_{n-1}, U_n)=0$, we obtain
$H^n\big(\Hom_A(\cpx{U},\,U_n)\big)=\End_A(U_n)$, and so
$\Hom_A(U,\,W)\otimes_{\Lambda}\End_A(U_n)=0$. Note that the left
$\Lambda$-module structure of $\End_A(U_n)$ is defined by the ring
homomorphism $\rho_n: \Lambda\lra \End_A(U_n)$ (see Lemma
\ref{Z-iso2} (1)). Since $\Ext^i_A(W, U_i)$ $=\Ext_A^{i+1}(W,
U_i)=0$ for all $0\leq i\leq n-1$, we see from Lemma \ref{Z-iso2}
(2) that $\rho_n$ is an isomorphism. This implies that
$$\Hom_A(U,\,W)\otimes_{\Lambda}\End_A(U_n)\simeq
\Hom_A(U,\,W)\otimes_{\Lambda}\Lambda\simeq\Hom_A(U,\,W)$$ and
therefore $\Hom_A(U,\,W)=0$. Since $_AW$ is an injective
cogenerator, we must have $U=0$. This is a contradiction. Thus
$\Ker({\bf H})$ is not homological in $\D{R}$. $\square$

\section{Counterexamples and open questions} \label{sect7}

In this section, we shall apply results in the previous sections to
give two examples which show that, in general, the category
$\Ker(_AT\otimesL_B-)$ for an $n$-tilting module $T$, or the
category $\Ker(\bf{H})$ for an $n$-cotilting module $U$ may not be
homological. At the end of this section, we mention a few open
questions related to some results in this paper.

Throughout this section, we assume that $A$ is a commutative,
noetherian, $n$-Gorensteion ring for a natural number $n$. Recall
that a ring is called $n$-\emph{Gorenstein} if the injective
dimensions of the regular left and right modules are at most $n$.

For an $A$-module $M$, we denote by $E(M)$ its injective envelope.
It is known that if $\mathfrak p$ and $\mathfrak q$ are two prime
ideals of $A$, then $\Hom_A(E(A/\mathfrak p), E(A/\mathfrak q))\neq
0$ if and only if $\mathfrak p \subseteq \mathfrak q$ (see
\cite[Theorem 3.3.8]{EJ}). In particular, $E(A/\mathfrak p)\simeq
E(A/\mathfrak q)$ if and only if $\mathfrak p=\mathfrak q$

\subsection{Higher $n$-tilting modules}\label{7.1}

In the following, we shall apply Corollary \ref{cor} to provide an
example of a good $n$-tilting $A$-module $T$ for which the category
$\Ker(_AT\otimesL_B-)$ in Theorem \ref{main-result} is not
homological.

For the $n$-Gorenstein ring $A$, it follows from  a classical result
of Bass that the regular module $_AA$ has a minimal injective
coresolution of the form:
$$
0\lra A\lra\bigoplus_{\mathfrak p\in {\mathcal P}_0}E(A/\mathfrak
p)\lra\cdots\lra \bigoplus_{\mathfrak p\in
\mathcal{P}_n}E(A/\mathfrak p)\lra 0, $$ where $\mathcal{P}_i$
stands for the set of all prime ideals of $A$ with height $i$ (see
\cite[Theorem 1, Theorem 6.2]{Bs}). It was pointed out in
\cite[Introduction]{TP} that the $A$-module
$$T:=\bigoplus_{0\leq i\leq
n}\,\bigoplus_{\mathfrak p\in \mathcal{P}_i}E(A/\mathfrak p)$$ is an
(infinitely generated) $n$-tilting module.

Clearly, the tilting module $_AT$ is good if we define
$T_i:=\bigoplus_{\mathfrak p\in \mathcal{P}_i}E(A/\mathfrak p)$.
Observe that, for $0\leq i<j\leq n$, we have $\Hom_A(E(A/\mathfrak
p), E(A/\mathfrak q))=0$ for $\mathfrak p\in P_j$ and $\mathfrak
q\in P_i$, and therefore $\Hom_A(T_j, T_i)=0$.

Now, we suppose that $n\geq 2$ and the injective dimension of $A$ is
exactly equal to $n$ (or equivalently, the Krull dimension of $A$ is
exactly $n$).

Note that $T_i\neq 0$ for all $2\leq i\leq n$ and that $T$ satisfies
the assumptions in Corollary \ref{cor} (2). Since the above
injective coresolution of $A$ is minimal, the module $_AT$ has
projective dimension equal to $n$ (see \cite[Proposition 3.5]{Bz1}).
By Corollary \ref{cor} (2), the category $\Ker(_AT\otimesL_B-)$ is
not homological in $\D B$. This means that for this tilting module
$T$, the subcategory $\Ker(_AT\otimesL_B-)$ cannot be realized as
the derived module category $\D C$ of a ring $C$ with a homological
ring epimorphism $B\ra C$. Thus, for higher $n$-tilting modules, the
answer to the question in Introduction is negative in general.

\subsection{Higher $n$-cotilting modules}\label{7.2}

Next, we apply Corollary \ref{new} to present an example of a good
$n$-cotilting $A$-module $U$, for which the category $\Ker({\bf H})$
in Corollary \ref{real-cotilt} is not homological in $\D{R}$.

Assume further that the ring $A$ is \emph{local} with the unique
maximal ideal $\mathfrak m$. In this case, $T_n$ is an injective
cogenerator for $A\Modcat$ since $\mathcal{P}_n$ is just the set
$\{\mathfrak m\}$. This follows from a general statement in
commutative algebra: If $S$ is a commutative noetherian ring, then
$\bigoplus_{\mathfrak m} E(S/\mathfrak{m})$ is an injective
cogenerator for $S$-Mod, where $\mathfrak m$ runs over all maximal
ideals of $S$.

Now, we take $$W:=T_n\quad \mbox{and}\quad  U:=\Hom_A(T,
W)=\bigoplus_{j=0}^n\Hom_A(T_j, W).$$ Since $_AT$ is an $n$-tilting
$A$-module, the module $_AU$ is an $n$-cotilting $A$-module.
Furthermore, applying $\Hom_A(-, W)$ to the minimal injective
coresolution of $_AA$, we get the following exact sequence of
$A$-modules:
$$
0\lra \Hom_A(T_n, W)\lra \Hom_A(T_{n-1}, W)\lra \cdots\lra
\Hom_A(T_1, W)\lra \Hom_A(T_0, W)\lra W\lra 0.
$$
This implies that the cotilting $A$-module $U$ is good if we define
$U_j:=\Hom_A(T_j, W)$ for $0\leq j\leq n$ (see the axiom $(C_3)'$ in
Definition \ref{cotilting}).

To see that $\Lambda:=\End_A(W)$ is a right noetherian ring, we note
that $W = E(A/\mathfrak m)$ and that $\Lambda$ is isomorphic to the
$\mathfrak m$-adic complete of $A$ (see \cite[Theorem 3.4.1
(6)]{EJ}). Since $A$ is noetherian, the ring $\Lambda$ is also
noetherian (see \cite[Corollary 2.5.16]{EJ}).

In the following, we shall prove that $_AU$ satisfies all the
assumptions in Corollary \ref{new}. In fact, it suffices to show
that, for any $m\geq 0$, we have

$(a)$ $\Ext_A^m(U_r, U_s)=0$ for $0\leq r<s\leq n$.

$(b)$ $\Ext_A^m(W, U_i)=0$ for $0\leq i\leq n-1$, and $\Ext_A^n(W,
U_n)\neq 0$.

\noindent The reason is the following: According to $(b)$, the
injective dimension of $U_n$ is at least $n$, and therefore exactly
$n$. This means that $_AU$ is a cotilting module of injective
dimension $n$. Moreover, from $(a)$ and $(b)$ we can conclude that
the assumptions in Corollary \ref{new} hold true for $U$. It then
follows from Corollary \ref{new} that, for this cotilting module
$U$, the category $\Ker({\bf H})$ in Corollary \ref{real-cotilt} is
not homological in $\D{R}$ with $R:=\End_A(U)$. In other words,
$\Ker({\bf H})$ cannot be realized as the derived module category
$\D S$ of a ring $S$ with a homological ring epimorphism $R\ra S$.

So, let us verify the above $(a)$ and $(b)$. First, we need the
following results about $n$-Gorenstein rings:

$(1)$ The flat dimension of the $A$-module $T_j$ is exactly $j$.

$(2)$ Any flat $A$-module $F$ admits a minimal injective
coresolution of the form
$$
0\lra {}_AF \lra I_0\lra I_1\lra \cdots\lra I_{n-1}\lra I_{n}\lra 0
$$
such that $I_j\in\Add(T_j)$ for all $0\leq j\leq n$.

$(3)$ Let $\mathfrak{p}$ and $\mathfrak{q}$ be prime ideals of $A$.
If $\mathfrak{p}\nsubseteq \mathfrak{q}$ or $\mathfrak{q}\nsubseteq
\mathfrak{p}$, then $\Tor_m^A(E(A/\mathfrak {p}),\, E(A/\mathfrak
{q}))=0$ for  all $m\geq 0$. Moreover, $\Tor_m^A(E(A/\mathfrak
{p}),\, E(A/\mathfrak {p}))\neq 0$ if and only if $m$ equals the
height of $\mathfrak{p}$ in $A$.

Here, $(1)$ and $(2)$ follow from \cite[Proposition 2.1 and Theorem
2.1]{X}, while $(3)$ is taken from \cite[Lemma 9.4.5 and Theorem
9.4.6]{EJ}.

\smallskip
Since the dual $A$-module $\Hom_A(F, W)$ of a flat $A$-module $F$ is
injective, we know from $(1)$ that the injective dimension of
$_AU_j$ is at most $j$. Since the dual $A$-module $\Hom_A(I, W)$ of
an injective $A$-module $I$ is always flat (see \cite[Corollary
3.2.16 (2)]{EJ}), we see that the $A$-module $U_j$ is flat since
$T_j$ is injective. It then follows from $(2)$ that $U_j$ admits a
minimal injective coresolution of the form
$$
0\lra U_j \lra I_{j, 0}\lra I_{j, 1}\lra \cdots\lra I_{j, j-1}\lra
I_{j, j}\lra 0
$$
with $I_{j,\, k}\in\Add(T_k)$ for all $0\leq k\leq j$.

Now, we show $(a)$. Actually, by Lemma \ref{tor-ext} (1), we have
$$\Ext_A^m(U_r, U_s)=\Ext_A^m\big(U_r, \,\Hom_A(T_s,
W)\big)\simeq\Hom_A\big(\Tor_m^A(T_s, U_r),\, W\big) \mbox{\; for
\;} m\ge 0.$$ Note that the flatness of $U_r$ implies that
$\Ext_A^m(U_r, U_s)=0$ for $m\ge 1$. It remains to show
$\Hom_A(U_r,\, U_s)=0$. For this aim, it is sufficient to show
$T_s\otimes_AU_r=0$. Since $T_s:=\bigoplus_{\mathfrak p\in
\mathcal{P}_s}E(A/\mathfrak p)$ and the functor $-\otimes_AU_r$
commutes with arbitrary direct sums, we have to prove $E(A/\mathfrak
p)\otimes_AU_r=0$ for every $\mathfrak p\in \mathcal{P}_s$. In fact,
since $r<s$ by assumption, we know that $\mathfrak{p}\nsubseteq
\mathfrak{q}$ for each $\mathfrak{q}\in\mathcal{P}_k$ with $0\leq
k\leq r$. It follows from $(3)$ that $\Tor_j^A\big(E(A/\mathfrak
p),\,E(A/\mathfrak q)\big)=0$ for all $j\geq 0$, and therefore
$$\Tor_j^A(E(A/\mathfrak p),\,T_k)\simeq \bigoplus_{\mathfrak q\in
\mathcal{P}_k}\Tor_j^A\big(E(A/\mathfrak p),\,E(A/\mathfrak
q)\big)=0.$$ Since $I_{r,\,k}\in \Add(T_k)$, we obtain
$\Tor_j^A(E(A/\mathfrak p),\,I_{r,\,k})=0$ for all $j\geq 0$. Now,
by applying the tensor functor $E(A/\mathfrak p)\otimes_A-$ to the
minimal injective coresolution of $U_r$, we can prove $E(A/\mathfrak
p)\otimes_AU_r=0$. Thus $T_s\otimes_AU_r=0$. This finishes the proof
of $(a)$.

Finally, we show $(b)$. Let $0\leq i\leq n-1$. Recall that $U_i =
\Hom_A(T_i,\,W)$. According to Lemma \ref{tor-ext} (1), we have
{\small $$\Ext_A^m\big(W,
\Hom_A(T_i,\,W)\big)\simeq\Hom_A\big(\Tor_m^A(T_i,\,
W),\,W\big)\simeq \Hom_A\big(\bigoplus_{\mathfrak p\in
\mathcal{P}_i}\Tor_m^A(E(A/\mathfrak p),\, W),\,W\big)\simeq
\prod_{\mathfrak p\in
\mathcal{P}_i}\Hom_A\big(\Tor_m^A(E(A/\mathfrak p),\,
W),\,W\big).$$\vspace{-0.2cm}}

\noindent Since the ideal $\mathfrak m$ is maximal (or of height
$n$), it holds that $\mathfrak {m}\nsubseteq \mathfrak{p}$ for every
$\mathfrak{p}\in \mathcal{P}_i$. Hence it follows from $(3)$ that
$\Tor_m^A(E(A/\mathfrak p),\, W)=0$, and therefore $\Ext_A^m(W,
U_i)=0$. Similarly, one can show that
$$\Ext_A^n(W, U_n)=\Ext_A^n(W,\,\Hom_A(W, W))\simeq
\Hom_A(\Tor_n^A(W,\, W),\,W).$$ Since $\Tor^A_n(W,
W)=\Tor^A_n\big(E(A/\mathfrak m),\, E(A/\mathfrak m)\big)\neq 0$ by
$(3)$ and since $_AW$ is an injective cogenerator, we infer that
$\Ext_A^n(W, U_n)\neq 0$. Thus $(b)$ follows.

Consequently, for the $n$-cotilting $A$-module $U$, the subcategory
$\Ker(\bf{H})$ is not homological in $\D R$.

\medskip
Let us end this paper by the following open questions related to our
results in this note.

\medskip
{\bf Question 1.} Let $A$ be a ring with identity. Is there a good
$n$-tilting $A$-module $T$ for $n\ge 2$ such that $T$ is not
equivalent to any classical tilting $A$-module and that
$\Ker(T\otimesL_B-)$ is homological?

\medskip
{\bf Question 2.} Is the converse of Corollary \ref{cor} (1) always
 true?

\medskip
For tilting modules over commutative noetherian $n$-Gorenstein
rings, Silvana Bazzoni even guesses a stronger answer: If
$\Ker(T\otimesL_{B}-)$ is homological in $\D B$, then $_AT$ should
be a $1$-tilting module, that is, the module $_AN$ in Corollary
\ref{cor} (1) should be zero.

\medskip
{\bf Question 3.} Given a good $1$-cotilting module $U$ over an
\textbf{arbitrary ring} $A$, is there a homological ring epimorphism
$\lambda: \End_A(U)\ra C$ and a recollement of the following form?
$$\xymatrix@C=1.2cm{\D{C}\ar[r]^-{D(\lambda_*)}
&\D{\End_A(U)}\ar[r]\ar@/^1.2pc/[l]\ar@/_1.2pc/[l] &\D{A}
\ar@/^1.2pc/[l]\ar@/_1.2pc/[l]}\vspace{0.2cm}$$

Note that this reccollement does not involve the derived categories
of the endomorphism rings of any injective cogenerators related to
$U$.

\medskip
{\bf Question 4.} Given an arbitrary ring $A$, how to parameterize
homological subcategories of $\D A$? Equivalently, how to classify
homological ring epimorphisms starting from $A$?

\medskip
{\bf Question 5.} Is the Ringel $R$-module $M$ in Lemma
\ref{cotilt1} always good?

A positive answer to this question would lead to a generalization of
Corollary \ref{real-cotilt}.

\medskip
{\bf Acknowledgement.} The research work of the corresponding author
C.C. Xi is partially supported by a grant of MOE (20110003110003),
China.

{\footnotesize
\bigskip Hongxing Chen, Beijing International Center for Mathematical
Research, Peking University, 100871 Beijing, People's Republic of
China

{\tt Email: chx19830818@163.com}

\bigskip
Changchang Xi, School of Mathematical Sciences, Capital Normal
University, 100048 Beijing, People's Republic of  China

{\tt Email: xicc@cnu.edu.cn} $\qquad$


\medskip
First version: April 6, 2012
\end{document}